\documentclass[twoside, 11pt]{article}

\usepackage{amsmath}
\usepackage{paralist}
\usepackage{amssymb,mathrsfs,psfrag,graphicx}
\usepackage{hyperref}
\def\ve{\varepsilon}

\def\x{\xi}
\def\t{\theta}

\def\k{\kappa}

\def\mb{\mathbf}
\def\a{\alpha}
\def\b{\beta}

\def\d{\delta}
\def\l{\lambda}
\def\f{\frac}
\def\p{\phi}
\def\r{\rho}

\def\s{\sigma}
\def\z{\zeta}
\def\o{\omega}
\def\di{\displaystyle}
\def\i{\infty}
\def\wt{\widetilde}

\newtheorem{theorem}{Theorem}[section]
\newtheorem{lemma}[theorem]{Lemma}

\newtheorem{proposition}{Proposition}[section]
\newtheorem{remark}{Remark}

\renewcommand{\theequation}{\thesection.\arabic{equation}}
\renewcommand{\thetheorem}{\thesection.\arabic{theorem}}
\renewcommand{\thelemma}{\thesection.\arabic{lemma}}

\begin{document}
	
	\title{\bf  Time-asymptotic stability of generic Riemann solutions for Boltzmann equation} \vskip 0.5cm
	\author{
		Yi Wang,\thanks{Institute of Applied Mathematics, Academy of Mathematics and Systems Science, Chinese Academy of Sciences, Beijing 100190, P. R. China
			and School of Mathematical Sciences, University of Chinese Academy of Sciences, Beijing 100049, P. R. China ({\tt wangyi@amss.ac.cn}). The work of Yi Wang is partially supported by NSFC grants (Grant No.s 12171459, 12288201, 12090014, 12421001) and CAS Project for Young Scientists in Basic Research, Grant No. YSBR-031.} \quad Qiuyang Yu \thanks{Institute of Applied Mathematics, Academy of Mathematics and Systems Science, Chinese Academy of Sciences, Beijing 100190, P. R. China
			and School of Mathematical Sciences, University of Chinese Academy of Sciences, Beijing 100049, P. R. China ({\tt yuqiuyang@amss.ac.cn}).}
	}
	\date{}
	\maketitle
	\begin{abstract}
		Time-asymptotic stability of generic Riemann solution, consisting of a rarefaction wave, a contact discontinuity and a shock, for the one-dimensional Boltzmann equation, has been a long-standing open problem in kinetic theory. In this paper, we proved that the composite waves of  generic Riemann profile including the inviscid self-similar rarefaction wave, the viscous contact wave (i.e., the viscous version of contact discontinuity) and the viscous shock profile with the time-dependent shift to both macroscopic and microscopic components are nonlinearly stable for the one-dimensional Boltzmann equation, by the first using the $a$-contraction method to the Boltzmann equation. Compared with the compressible Navier-Stokes-Fourier equations, the new difficulties here lie in the microscopic effects of the Boltzmann shock profile and their interactions and/or couplings with the rarefaction wave, viscous contact wave and the macroscopic components from the macro-micro decomposition of the Boltzmann equation.

		
	\end{abstract}
	\tableofcontents
	\section{Introduction}
	\renewcommand{\theequation}{\arabic{section}.\arabic{equation}}
	\setcounter{equation}{0}

	The one-dimensional Boltzmann equation takes the form
	\begin{align}\label{(B)}
		f_t + \xi_1 f_x= Q(f,f),
	\end{align}
	where  $\xi=(\xi_1,\xi_2,\xi_3)\in {\mathbb R}^3$, $x\in{\mathbb R}^1$, $t\in{\mathbb R}_+$ and $f(t,x,\xi)$ is the density distribution function of particles at the time $t$ with the location
	$x$  and the velocity $\xi$. Equation (\ref{(B)}) was first established by Boltzmann \cite{Boltzmann} in 1872 to describe the motion of rarefied gases and it is a fundamental equation in statistics physics. For the hard sphere model, the collision operator $Q(f,f)$ is bilinear as
	\begin{equation}\label{B}
		\begin{aligned}
			\di Q(g,h)(\xi)& =\di 
			\int_{{\mathbb R}^3}\!\int_{{\mathbb S}^2_+} |(\xi-\xi_*)\cdot \Omega|\Big(
			g(\xi')h(\xi_*')- g(\xi)h(\xi_*) \Big) d\Omega d \xi_*\\[1mm]
			&:=\di Q_+(g,h)(\xi)-Q_-(g,h)(\xi),
		\end{aligned}
	\end{equation}
	where ${\mathbb S}^2_+:=\left\lbrace \Omega\in {\mathbb S}^2|(\xi-\xi_*)\cdot \Omega\geqslant0\right\rbrace $ with ${\mathbb S}^2$ being a two-dimensional unit sphere and $\xi',\xi_*'$ are the velocities after an elastic collision of two particles with velocities $\xi$ and $\xi_*$ before the collision. Note that the collision operator $Q(g,h)(\xi)$ can be split into the gain and loss terms, namely,
	\begin{equation}\label{gain}
		Q_+(g,h)(\xi):=\int_{{\mathbb R}^3}\!\int_{{\mathbb S}^2_+} |(\xi-\xi_*)\cdot \Omega|
		g(\xi')h(\xi_*') d\Omega d \xi_*
	\end{equation}
	being the gaining of the particle with velocity $\xi$ from the collision of the two particles with velocities $\xi'$ and $\xi_*'$ and 
	\begin{equation}\label{loss}
		Q_-(g,h)(\xi):=g(\xi)\int_{{\mathbb R}^3}\!\int_{{\mathbb S}^2_+} |(\xi-\xi_*)\cdot \Omega|
		h(\xi_*) d\Omega d \xi_*
	\end{equation}
	being the loss of the particle with velocity $\xi$ due to the collision of  the two particles with velocities $\xi$ and $\xi_*$. The conservations of the momentum and the energy of the unit particles  yield the following relations between the velocities before and after the elastic collision:
	$$
	\xi'= \xi -[(\xi-\xi_*)\cdot \Omega] \; \Omega, \qquad
	\xi_*'= \xi_* + [(\xi-\xi_*)\cdot \Omega] \; \Omega.
	$$
	
	The Boltzmann equation \eqref{(B)} is closely associated with fluid dynamical systems, such as the compressible Euler equations and the Navier-Stokes-Fourier equations. If a gas is in thermal equilibrium (i.e., the density distribution function is a local Maxwell distribution), then the compressible Euler equations can be derived from the Boltzmann equation through the celebrated Hilbert expansion \cite{Hilbert}. The compressible Euler equations, serving as a typical example for hyperbolic conservation laws, have been garnered considerable attention.  The main features of the compressible Euler equations are the formation of the shock singularities, no matter how smooth or small the initial values are.  If the piece-wise constant discontinuous Riemann data is given, then the corresponding entropic Riemann solutions contain two nonlinear waves, i.e., shock and rarefaction waves in the genuinely nonlinear characteristic fields, and a linearly degenerate wave, called contact discontinuity. Generic Riemann solutions consisting those three elementary wave patterns are quite fundamental for the solutions to the general initial value problems of Euler equations. The local, global and large-time behaviors of the general small BV (i.e., bounded variation) solution are fully determined by Riemann solutions for compressible Euler equations. As for the compressible Navier-Stokes-Fourier equations, they can be derived from the Chapman-Enskog expansion \cite{Chapman-Cowling} of the Boltzmann equation by considering both viscosities and heat-conductivity.
	
	It can be expected that wave phenomena analogous to macroscopic fluid dynamics also exist in Boltzmann equation. The existence of Boltzmann shock profile was first proved by Nicolaenko and Thurber for the hard sphere model \cite{Nicolaenko-Thurber}  in 1975 and generalized by Caflisch and Nicolaenko \cite{Caflisch-Nicolaenko} to hard potential case in 1982. The positivity and nonlinear stability of the Boltzmann shock profile with zero mass condition were proven by Liu and Yu \cite{Liu-Yu-2004}  by using the macro-micro decomposition introduced in \cite{Liu-Yu-2004} and further elaborated by Liu, Yang and Yu \cite{Liu-Yang-Yu} and then the zero mass condition is removed by Yu \cite{Yu-2010} using the point-wise Green function around Boltzmann shock profile. The compressibility of the Boltzmann shock profile, crucially needed in the stability analysis, can be found in  Liu and Yu \cite{Liu-Yu-2013} and  Pogan and Zumbrun \cite{Pogan-Zumbrun-2018} by different methods.  Wang and Wang \cite{Wang-Wang-2015} proved the nonlinear stability of the composition of two Boltzmann shock profiles without zero mass condition by using the weighted energy methods. On the other hand, Liu, Yang, Yu and Zhao \cite{Liu-Yang-Yu-Zhao} proved that the inviscid self-similar rarefaction wave is time-asymptotically stable for Boltzmann equation, and Huang and Yang \cite{Huang-Yang} and
	Huang, Xin and Yang \cite{Huang-Xin-Yang} proved the meta-stability of viscous contact wave, which can be viewed as the viscous version of the inviscid contact discontinuity, to the nonlinear Boltzmann equation with and without zero mass condition respectively. Therefore, the time-asymptotic stability of single wave pattern to the Boltzmann equation have been well-established since the invention of macro-micro decomposition in \cite{Liu-Yu-2004} and \cite{Liu-Yang-Yu}. However, the time-asymptotic stability of composite waves of different types of wave patterns, in particular the case of generic Riemann solution consisting of all three types of elementary waves, i.e., shock, rarefaction wave and contact discontinuity, for the one-dimensional Boltzmann equation, has been a long-standing open problem in kinetic theory. In the present paper, we aim to resolve this problem and prove that the composite waves of  generic Riemann profile including the inviscid self-similar rarefaction wave, the viscous contact wave (i.e., the viscous version of contact discontinuity) and the viscous shock profile with the time-dependent shift to both macroscopic and microscopic components are nonlinearly stable by the first using the $a$-contraction method to the Boltzmann equation.

	By Chapman-Enskog expansion \cite{Chapman-Cowling} and the macro-micro decomposition in \cite{Liu-Yu-2004} and \cite{Liu-Yang-Yu}, Boltzmann equation can be decomposed into the macroscopic part which satisfies the compressible Navier-Stokes-Fourier equations coupled with the microscopic equation. Therefore, the time-asymptotic stability of wave patterns to Boltzmann equation is inspired by but far from enough the compressible Navier-Stokes-Fourier equations and the viscous conservation laws. In the past decades, there are plenty of literatures and much progress on the time-asymptotic stability of basic wave patterns to the viscous conservation laws since the pioneer work of Il'in and Oleinik \cite{Il'in-Oleinik} for Burgers equation in 1960. Then the stability of single viscous shock, rarefaction wave and viscous contact wave are proven and understood rather satisfactorily by direct or weighted energy methods, spectral methods, point-wise Green function methods, $L^1$-stability and even the combined methods mentioned above. Shock wave is a compressed wave such that the classical  $L^2$-relative entropy methods can not be directly utilized and then the anti-derivative variables for the perturbation around the viscous shock profile is introduced by Matsumura and Nishihara \cite{Matsumura-Nishihara-1985} and Goodman \cite{Goodman} independently in 1985-1986 to fully use the compressibility of viscous shock. On the other hand, rarefaction wave is expanding and the direct $L^2$-relative entropy methods around the rarefaction wave can be successfully applied to obtain its stability by Matsumura and Nishihara \cite{Matsumura-Nishihara-1986}. While the viscous contact wave is proven nonlinearly meta-stable by using anti-derivative techniques or direct $L^2$-relative entropy methods with suitably weighted estimates by Huang, Matsumura and Xin \cite{Huang-Matsumura-Xin} and Huang, Xin and Yang \cite{Huang-Xin-Yang}.  It should be emphasized that the stability proof-frameworks for these three individual wave patterns are quite different and sharply incompatible with each other due to the distinct intrinsic properties of the three waves. 
	Therefore, it is highly nontrivial to prove the stability for the composite waves of different types of wave patterns. 
	
	In 2010, Huang, Li and Matsumura \cite{Huang-Li-Matsumura} first proved the stability of the composite wave of two rarefaction waves and a viscous contact wave to one-dimensional compressible Navier-Stokes-Fourier equations by establishing a new heat-kernel inequality. Very recently, by using $a$-contraction method invented by Kang and Vasseur \cite{Kang-Vasseur-2017} with the time-dependent shift to the viscous shock wave,  Kang, Vasseur and Wang \cite{Kang-Vasseur-Wang-2023} successfully proved the time-asymptotic stability of the composite wave of viscous shock and rarefaction wave to barotropic Navier-Stokes equations and then the generic Riemann profiles containing rarefaction wave, viscous contact wave and viscous shock to full compressible Navier-Stokes-Fourier equations in \cite{Kang-Vasseur-Wang-2024}. 
	
	For the time-asymptotic stability of the composite waves of different types of wave patterns to Boltzmann equation, in particular the generic Riemann solution case,  besides all the difficulties encountered for Navier-Stokes-Fourier equations as in \cite{Kang-Vasseur-Wang-2024}, the new difficulties lie in the microscopic effects of the Boltzmann shock profile and their interactions and/or couplings with the rarefaction wave, viscous contact wave and the macroscopic components from the macro-micro decomposition of the Boltzmann equation. For the Boltzmann shock profile, even it can be well approximated  by Navier-Stokes-Fourier shock, the microscopic effect is essential for its time-asymptotic stability and persists for all time, which is quite different from the other two wave patterns, that is, rarefaction wave and contact discontinuity. For the stability of either rarefaction wave or viscous contact wave to Boltzmann equation,  the microscopic effect disappears time-asymptotically, even though it affects the corresponding solution behaviors in any finite time. Motivated by \cite{Kang-Vasseur-Wang-2024} for Navier-Stokes-Fourier equations, $a$-contraction method for the time-asymptotic stability of Boltzmann shock profile is needed for the consistence of its stability frameworks towards rarefaction wave and viscous contact wave. Therefore, time-dependent shift should be equipped to both macroscopic and microscopic components of Boltzmann shock profile, which yields the new solution behaviors beyond the Navier-Stokes-Fourier equations and brings main difficulties coming from their interactions and/or couplings with the rarefaction wave, viscous contact wave and the macroscopic components of the Boltzmann equation.

	Now we review the $a$-contraction method for the $L^2$-stability of the conservation laws. For the hyperbolic inviscid conservation laws, $L^1$-stability is extensively used and successfully applied to prove the global existence and uniqueness of the solution \cite{Kruzkov,Liu-Yang-1999,Bressan-Liu-Yang}, while $L^2$-relative entropy norm is natural from the viewpoint of the physical energy.  However, it can be shown that $L^2$-relative entropy around the inviscid shock is unstable even for the inviscid Burgers equation. With suitable time-dependent shift and weight function $a$, the inviscid extreme shock can be proven to be nonlinearly stable under the $L^2$-relative entropy perturbation \cite{Leger, Leger-Vasseur}. For the viscous conservation laws, the time-dependent shift is first applied to obtain $L^2$-stability of viscous shock profile in \cite{Kang-Vasseur-2017} without using the classical anti-derivative techniques and then is extended to the barotropic compressible Navier-Stokes equations for $L^2$-contraction and stability of weak viscous shock with both time-dependent shift and suitable weight function $a$ \cite{Kang-Vasseur-2021}. Since the $a$-contraction method for $L^2$-stability of viscous shock profile is energy based and is compatible with the stability proof framework of rarefaction wave and viscous contact wave, Kang, Vasseur and Wang \cite{Kang-Vasseur-Wang-2023, Kang-Vasseur-Wang-2024} proved the time-asymptotic stability of  generic Riemann profile for both barotropic Navier-Stokes equations and full Navier-Stokes-Fourier equations.

	In this paper, the time-asymptotic stability of generic Riemann profile, including rarefaction wave, viscous contact wave and Boltzmann shock profile, is investigated for Boltzmann equation. 
	Through the micro-macro decomposition, Boltzmann equation can be rewritten as the macroscopic Navier-Stokes-Fourier type equations coupled with the microscopic equation.  Thanks to \cite{Liu-Yu-2013} and \cite{Pogan-Zumbrun-2018}, the Boltzmann shock profile can be parameterized by a new variable $\eta$ which satisfies the Burgers-like equation \eqref{eta}, as shown in Appendix. Then by using $a$-contraction method, which is energy based and can seamlessly handle the superposition of waves of different kinds, the time-asymptotic stability of generic Riemann profile for Boltzmann equation is proven. As mentioned before, the macroscopic part can be treated similarly as in \cite{Kang-Vasseur-Wang-2024} for Navier-Stokes-Fourier equations. The main task here is to handle the microscopic part and its interactions and/or couplings with the macroscopic part. Due to the micro H-theorem of Boltzmann equation,  the dissipative properties of the linearized collision operator $\mb L$ around the equilibrium Maxwellian plays an important role for the analysis of the microscopic part. However, it is far from enough. Since the microscopic effect is essential for the time-asymptotic stability of Boltzmann shock profile, the time-dependent shift $\mb X(t)$ should also be equipped to the microscopic part of Boltzmann shock profile. When conducting the energy analysis in the microscopic level (see \eqref{M.10} and \eqref{Gt}) and highest order estimates (see \eqref{M.45}), we need to control both $\di \int_0^T |\dot{\mb X}(t)|^2dt$ and, in particular, $\di \int_0^T |\ddot{\mb X}(t)|^2dt$, which is quite different from the macroscopic Navier-Stokes-Fourier equations. Eventually, our time-asymptotic stability results are obtained based on standard local existence and uniform-in-time a priori estimates through the continuity argument.
	
	The Boltzmann equation has also been extensively studied in other important aspects, such as the renormalized solution, regularity of solutions, fluid dynamic limits, and global existence around a global Maxwellian, etc.; see \cite{Caflisch, DiPerna-Lions, Golse, Golse-Saint, Guo, Huang-Wang-Wang-Yang, Huang-Wang-Yang-2010-1, Huang-Wang-Yang-2010-2,  Lions-Masmoudi, Nishida, Xin-Zeng, Yu-2005} and the references therein.

	For a solution $f(t,x, \xi)$ of (\ref{(B)}), there are five conserved macroscopic quantities: the density $\rho(t,x)$, the momentum $m(t,x)=\rho u(t,x)$, and the total energy $E(t,x)=\rho(e+\frac{1}{2}|u| ^2) (t,x)$, given by
	\begin{align*}
		\left(\rho,\rho u_i,\rho\left( e+\frac{1}{2}\left|u \right| ^2\right)  \right)(t,x) =\int_{\mathbb R^3}(\varphi_0,\varphi_i,\varphi_4)(\xi)f(t,x,\xi)d\xi, \quad i=1,2,3,
	\end{align*}
	where $\varphi_i(\xi)(i=0,1,2,3,4)$ are the collision invariants given by
	\begin{equation}\label{collision-invar}
		\varphi_0(\xi) = 1,~~~
		\varphi_i(\xi) = \xi_i~  (i=1,2,3),~~~
		\varphi_4(\xi) = \f{1}{2} |\xi|^2,
	\end{equation}
	that satisfy
	$$
	\int_{{\mathbb R}^3} \varphi_i(\xi) Q(g,g)(\xi) d \xi =0,\quad {\textrm
		{for} } \ \  i=0,1,2,3,4.
	$$
	The local Maxwellian $\mb{M}$ associated to solution $f( t,x, \xi)$ to the Boltzmann equation
	(\ref{(B)}) is defined in terms of the conserved fluid variables:
	\begin{equation}
		\mb{M}=\mb{M}_{[\rho,u,\t]} (t,x,\xi) = \f{\rho(t,x)}{\sqrt{ (2 \pi
				R \t(t,x))^3}} e^{-\f{|\xi-u(t,x)|^2}{2R\t(t,x)}}. \label{localmaxwellian}
	\end{equation}
	Here $\t$ is the temperature which is related to the internal energy $e(t,x)=\frac{3}{2}R\t(t,x)$ with $R$ being a positive gas constant, and $u(t,x)=(u_1,u_2,u_3)(t,x)$ is the fluid velocity. 
	
	It is well known that when the gas is in local thermo-equilibrium, i.e., $f=\mb{M}$, the Boltzmann equation (\ref{(B)}) is reduced to the compressible Euler equations that consist of conservation of mass, momentum, and energy:
	\begin{equation}\label{Euler}
		\left\{
		\begin{array}{l}
			\di\rho_t+(\rho u_1)_x=0,\\[1mm]
			\di(\rho u_1)_t+(\rho u_1^2+p)_x=0,\\[1mm]
			\di(\rho u_i)_t+(\rho u_1u_i)_x=0,~i=2,3,\\[2mm]
			\di\left[ \rho\left( e+\f{|u|^2}{2}\right) \right] _t+\left[ \rho u_1\left( e+\f{|u|^2}{2}\right) +pu_1\right] _x=0,
		\end{array}
		\right.
	\end{equation}
	where $p=R\rho \t$ is the pressure. From now on, the inner product of $g_1$ and $g_2$ in
	$L^2_{\xi}({\mathbb R}^3)$ with respect to a given Maxwellian $\tilde{\mb{M}}$ is denoted by:
	\begin{equation}\label{product}
		\langle g_1,g_2\rangle_{\tilde{\mb{M}}}\equiv \int_{{\mathbb R}^3}
		\f{1}{\tilde{\mb{M}}}g_1(\xi)g_2(\xi)d \xi.
	\end{equation}
	If $\tilde{\mb{M}}$ is the local Maxwellian $\mb{M}$ defined in (\ref{localmaxwellian}), the macroscopic space is spanned by the following five pairwise orthogonal base,
	\begin{equation}\label{orthogonal-base}
		\left\{
		\begin{array}{l}
			\chi_0(\xi) \equiv {\di\f1{\sqrt{\rho}}\mb{M}}, \\[3mm]
			\chi_i(\xi) \equiv {\di\f{\xi_i-u_i}{\sqrt{R\t\rho}}\mb{M}} \ \ {\textrm {for} }\ \  i=1,2,3, \\[3mm]
			\chi_4(\xi) \equiv
			{\di\f{1}{\sqrt{6\rho}}\left( \f{|\xi-u|^2}{R\t}-3\right) \mb{M}},\\[4mm]
			\langle\chi_i,\chi_j\rangle=\delta_{ij}, ~i,j=0,1,2,3,4.
		\end{array}
		\right.
	\end{equation}
	For brevity, if $\tilde{\mb{M}}$ is the local Maxwellian $\mb{M}$, we will simply use $\langle\cdot,\cdot\rangle$ to denote $\langle\cdot,\cdot\rangle_{\mb{M}}$.

	By using the above base, the macroscopic projection $\mb{P}_0$ and microscopic projection
	$\mb{P}_1$ can be defined as
	\begin{equation*}
		\mb{P}_0g = {\di\sum_{j=0}^4\langle g,\chi_j\rangle\chi_j},\qquad
		\mb{P}_1g= g-\mb{P}_0g.
	\end{equation*}
	A function $g(\xi)$ is called microscopic, or non-fluid, if
	$$
	\int_{\mathbb R^3} g(\xi)\varphi_i(\xi)d\xi=0,~i=0,1,2,3,4,
	$$
	where again $\varphi_i(\xi)(i=0,1,2,3,4)$ represent the collision invariants.
	
	Under the above projection, the solution of the Boltzmann equation $f(t,x,\xi)$ can be decomposed into the macroscopic (fluid) component, i.e., the local Maxwellian $\mb{M}(t,x,\xi)$ defined in (\ref{localmaxwellian}), and the microscopic (non-fluid) component $\mb{G}(t,x,\xi)$,
	\begin{equation}\label{macro-micro}
		f(t,x,\xi)=\mb{M}(t,x,\xi)+\mb{G}(t,x,\xi),~~~
		\mb{P}_0f=\mb{M},~~~
		\mb{P}_1f=\mb{G},
	\end{equation}
	and the Boltzmann equation (\ref{(B)}) becomes
	\begin{equation}\label{decompo}
		(\mb{M}+\mb{G})_t+\xi_1(\mb{M}+\mb{G})_x
		=Q(\mb{M},\mb{G})+Q(\mb{G},\mb{M})+Q(\mb{G},\mb{G}).
	\end{equation}
	By integrating the product of the equation (\ref{decompo}) and the collision invariants $\varphi_i(\xi)(i=0,1,2,3,4)$  with respect to $\xi$ over ${\mathbb R}^3$, one has the following system for the fluid variables $(\rho, u, \theta)$:
	\begin{equation}\label{EulerwithG}
		\left\{
		\begin{array}{lll}
			\di \rho_{t}+(\rho u_1)_x=0, \\[1mm]
			\di (\rho u_1)_t+(\rho u_1^2
			+p)_x=-\int_{\mathbb R^3}\xi_1^2\mb{G}_xd\xi,  \\[3mm]
			\di (\rho u_i)_t+(\rho u_1u_i)_x=-\int_{\mathbb R^3}\xi_1\xi_i\mb{G}_xd\xi,~ i=2,3,\\[3mm]
			\di \left[ \rho\left( e+\f{|u|^2}{2}\right) \right] _t+\left[ \rho
			u_1\left( e+\f{|u|^2}{2}\right) +pu_1\right] _x=-\int_{\mathbb R^3}\f12\xi_1|\xi|^2\mb{G}_xd\xi.
		\end{array}
		\right.
	\end{equation}
	Note that the above fluid-type system is not self-contained and one more  equation for the microscopic component ${\mb{G}}$ is needed, which can be derived by applying the microscopic projection operator $\mb{P}_1$ to \eqref{decompo}:
	\begin{equation}
		\mb{G}_t+\mb{P}_1(\xi_1\mb{M}_x)+\mb{P}_1(\xi_1\mb{G}_x)
		=\mb{L}_\mb{M}\mb{G}+Q(\mb{G}, \mb{G}).
		\label{(1.11)}
	\end{equation}
	Here $\mb{L}_\mb{M}$ is the linearized collision operator of $Q(f,f)$ with respect to the local Maxwellian $\mb{M}$ given by
	$$
	\mb{L}_\mb{M} g:=Q(\mb{M}, g)+ Q(g,\mb{M}).
	$$
	Note that  the null space $\mathfrak{N}$ of $\mb{L}_\mb{M}$ is spanned by the macroscopic variables:
	$$
	\chi_j(\xi), ~j=0,1,2,3,4.
	$$
	Furthermore, there exists a positive constant $\sigma_0(\rho,u,\t)>0$ such that for any function $g(\xi)\in \mathfrak{N}^\bot$, cf. \cite{Grad},
	$$
	\langle g,\mb{L}_\mb{M}g\rangle \leqslant -\sigma_0 \langle
	(1+|\xi|)g,g\rangle.
	$$
	Consequently, the linearized collision operator $\mb{L}_\mb{M}$ is a dissipative operator on $\mathfrak{N}^\bot$, and its inverse	$\mb{L}_\mb{M}^{-1}$  is also a bounded operator on $\mathfrak{N}^\bot$. Their more detailed properties will be discussed in Section 4.
	
	It follows from (\ref{(1.11)}) that
	\begin{equation}\label{(1.12)}
		\mb{G}= \mb{L}_\mb{M}^{-1}[\mb{P}_1(\xi_1\mb{M}_x)] +\Pi,
	\end{equation}
	with
	\begin{equation}\label{(1.13)}
		\Pi:=\mb{L}_\mb{M}^{-1}[\mb{G}_t+\mb
		{P}_1(\xi_1\mb{G}_x)-Q(\mb{G}, \mb{G})].
	\end{equation}
	Plugging  (\ref{(1.12)}) into (\ref{EulerwithG}) gives
	\begin{equation}\label{(1.14)}
		\left\{
		\begin{array}{l}
			\di \rho_{t}+(\rho u_1)_x=0,\\
			\di (\rho u_1)_t+(\rho u_1^2 +p)_x=\f{4}{3}(\mu(\t)
			u_{1x})_x-\int_{\mathbb R^3}\xi_1^2\Pi_xd\xi,  \\[3mm]
			\di (\rho u_i)_t+(\rho u_1u_i)_x=(\mu(\t)
			u_{ix})_x-\int_{\mathbb R^3}\xi_1\xi_i\Pi_xd\xi,~~~~ i=2,3,\\[4mm]
			\di \left[ \rho\left( \t+\f{|u|^2}{2}\right) \right]_t+\left[ \rho
			u_1\left( \t+\f{|u|^2}{2}\right) +pu_1\right] _x=(\k(\t)\t_x)_x+\f{4}{3}(\mu(\t)u_1u_{1x})_x\\[4mm]
			\di\qquad\qquad +\sum_{i=2}^3(\mu(\t)u_iu_{ix})_x
			-\int_{\mathbb R^3}\f12\xi_1|\xi|^2\Pi_xd\xi,
		\end{array}
		\right.
	\end{equation}
	where the viscosity coefficient $\mu(\t)>0$ and the heat conductivity coefficient $\k(\t)>0$ are smooth functions of the temperature $\t$ \cite{Kawashima-Matsumura-Nishida}. Here,  we normalize the gas constant $R$ to be $\f{2}{3}$ so that $e=\t$ and $p=\f23\rho\t$.
	
	Since the problem considered in this paper is  one-dimensional in the space variable $x\in {\mathbb R}$, it is more convenient to rewrite the equation \eqref{(B)} and the system \eqref{Euler} in the {\it Lagrangian} coordinates. For this, set the coordinate transformation
	\begin{equation}\label{Lag}
		(t,x)\mapsto \left( t,\int_{(0,0)}^{(t,x)} \rho(\tau,y)dy-(\r
		u_1)(\tau,y)d\tau\right) ,
	\end{equation}
	where $\di \int_A^B fdy+gd\tau$ represents a line integration from point $A$ to point $B$ on
	${\mathbb R}_+\times {\mathbb R}$. Here, the line integration in \eqref{Lag} is independent of the path and then unique because of the mass conservation law.
	
	We will still denote the {\it Lagrangian} coordinates by $(t,x)$ for the simplicity of notations and let the volume function $v:=\frac1\rho$. Then \eqref{(B)} and \eqref{Euler} in the Lagrangian coordinates become, respectively,
	\begin{equation}\label{Lag-B}
		f_t-\f{u_1}{v}f_x+\f{\xi_1}{v}f_x=Q(f,f),
	\end{equation}
	and
	\begin{equation}\label{(1.16)}
		\left\{
		\begin{array}{llll}
			\di v_{t}-u_{1x}=0,\\[1mm]
			\di u_{1t}+p_x=0,\\[1mm]
			\di u_{it}=0, ~i=2,3,\\[2mm]
			\di \Big( \t+\f{|u|^{2}}{2}\Big)_t+ (pu_1)_x=0.\\
		\end{array}
		\right.
	\end{equation}
	Moreover, (\ref{EulerwithG})-(\ref{(1.14)}) take the form of
	\begin{equation*}
		\left\{
		\begin{array}{llll}
			\di v_t-u_{1x}=0,\\[1mm]
			\di u_{1t}+p_x=-\int_{\mathbb R^3}\xi_1^2\mb{G}_xd\xi,\\[3mm]
			\di u_{it}=-\int_{\mathbb R^3}\xi_1\xi_i\mb{G}_xd\xi,
			~i=2,3,\\[5mm]
			\di\left( \t+\f{|u|^{2}}{2}\right) _{t}+ (pu_1)_x=-\int_{\mathbb R^3}\f12\xi_1|\xi|^2\mb{G}_xd\xi,
		\end{array}
		\right.
	\end{equation*}
	\begin{equation}\label{Lag-G}
		\mb{G}_t-\f{u_1}{v}\mb{G}_x+\f{1}{v}\mb{P}_1(\xi_1\mb{M}_x)+\f{1}{v}\mb{P}_1(\xi_1\mb{G}_x)=\mb{L}_\mb{M}\mb{G}+Q(\mb{G},\mb{G}),
	\end{equation}
	with
	\begin{equation*}
		\mb{G}=\mb{L}^{-1}_\mb{M}\left( \f{1}{v} \mb{P}_1(\xi_1
		\mb{M}_x)\right) +\Pi_1,
	\end{equation*}
	\begin{equation}\label{(1.21)}
		\Pi_1:=\mb{L}_\mb{M}^{-1}\left[ \mb{G}_t-\f{u_1}v\mb{G}_x+\f{1}{v}\mb{P}_1(\xi_1\mb{G}_x)-Q(\mb{G},\mb{G})\right] ,
	\end{equation}
	and
	\begin{equation}\label{(1.22)}
		\left\{
		\begin{array}{llll}
			\di v_t-u_{1x}=0,\\[1mm]
			\di u_{1t}+p_x=\f{4}{3}\left( \f {\mu(\t)}vu_{1x}\right) _{x}-\int_{\mathbb R^3}\xi_1^2\Pi_{1x}d\xi,\\[5mm]
			\di u_{it}=\left( \f{\mu(\t)}{v}u_{ix}\right) _x-\int_{\mathbb R^3}\xi_1\xi_i\Pi_{1x}d\xi,
			~i=2,3,\\[5mm]
			\di\left( \t+\f{|u|^{2}}{2}\right) _{t}+
			(pu_1)_{x}=\left( \f{\k(\t)}{v}\t_x\right) _x+\f{4}{3}\left( \f{\mu(\t)}{v}u_1u_{1x}\right) _x\\[5mm]
			\di\qquad\qquad+\sum_{i=2}^3\left( \f{\mu(\t)}{v} u_iu_{ix}\right) _x
			-\int_{\mathbb R^3}\f12\xi_1|\xi|^2\Pi_{1x}d\xi.
		\end{array}
		\right. 
	\end{equation}
	
	In the present paper, we consider the Boltzmann equation \eqref{Lag-B} with the distinct far-fields  initial data
	\begin{equation}\label{data1}
		f(0,x,\xi)=f_0(x,\xi)\rightarrow \mathbf{M}_{[v_\pm,u_\pm,\t_\pm]}(\xi),\qquad {\rm as}~~x\rightarrow \pm\infty,
	\end{equation}
	where $v_\pm, \t_\pm>0$ and $u_\pm = (u_{1\pm}, 0, 0)$. We focus our attention on the generic case when the solution to the Riemann problem	is the  superposition of a 1-rarefaction wave and a 3-shock wave with a 2-contact discontinuity in between. To this end, let us recall the Riemann problem for the compressible Euler equation (\ref{(1.16)}) with the Riemann initial data
	\begin{equation}\label{Rdata}
		\begin{array}{ll}
			(v,u,\t)(0,x)=\left\{\begin{array}{ll}
				(v_{-},u_-,\t_{-}),~x<0,\\
				(v_{+},u_+,\t_{+}),~x>0,
			\end{array}\right.
		\end{array}
	\end{equation}
	where $u=(u_1,u_2,u_3)$, $u_\pm = (u_{1\pm}, 0, 0)$. It is well known that the Euler system \eqref{(1.16)} for $(v, u_1, \t)$ has three eigenvalues: $\lambda_1=-\sqrt{\frac{5p}{3v}}$, $\lambda_2=0,$ $\lambda_3=\sqrt{\frac{5p}{3v}}$, where the second characteristic field is linearly degenerate and the others two are genuinely nonlinear. The generic Riemann solution of \eqref{(1.16)}, \eqref{Rdata} consists of three basic wave patterns: rarefaction wave, contact discontinuity, and shock wave. To be more specific, given the right end state $(v_+,u_{1+},\t_+)$,  the following wave curves for the left end state  $(v,u_1,\t)$ in the phase space are defined with $v>0$ and $\t>0$ for the Euler system \eqref{(1.16)}.
	
	$\bullet$ $i$-Rarefaction wave curve $(i=1,3)$:
	\begin{equation*}
		R_i (v_+, u_{1+}, \theta_+):=\Bigg{ \{} (v, u_1, \theta)\Bigg{ |}v<v_+ ,~u_1=u_{1+}-\int^v_{v_+}
		\lambda_i(\eta,
		s_+) \,d\eta,~ s(v, \theta)=s_+\Bigg{ \}},
	\end{equation*}
	where $s_+=s(v_+,\t_+)$ and $\l_i=\l_i(v,s)$ is the $i$-th
	characteristic speed of  \eqref{(1.16)}.
	
	$\bullet$ 2-Contact discontinuity curve:
	\begin{equation*}
		CD_2(v_+,u_{1+},\t_+):= \{(v,u_1,\t)  |  u_1=u_{1+}, p=p_+, v
		\not\equiv v_+
		\}. 
	\end{equation*}
	
	$\bullet$ $i$-Shock wave curve $(i=1,3)$:
	\begin{equation*}
		S_i (v_+, u_{1+}, \theta_+):=\Bigg{ \{} (v, u_1, \theta)\Bigg{ |}
		\begin{array}{ll}
			-\sigma_i(v_+-v)-(u_{1+}-u_1)=0,\\[1mm]
			-\sigma_i(u_{1+}-u_1)+(p_+-p)=0,\\[1mm]
			-\sigma_i(E_+-E)+(p_+u_{1+}-pu_1)=0,
		\end{array}
		{\rm and}~~\l_{i+}<\sigma_i<\l_{i}\Bigg{ \}},
	\end{equation*}
	where $E=\t+\f{|u|^2}{2}, p=\f{2\t}{3v},E_+=\t_++\f{|u_+|^2}{2},
	p_+=\f{2\t_+}{3v_+}$, $\l_{i\pm}=\l_i(v_\pm,\t_\pm)$ and $\sigma_i$ is
	the $i$-shock speed.
	
	In this paper, we are interested in the case when the end state $(v_-,u_{1-},\t_-) \in$ $R_1$-$CD_2$-$S_3(v_+,u_{1+},\t_+)$. In such a case, cf. \cite{Smoller}, there exist uniquely two intermediate states $(v_*, u_{1*},\t_*)$ and $(v^*, u_1^*,\t^*)$ such that $(v_-, u_{1-},\t_-)\in R_1(v_*,u_{1*},\t_*)$, $(v_*, u_{1*},\t_*)\in CD_2(v^*, u_1^*,\t^*)$  and
	$(v^*, u_1^*,\t^*)\in S_3(v_+,u_{1+},\t_+)$, at least locally. Our stability result is, roughly speaking, that if the state $(v_-,u_{1-},\t_-)$ $\in$ $R_1$-$CD_2$-$S_3(v_+,u_{1+},\t_+)$, then the solution to the Boltzmann equation \eqref{Lag-B}, \eqref{data1} tends to the superposition of the inviscid rarefaction wave, the viscous contact wave, and the Boltzmann shock profile time-asymptotically, provided that 
	 the conditions in Theorem \ref{maintheorem} hold.

	The rest of the paper will be arranged as follows. In Section 2, we will construct the solution profiles to the Boltzmann equation corresponding to the basic wave patterns to the inviscid Euler system and state the main theorem thereafter. The reformulation of the problem and the main idea for the stability proof will be given in Section 3. Section 4 is dedicated to the properties of the collision operator for later use. The lower and higher order a priori estimates  will be presented in Section 5 and Section 6 respectively. In Appendix, we give the proof of the key Lemma \ref{Lemma-shock} and the local-in time existence Proposition \ref{localexistence} for the solution to the
		Boltzmann equation \eqref{tran-Lag-B} in Lagrangian coordinates.
	
	\emph{Notations.} Throughout this paper, generic positive constants are denoted by $c$ and $C$ if without confusions, which are independent of the small constants $\ve_0, \ve_1, \d_0, \delta_1, \delta_R,  \delta_C, \delta_S$, and the time $T$. For function spaces, $L^{p}(\Omega),1\leqslant p\leqslant\infty$ denotes the usual Lebesgue space on $\Omega
	\subset\mathbb R$ or $\mathbb R^3$ with its norm given by
	$$
	\|f\|_{L^{p}(\Omega)}:=\left(\int_{\Omega}|f(x)|^{p}dx\right)^{\frac{1}{p}},
	\quad 1\leqslant p<\infty, \quad  \parallel
	f\parallel_{L^{\infty}(\Omega)}:=\mbox{ess.sup}_{\Omega} |f(x)|.
	$$
	$W^{k,p}(\Omega)$ denotes the $k^{th}$ order Sobolev space with its norm
	$$
	\|f\|_{W^{k,p}(\Omega)}:=\left(\sum ^{k}_{j=0}
	\parallel \partial^{j}_{x}f\parallel^{p}_{L^{p}(\Omega)}\right)^{\frac{1}{p}}, \quad 1\leqslant p<\infty.
	$$
	And if $p=2$, we note $H^{k}(\Omega):=W^{k,2}(\Omega)$. When the domain $\Omega$ represents the whole space $\mathbb R$ or $\mathbb R^3$, it will be often abbreviated if without confusions. For any function $f : \mathbb R\to \mathbb R$, and $g : \mathbb R_+\times \mathbb R\to \mathbb R$ and any time-dependent shift $\mb X(t)$ (to appear in \eqref{X(t)}),  we denote
	$$f^{- \mb X}(y):=f(y-\mb X(t)), ~~g^{- \mb X}(t,y):=g(t,y-\mb X(t)).$$ 
	Finally, we denote the equivalence of the quantities $A$ and $B$ by $A\sim B$, which means that
	$$0<c\leqslant \left| \frac{A}{B}\right| \leqslant C<+\i$$
holds for two uniform positive constants $c$ and $C$.

	\section{Preliminaries and main result}
	\setcounter{equation}{0}

	In this section, we will first describe the three wave patterns considered in the paper, with the main result provided thereafter.  We start with the rarefaction wave.
	
	\subsection{Approximate rarefaction wave}
	
	First we consider the inviscid rarefaction wave. If $(v_-, u_{1-}, \theta_-) \in R_1 (v_*, u_{1*}, \theta_*)$, then there exists  a self-similar $1$-rarefaction wave fan
	\begin{equation*}
		(v^{r}, u^{r},
		\t^{r})=(v^{r}, u^{r},
		\t^{r})\left( \frac{x}{t}\right) ,
	\end{equation*}
	which is a global entropic solution to the following Riemann problem \cite{Smoller}
	\begin{eqnarray*}
		\left\{
		\begin{array}{l}
			\di  v_{t}- u_{1x}= 0,\\[1mm]
			\di u_{1t} +  p_{x} = 0 ,\\[1mm]
			\di u_{it}=0,~i=2,3,\\[1mm]
			\di \Big(\t+\frac{1}{2}u_1^2 \Big)_t + (pu_1)_x =0,\\[1mm]
			\di  (v, u, \t)(0,x)=\left\{
			\begin{array}{l}
				\di (v_-, u_{-}, \t_-),  ~ x<
				0 ,\\[1mm]
				\di  (v_*, ~u_{*}, ~\theta_*),~ x> 0,
			\end{array}
			\right.
		\end{array} \right.
	\end{eqnarray*}
	where $u_{-}=(u_{1-},0,0)$, $u_{*}=(u_{1*},0,0)$.
	
	Since there is no exact rarefaction wave profile for either the Navier-Stokes equations or
	the Boltzmann equation, the following approximate rarefaction wave profile satisfying the Euler equations was introduced in \cite{Matsumura-Nishihara-1986,Xin-1993}. The construction and properties of approximate rarefaction wave are based on the following inviscid Burgers equation
	\begin{eqnarray}\label{(2.11)}
		\left\{
		\begin{array}{l}
			\di w_{t}+ww_{x}=0,\\
			\di w( 0,x
			)=w_1(x)=\f{w_++w_-}{2}+\f{w_+-w_-}{2}\tanh x.
		\end{array}
		\right.
	\end{eqnarray}
	Note that the solution $w^R(t,x)$ of the problem (\ref{(2.11)}) is given by
	$$
	w^R(t,x)=w_1(x_0(t,x)),\qquad x=x_0(t,x)+w_1(x_0(t,x))t.
	$$
	
	The smooth approximate rarefaction wave profile denoted by $(v^{R},u^{R},\theta^{R})(t,x)$ can be defined by
	\begin{eqnarray}\label{(2.2)}
		\left\{
		\begin{array}{l}
			\di w_-=\l_{1-}:=\l_1(v_-,\t_-), ~ w_*=\l_{1*}:=\l_1(v_*,\t_*), \\[2mm]
			\di \l_1(v^{R}(t,x),\t^{R}(t,x))=w^R(t+1,x),\\[1mm]
			\di u^{R}_1(t,x)=u_{1*}-\int^{v^{R}}_{v_*}  \l_1(v,s_*)
			dv,\\[2mm]
			\di s(v^{R}(t,x),\t^{R}(t,x))=s_*,\\[2mm]
			\di u^{R}_i(t,x)\equiv0,~i=2,3,
		\end{array} \right.
	\end{eqnarray}
	where $s_*=s(v_*,\t_*)$. 	One can easily check that the above approximate rarefaction wave $(v^R, u^R, \theta^R)$ satisfies the system
	\begin{eqnarray}
		\begin{cases}
			v^{R}_t-u^{R} _{1x} = 0,     \cr
			u^{R}_{1t}+p^{R}_x
			=0,\cr
			u^{R}_{i}=0,~i=2,3,\cr
			\t^R_t
			+ p^{R} u^{R}_{1x}
			=0,
		\end{cases}\label{rarefaction-equ}
	\end{eqnarray}
	where $ p^{R}=p( v^{R}, \t^{R})=\f{2\t^{R}}{3v^{R}}$.
	
	The properties of the approximate rarefaction wave can be summarized as follows.
	\begin{lemma} \label{rarefaction}
		(see \cite{Xin-1993}) Let $\delta_R$ denotes the rarefaction wave strength as $\delta_R := |v_* - v_-|\sim |u_{1*}-u_{1-}|\sim |\theta_*-\theta_-|$. The smooth approximate 1-rarefaction wave $(v^R, u_1^R, \theta^R)(t,x)$ defined in \eqref{(2.2)} satisfies the following properties.
		\begin{itemize}
			\item [(1)] $ u^R_{1x} = \frac{3 v^R}{4}w_x > 0$, $ v^R_x = \frac{3v^R}{\sqrt{10\theta^R}} u^R_{1x}>0$ and $\theta^R_x=-\frac{2\theta^R}{3v^R}v^R_x<0$,  $\forall x \in \mathbb {R},\ t \geqslant 0$.
			
			\item [(2)] The following estimates hold for all $t \geqslant 0$ and $p \in [1, + \infty]$:
			\begin{align*}
				&\|(v^R_x, u^R_{1x}, \theta^R_x)\|_{L^p} \leqslant C \min\{\delta_R, \delta_R^{1/p}(1+t)^{-1+1/p}\}, \\[1mm]
				&\|(v^R_{xx}, u^R_{1xx}, \theta^R_{xx})\|_{L^p} \leqslant C \min\{\delta_R, (1+t)^{-1}\}, \\[1mm]
				&\|(v^R_{xxx}, u^R_{1xxx}, \theta^R_{xxx})\|_{L^p} \leqslant C \min\{\delta_R, (1+t)^{-1}\}, \\[1mm]
				& |u^R_{1xx}| \leqslant C |u^R_{1x}|,\quad |\theta^R_{xx}| \leqslant C |\theta^R_{x}| ,\quad \forall x\in\mathbb {R}.
			\end{align*}
			
			\item [(3)] For $x\geqslant \lambda_{1*}(1+t), t\geqslant 0,$ it holds that
			\begin{align*}
				&|(v^R, u_1^R, \theta^R)(t,x)-(v_*, u_{1*}, \theta_*)| \leqslant C\delta_R \ e^{-2|x-\lambda_{1*}(1+t)|}, \\[1mm]
				&|(v^R_x, u^R_{1x}, \theta^R_x)(t,x)|\leqslant C\delta_R\ e^{-2|x-\lambda_{1*}(1+t)|}.
			\end{align*}
			
			\item [(4)] For $x\leqslant\lambda_{1-}(1+t), t\geqslant 0,$ it holds that
			\begin{align*}
				&|(v^R, u_1^R, \theta^R)(t,x)-(v_-,u_{1-}, \theta_-)| \leqslant C\delta_R \ e^{-2|x-\lambda_{1-}(1+t)|}, \\[1mm]
				&|(v^R_x, u^R_{1x}, \theta^R_x)(t,x)|\leqslant C\delta_R\ e^{-2|x-\lambda_{1-}(1+t)|}.
			\end{align*}
			
			\item [(5)] The smooth approximate rarefaction wave and the inviscid rarefaction wave are equivalent time-asymptotically:
			\begin{equation}\label{Randr}
				\di \lim_{t \to +\infty} \sup_{x \in \mathbb {R}} \left| (v^R, u^R_1, \theta^R)(t,x) - \left( v^r, u_1^r, \theta^r\right) \left(  \frac xt\right) \right| = 0.
			\end{equation}
		\end{itemize}
	\end{lemma}

	\subsection{Viscous contact wave}
	For the Euler system \eqref{(1.16)} with a Riemann initial data
	\begin{equation}\label{(2.26)}
		(v,u,\t)(0,x)=\left\{
		\begin{array}{l}
			(v_*,u_*,\t_*),~~~x<0,\\[1mm]
			(v^*,u^*,\t^*),~~~x>0,
		\end{array}
		\right.
	\end{equation}
	where $u_{*}=(u_{1*},0,0)$, $u^{*}=(u_1^{*},0,0)$. It is known (cf. \cite{Smoller}) that the Riemann problem (\ref{(1.16)}), (\ref{(2.26)}) admits a contact discontinuity solution
	\begin{equation*}
		(v^{c},u^{c},\t^{c})(t,x)=\left\{
		\begin{array}{l}
			(v_*,u_*,\t_*),~~~x<0,\\[1mm]
			(v^*,u^*,\t^*),~~~x>0,
		\end{array}
		\right.
	\end{equation*}
	provided that $(v_*,u_{1*},\t_*)\in CD_2(v^*,u^*_1,\t^*)$, that is,
	\begin{equation*}
		u_{1}^*=u_{1*},\qquad p_*:=\f{2\t_*}{3v_*}=\f{2\t^*}{3v^*}=:p^*.
		\label{(2.28)}
	\end{equation*}
	
	The corresponding viscous contact wave $\left( v^{C}, u^{C},\theta^{C}\right) (t,x)$ for the Boltzmann equation can be constructed by using a similar technique for the Navier-Stokes-Fourier equations as in \cite{Huang-Xin-Yang}. Precisely, let
	\begin{equation}\label{vc}
		\begin{array}{ll}
			\di  v^{C}\Big( \frac{x}{\sqrt{1+t}}\Big) 
			:=  \frac{2\Theta^{\rm sim}}{3p_*},\\[4mm]
			\di u_1^{C}\Big( t,\frac{x}{\sqrt{1+t}}\Big) 
			:= u_* +\frac{3\kappa(\Theta^{\rm sim})\Theta^{\rm sim}_x}{5\Theta^{\rm sim}},\\[4mm]
			\di u_i^{C}\Big( t,\frac{x}{\sqrt{1+t}}\Big) 
			\equiv0, ~i=2,3,\\[4mm]
			\di \theta^{C}\Big( \frac{x}{\sqrt{1+t}}\Big) 
			:= \Theta^{\rm sim},
		\end{array}
	\end{equation}
	where $\Theta^{\rm sim}=\Theta^{\rm sim}\left(\frac{x}{\sqrt{1+t}}\right)$ is the unique
	self-similar solution to the following nonlinear diffusion equation (cf.  \cite{Atkinson-Peletier})
	\begin{eqnarray*}
		\begin{cases}
			\di\Theta_t= \frac{9p_*}{10}\left(\frac{\kappa(\Theta)\Theta_x}{\Theta}\right)_x, \cr
			\di\Theta(t,-\infty)=\theta_{*},\ \  \Theta(t,+\infty)=\theta^{*}.
		\end{cases}
	\end{eqnarray*}
	
	The properties of the viscous contact wave defined above can be summarized as follows.
	
	\begin{lemma}\label{contact}
		(see \cite{Huang-Xin-Yang}) Let $\delta_C$ denotes the rarefaction wave strength as $\delta_C:= |v^* - v_*|\sim |\theta^*-\theta_*|$. The viscous contact wave  $\left( v^C, u_1^{C},\theta^{C}\right) (t,x)$ defined in \eqref{vc} satisfies
		\begin{equation}\label{vc-pe}
			\begin{array}{ll}
				\di\big(v^C-v_*, u_1^C-u_{1*},\theta^{C}-\theta_* \big)= O(1)\delta_{_C}e^{ -\frac{c_0 x^2}{1+t}}, \quad \forall x<0,\\[2mm]
				\di\big(v^C-v^*, u_1^C-u_1^*,\theta^{C}-\theta^* \big)= O(1)\delta_{_C}e^{ -\frac{c_0 x^2}{1+t}}, \quad \forall x>0,\\[2mm]
				\di (\partial_x^nv^{C},\partial_x^n\theta^{C})(t,x)
				=O(1)\delta_{_C}(1+t)^{-\frac{n}{2}}e^{ -\frac{c_0 x^2}{1+t}}, \qquad\forall x\in \mathbb{R}, \quad n=1,2,\cdots,\\[2mm]
				\di \partial_x^n u_1^{C}(t,x)
				=O(1)\delta_{_C}(1+t)^{-\frac{1+n}{2}}e^{ -\frac{c_0 x^2}{1+t}},\qquad\qquad\ \   \forall x\in \mathbb{R},\quad n=1,2,\cdots,\\
			\end{array}
		\end{equation}
		where $c_0>0 $ is a generic constant. 
	\end{lemma}
	
	The viscous contact wave  $\left( v^C, u^{C},\theta^{C}\right) (t,x)$ defined in \eqref{vc} satisfies the system
	\begin{eqnarray}\label{contact-equ}
		\begin{cases}
			v^{C}_t-u^{C}_{1x} = 0, \cr
			\di u_{1t}^{C}+ p^{C}_x
			=\frac{4}{3}\Big( \frac{\mu(\theta^{C})u^{C}_{1x}}{ v^{C}} \Big) _x+Q^C_1, \cr
			u_{i}^{C}=0,~i=2,3,\\
			\di	\theta^{C}_t+  p^{C} u^{C}_{1x}
			=\Big( \frac{\kappa(\theta^{C})\theta^{C}_{x}}{v^{C}} \Big) _x+\frac{4}{3}\mu(\theta^{C})\frac{(u^{C}_{1x})^2}{v^{C}}+Q^C_2,
		\end{cases}
	\end{eqnarray}
	where $ p^{C}=\frac{2\theta^{C}}{3v^C}$ and
	\begin{equation}\label{QC12}
		\begin{aligned}
			&Q^C_1=u_{1t}^C-\frac{4}{3}\Big( \frac{\mu(\theta^{C})u^{C}_{1x}}{ v^{C}} \Big) _x=O(1)\delta_{C}(1+t)^{-\frac32}e^{-\frac{c_0 x^2}{1+t}},\\[2mm]
			& Q^C_2=-\frac{4}{3}\mu(\theta^{C})\frac{(u^{C}_{1x})^2}{v^{C}}=O(1)\delta_{C}(1+t)^{-2}e^{-\frac{2c_0 x^2}{1+t}},
		\end{aligned}
	\end{equation}
	as $x\rightarrow\pm \infty$ due to Lemma \ref{contact}.

	\subsection{Boltzmann shock profile}
	In this subsection, we consider the shock profile $F^{S}(x-\sigma t,\xi)$ of the Boltzmann equation (\ref{Lag-B}) in Lagrangian coordinates. The construction of the shock profile for the Boltzmann equation is quite different from that for the Navier-Stokes equations due to the microscopic effects. The existence of Navier-Stokes shock profiles can be reduced to the analysis of an autonomous ODE system and Gilbarg \cite{Gilbarg} proved the existence and uniqueness (up to a translation) of the Navier-Stokes profile even with large amplitude. Note that the general center manifold theory can also be applied to obtain the existence of weak Navier-Stokes shock profile for the autonomous ODE system.  While the Boltzmann shock profile satisfies an infinite-dimensional differential equation and even for the weak Boltzmann shock, the standard center manifold theory based on the spectral information can not be applied straightforwardly.

	By using the bifurcation arguments and Lyapunov-Schmidt reduction, Nicolaenko and Thurber \cite{Nicolaenko-Thurber} first proved the existence of weak Boltzmann shock profile 
	for the hard sphere model and then Caflisch and Nikolaenko \cite{Caflisch-Nicolaenko} generalized it to the hard potential case in 1982. In 2004, Liu-Yu \cite{Liu-Yu-2004} provided a new approach for the existence and positivity of shock profiles based on the micro-macro decomposition and weighted energy method, which can also be utilized in studying the hydrodynamic limit of weak shock to the compressible Euler equations \cite{Yu-2005}. In 2009, M\'etivier-Zumbrun \cite{Metivier-Zumbrun} proved the existence of weak positive Boltzmann shock profile by a weighted $H^s$-contraction mapping argument. 
	
	Very recently, Liu-Yu \cite{Liu-Yu-2013} and Pogan-Zumbrun \cite{Pogan-Zumbrun-2018,Pogan-Zumbrun-2019} analyze the center manifolds for the steady Boltzmann shock profile by different methods. Liu-Yu \cite{Liu-Yu-2013} develop a time-asymptotic method based on the pointwise Green’s functions for the construction of the invariant manifolds, while Pogan-Zumbrun \cite{Pogan-Zumbrun-2018,Pogan-Zumbrun-2019} present dynamical systems tools for degenerate evolution equations including the steady Boltzmann equation. Both methods can imply the key monotonicity property of characteristic function along the Boltzmann shock profile from the dynamics of the Burgers-like equation governing the flow on the invariant manifolds.
	
	Note that the results mentioned above are proven in Eulerian coordinates. In fact, under the Lagrangian transformation (\ref{Lag}), the shock profile in the Eulerian coordinates can be transformed to that in Lagrangian coordinates; see \cite{Huang-Wang-Wang-Yang} for details. 
	
	Now we clarify some basic facts of Boltzmann shock profile $F^{S}(x-\s t,\xi)$ in Lagrangian coordinates. First of all, set $y:=x-\s t$, then $F^{S}(y,\x)$ satisfies
	\begin{equation}\label{BS}
		\left\{
		\begin{array}{ll}
			\di -\s(F^{S})^\prime-\frac{u_1^S}{v^S}(F^{S})^\prime+\frac{\xi_1}{v^S}( F^{S})^\prime=Q(F^{S},F^{S}),\\[3mm]
			F^{S}(-\i,\xi)=\mathbf{M}_{[v^*,u^*,\t^*]}(\xi), \quad F^{S}(+\i,\xi)=\mathbf{M}_{[v_+,u_+,\t_+]}(\xi),
		\end{array}
		\right.
	\end{equation}
	where $(\cdot)^\prime =\f{d(\cdot)}{dy}$, $u^*=(u_{1}^*,0,0)$,  and $u_+=(u_{1+},0,0)$. Note that $(v^{*},u^{*}_1,\t^*)$ and $(v_{+},u_{1+},\t_+)$ satisfy the Rankine-Hugoniot condition
	
	\begin{equation}\label{rh}
		\left\{
		\begin{array}{ll}
			-\sigma(v_+-v^*)-(u_{1+}-u_1^*)=0,\\[1mm]
			-\sigma(u_{1+}-u_1^*)+(p_+-p^*)=0,\\[1mm]
			-\sigma(E_+-E^*)+(p_+u_{1+}-p^*u_1^*)=0,
		\end{array}
		\right.
	\end{equation}
	and Lax entropy condition
	\begin{equation}
		\l_{3+}<\s<\l_{3}^{*}:=\s^*,
		\label{Lax-E}
	\end{equation}
	with $\s$ being the 3-shock wave speed, $E^*=\t^*+\f{|u^*|^2}{2},$ $p^*=\f{2\t^*}{3v^*}$, $E_+=\t_++\f{|u_+|^2}{2}$, 
	$p_+=\f{2\t_+}{3v_+}$, $\s^*:=\lambda_3^*=\sqrt{\frac{5p^*}{3v^*}}$, and $\lambda_{3+}=\sqrt{\frac{5p_+}{3v_+}}$. It follows from \eqref{rh} that
	\begin{equation}\label{sigmadef}
		\sigma=\sqrt{-\frac{p_+-p^*}{v_+-v^*}}.
	\end{equation}
	By \eqref{Lax-E}, we can see that
	\begin{equation}\label{sm1}
		|\sigma-\sigma^*|\leqslant |\l_{3+}-\l_{3}^{*}|=O(\delta_S).
	\end{equation}
	
	Similar to (\ref{macro-micro}), by the micro-macro decomposition  around the local Maxwellian
	$\mb{M}^{S}$, we have
	$$
	F^{S}(x-\s t,\x)=\mathbf{M}^{S}(x-\s t,\x)+\mathbf{G}^{S}(x-\s t,\x),
	$$
	where
	$$
	\mathbf{M}^{S}(x-\s t,\x):=\f{1}{v^{S}(x-\s t)\sqrt{(2\pi
			R\t^{S}(x-\s t))^3}}e^{-\f{|\x-u^{S}(x-\s t)|^2}{2R\t^{S}(x-\s t)}}.
	$$
	With respect to the inner product $\langle\cdot,\cdot\rangle_{\mathbf{M}^{S}}$ defined in
	\eqref{product},  we can now define the macroscopic projection $\mb{P}^{S}_0$ and microscopic projection $\mb{P}^{S}_1$ by
	\begin{equation*}
		\mb{P}^{S}_0g = {\di\sum_{j=0}^4\langle g,\chi^{S}_j\rangle_{\mathbf{M}^{S}}\chi^{S}_j,}
		\qquad
		\mb{P}^{S}_1g= g-\mb{P}^{S}_0g,
	\end{equation*}
	where $\chi^{S}_j~(j=0,1,2,3,4)$ are the corresponding pairwise
	orthogonal base defined in \eqref{orthogonal-base} by replacing
	$(\r,u,\t, \mb{M})$ by $(\frac{1}{v^{S}},u^{S},\t^{S}, \mb{M}^{S})$.
	
	Under the above micro-macro decomposition, the solution
	$F^{S}=F^{S}(x-\s t,\xi)$ satisfies
	$$
	\mb{P}^{S}_0F^{S}=\mb{M}^{S},~~\mb{P}^{S}_1F^{S}=\mb{G}^{S},
	$$
		In fact, from the invariance of the equation \eqref{BS} by changing
		$\x_i$ with $-\x_i$ and the fact that $u_{i}^*=u_{i+}=0$, we have $\di
		u^{S}_i=\int\xi_1\xi_i\Pi^{S}_xd\xi\equiv0$ for $i=2,3.$ From now on, we denote $u^S:=(u_1^S,0,0)$. Then the macroscopic part satisfies the equation
		\begin{equation}\label{shock-profile}
			\left\{
			\begin{array}{lll}
				\displaystyle -\s v^{S}_y-u_{1y}^{S}=0,\\[1mm]
				\displaystyle
				-\s u_{1y}^{S}+p^{S}_y=\f43\Big( \frac{\mu(\t^{S})u_{1y}^{S}}{v^{S}}\Big) _y-\int\x_1^2\Pi^{S}_{1y}d\x,\\[3mm]
				\di u_{i}^{S}\equiv 0,~~~i=2,3,\\[2mm]
				\displaystyle
				-\s \Big( \t^{S}+\f{|u^{S}|^2}{2}\Big) _y+(p^{S}u_1^{S})_y=\Big( \f{\k(\t^{S})\t^{S}_y}{v^{S}}\Big) _y+\f43\Big( \frac{\mu(\t^{S})u_1^{S}u_{1y}^{S}}{v^{S}}\Big) _y
				-\int\x_1\f{|\x|^2}{2}\Pi^{S}_{1y}d\x,
			\end{array}
			\right.
		\end{equation}
		and the non-fluid component ${\mb{G}}^{S}$ satisfies the equation
		$$
		-\s\mb{G}^{S}_y-\f{u_1^{S}}{v^{S}}\mb{G}^{S}_y+\f{1}{v^{S}}\mb{P}^{S}_1(\xi_1\mb{M}^{S}_y)
		+\f{1}{v^{S}}\mb{P}^{S}_1(\xi_1\mb{G}^{S}_y)
		=\mb{L}_{\mb{M}^{S}}\mb{G}^{S}+Q(\mb{G}^{S},
		\mb{G}^{S}),
		$$
		where $\mb{L}_{\mb{M}^{S}}$ is the linearized collision operator of
		$Q(F^{S},F^{S})$
		with respect to the local Maxwellian $\mb{M}^{S}$:
		$$
		\mb{L}_{\mb{M}^{S}} g:=Q(\mb{M}^{S}, g)+
		Q(g,\mb{M}^{S}).
		$$
		Thus
		\begin{align}\label{GS}
			\di\mb{G}^{S}= \mb{L}_{\mb{M}^{S}}^{-1}\Big[\f{1}{v^{S}}\mb{P}^{S}_1(\xi_1\mb{M}^{S}_y)\Big] +\Pi^{S}_1,
		\end{align}
		\begin{equation}\label{Pi1S}
			\di\Pi^{S}_1:=\mb{L}_{\mb{M}^{S}}^{-1}\Big[-\s\mb{G}^{S}_y-\f{u_1^{S}}{v^{S}}\mb{G}^{S}_y
			+\f{1}{v^{S}}\mb{P}^{S}_1(\xi_1\mb{G}^{S}_y)-Q(\mb{G}^{S},\mb{G}^{S})\Big].
		\end{equation}
		Integrating the system \eqref{shock-profile} over $(-\i,y]$ gives that
		\begin{equation*}
			\left\{
			\begin{array}{ll}
				\di -\sigma (v^S-v^*)-( u_1^S-u_1^*)=0,\\[3mm]
				\di \frac{4}{3}\mu(\t^S)\frac{( u_1^S)^\prime}{v^S}=-\sigma( u_1^S-u^*_1)+(p^S-p^*)+\int\xi_1^2\Pi_1^Sd\x,\\[3mm]
				\di \kappa(\t^S)\frac{( \theta^S)^\prime}{ v^S}+\f43\mu(\t^S)\frac{u_1^S( u_1^S)^\prime}{ v^S}=-\sigma\left(  \t^S+\f12|u^S|^2-\t^*-\f12|u^*|^2\right)\\[3mm]
				\di \qquad\qquad\qquad\qquad\qquad+(p^S u_1^S-p^*u_1^*)+\int\xi_1\frac{|\xi|^2}{2}\Pi^S_1d\xi,
			\end{array}
			\right.
		\end{equation*}
		which implies the plane dynamical system
		\begin{equation}\label{VS2}
			\left\{
			\begin{array}{ll}
				\di -\f43\mu(\t^S)\sigma\frac{( v^S)^\prime}{v^S}=(p^S-p^*)+\sigma^2(v^S-v^*)+\int\xi_1^2\Pi_1^Sd\x,\\[3mm]
				\di -\kappa(\t^S)\frac{( \theta^S)^\prime}{\sigma v^S}=(\theta^S-\theta^*)+p^*(v^S -v^*)-\frac12\sigma^2(v^S-v^*)^2\\[3mm]
				\di \qquad\qquad\qquad\qquad\qquad-\frac{1}{\s}\int\xi_1\left( \frac{1}{2}|\xi|^2-u^S_1 \xi_1\right) \Pi^S_1d\xi.
			\end{array}
			\right.
		\end{equation}

		Now we can state some important properties of Boltzmann shock profile $F^{S}(x-
		\sigma t,\x)$ that are given or can be induced by \cite{Liu-Yu-2013, Pogan-Zumbrun-2018}.

		\begin{lemma}\label{Lemma-shock}
			Let $\delta_S$ denote the shock wave strength as $\delta_S := |v_+ - v^*|\sim |u_{1+}-u_1^*|\sim |\t_+ - \t^*|$. If $(v^*,u_1^*,\t^*)\in S_3(v_+,u_{1+},\t_+)$ and the shock wave strength $\d_{S}$ is suitably small, then there exists a unique shock profile $F^{S}(y,\x)$ with $y=x-\sigma t$ up to a shift, to the Boltzmann equation \eqref{Lag-B}. Moreover, there exist positive constants $c$ and $C$ such that
			the following properties hold:
			\begin{itemize}
				\item [(1)] The macroscopic variables $(v^S,u_1^S,\t^S)$ satisfy:
				\begin{align}
					\begin{aligned}\label{shock-macro}
						& v^S_y >0, \qquad u^S_{1y}<0, \qquad \theta^S_y<0,\quad  \forall y\in\mathbb{R},\\
						& | (v^S -v^*, u_1^S -u^*, \theta^S -\theta^*)|\leqslant C\delta_S\ e^{-c\delta_S |y|}, \quad y<0,\\[1mm]
						& | (v^S -v_+, u_1^S -u_+, \theta^S -\theta_+)|\leqslant C\delta_S\ e^{-c\delta_S |y|}, \quad y>0,\\[1mm]
						&|( v^S_y, u^S_{1y},\theta^S_y)|\leqslant C\delta_S|(v^S-v^*,u_1^S -u^*, \theta^S -\theta^*)|, \\[1mm]
						&|( v^S_y, u^S_{1y},\theta^S_y)|\leqslant C\d_S^2\ e^{-c\delta_S |y|},\\[1mm]
						& \di |\partial^k_y(v^{S},u^{S}_{1},\t^{S})|\leqslant
						C\d_{S}^{k-1}|(v^{S}_y,u^S_{1y},\theta^S_y)|, \quad k\geqslant 2.
					\end{aligned}
				\end{align}		
				\item[(2)] The microscopic variable $\mb{G}^{S}$ and $\Pi_1^S$ satisfy:
				\begin{align}
					\begin{aligned}\label{shock-micro}
						& \di\left( \int
						\f{(1+|\xi|)|\mb{G}^{S}|^2}{\mb{M}_0}d\x\right) ^{\f12}\leqslant
						C\d_{S}^2 e^{-c\d_{S}|y|},\\
						& \di \left( \int
						\f{(1+|\xi|)|\partial^k_y\mb{G}^{S}|^2}{\mb{M}_0}d\x\right) ^{\f12}\leqslant
						C\d_{S}^k\left( \int
						\f{(1+|\xi|)|\mb{G}^{S}|^2}{\mb{M}_0}d\x\right) ^{\f12},~~k\geqslant 1 , k\in\mathbb{N},\\[1mm]
						&\di \left| \int\x_1\varphi_i(\x)\Pi_1^{S}d\x\right|  \leqslant
						C\d_{S}\left( \int
						\f{(1+|\xi|)|\mb{G}^{S}|^2}{\mb{M}_0}d\x\right) ^{\f12},~~i=1,4,\\[1mm]
						&\di\left| \int\x_1\varphi_i(\x)\Pi^{S}_{1y}d\x\right| \leqslant
						C\d_{S}^2\di\left( \int
						\f{(1+|\xi|)|\mb{G}^{S}|^2}{\mb{M}_0}d\x\right) ^{\f12},~~i=1,4,\\[1mm]
						&\di \int\x_1\varphi_i(\x)\Pi_1^{S}d\x =0,~~i=2,3,
					\end{aligned}
				\end{align}		
				where $\mb{M}_0$ is a global Maxwellian which is close to the shock profile
				with its precise definition given in \cite[Theorem 21]{Liu-Yu-2013}, and $\varphi_i(\x)~(i=1,2,3,4)$ are the collision invariants defined in \eqref{collision-invar}.
				
				\item[(3)] The relation between macroscopic variables $(v^S,u_1^S,\t^S)$ and microscopic variable $\mb{G}^{S}$ can be expressed as:
				\begin{equation}\label{equivalence}
					\begin{array}{ll}
						\di v^{S}_y\sim u^{S}_{1y}\sim
						\t^{S}_y\sim \left( \int
						\f{(1+|\x|)|\mb{G}^{S}|^2}{\mb{M}_0}d\x\right) ^{\f12}.
					\end{array}
				\end{equation}
				
				\item[(4)] The more precise relation between the macroscopic variables $(v^S,u_1^S,\t^S)$ can be expressed as:
				\begin{equation}\label{shock-vu}
					\left|  u^S_{1y} + \sigma^* v^S_y\right| \leqslant C \delta_Sv^S_y, \quad\forall y\in\mathbb{R},
				\end{equation}
				and
				\begin{equation}\label{theta-s}
					\left|  \theta^S_y + p^*v^S_y \right|  \leqslant C \delta_Sv^S_y,\quad\forall y\in\mathbb{R}.
				\end{equation}

				\item[(5)] There exists a positive constant $C$ such that
				\begin{equation}\label{2ndorder}
					\begin{array}{ll}
						\di \left| \frac{p^S-p_+}{v^S-v_+}-\frac{p^S-p^*}{v^S-v^*}- \frac{ 5p^*}{9(v^*)^2}\bigg(\frac{10\mu(\t^*)-9\k(\t^*)}{10\mu(\t^*)+3\k(\t^*)} +3\bigg) \d_S\right| \leqslant C\delta_S^2.
					\end{array}
				\end{equation}
			\end{itemize}
		\end{lemma}
		\begin{remark}
			The expansion estimate \eqref{2ndorder} includes the microscopic effects of Boltzmann shock profile and is quite different from Navier-Stokes-Fourier shock profile in \cite{Kang-Vasseur-Wang-2024} and crucial in the following energy anlaysis.
		\end{remark}
		\begin{remark}
			Roughly speaking, \eqref{equivalence} means that the microscopic part of Boltzmann shock profile is equivalent to the first-order derivative of the macroscopic part. That is, Boltzmann shock profile can be seen as a small perturbation of Navier-Stokes-Fourier shock profile, even though the higher order  microscopic part of Boltzmann shock profile is essential for its nonlinear stability. 
		\end{remark}	
		
		The proof of Lemma \ref{Lemma-shock} will be given in Appendix A.1.
		
		\subsection{Main result}
		Now we state our main result as follows.
		\begin{theorem}\label{maintheorem}
			For each $(v_+,u_+,\theta_+)$ with $v_+,\t_+>0$, let $(v_-,u_-,\t_-)$ satisfies that the Riemann problem of Euler equation \eqref{(1.16)} consists of a 1-rarefaction wave, a 2-contact discontinuity, and a 3-shock wave with two intermediate states $(v_*,u_*,\t_*)$ and $(v^*,u^*,\t^*)$. Then there exist two positive constants $\delta_0,\varepsilon_0$ and a global Maxwellian $\mb M_{_\#}=\mb M_{[v_{_\#},u_{_\#},\t_{_\#}]}$ with $v_{_\#}, \t_{_\#}>0$ such that if the wave strength $\delta_R+\delta_C+\delta_S\leqslant\delta_0$ and the initial data $f_0(x,\xi)$ satisfies $f_0(x,\xi)\geqslant0$ and
			\begin{equation}\label{epsilon0}
				\begin{array}{ll}
					\di \sum_\pm\left\| f_0(x,\xi)-\mb M_{[v_\pm,u_\pm,\t_\pm]}(\xi) \right\|_{L^2_x\left(\mathbb R_\pm,L^2_\xi\left(\frac{1}{\sqrt{\mb M_{_\#}}} \right)  \right) }+\left\| (f_{0x},f_{0t})(x,\xi) \right\|_{H^1_x\left(L^2_\xi\left(\frac{1}{\sqrt{\mb M_{_\#}}} \right)  \right) }\leqslant\ve_0,
				\end{array}
			\end{equation}
			then the Boltzmann equation \eqref{Lag-B}, \eqref{data1} admits a unique global solution $f(t,x,\xi)\geqslant 0$ for all $t\in\mathbb R_+$. Moreover, there exists an absolutely
			continuous shift $\mb X(t)$, such that
			\begin{equation}\label{mainlimit}
				\begin{array}{ll}
					\di \lim\limits_{t\to +\infty}\bigg\| f(t,x,\xi)-\Big(\mb M_{\left[ v^r\left(\f xt \right) ,u^r\left(\f xt \right),\t^r\left(\f xt \right)\right] }(\xi)+ \mb M_{\left[ v^C\left(\frac{x}{\sqrt{1+t}} \right) ,u^C\left(t,\frac{x}{\sqrt{1+t}} \right),\t^C\left(\frac{x}{\sqrt{1+t}} \right)\right] }(\xi) \\[6mm]
					\qquad+F^S\big(x-\s t-\mb X(t),\xi\big)-\mb M_{[v_*,u_*,\t_*]}(\xi)-\mb M_{[v^*,u^*,\t^*]}(\xi)\Big)\bigg\|_{L^\infty_x\left( L^2_\xi\left( \frac{1}{\sqrt{\mb M_{_\#}}} \right)  \right)  } =0,
				\end{array}
			\end{equation}
			and
			\begin{equation}\label{Xdotlimit}				
				\lim\limits_{t\to +\infty} \dot{\mb X}(t) =0.
			\end{equation}
			Here the norm $\|\cdot\|_{L_\xi^2\left( \f{1}{\sqrt{\mb{M}_{_\#}}}\right) }$ means
			$\|\f{\cdot}{\sqrt{\mb{M}_{_\#}}}\|_{L_\xi^2(\mathbb R^3)}$.
		\end{theorem}
		\begin{remark}
			Theorem \ref{maintheorem} states that if the two far-field states $(v_\pm,u_\pm,\t_\pm)$ are connected by the generic Riemann solution consisting a rarefaction wave, a contact discontinuity, and a shock wave, to the compressible Euler equations,  then the global solution to Boltzmann equation \eqref{Lag-B}  converges to the superposition wave of inviscid self-similar rarefaction wave, the viscous contact wave, and the Boltzmann shock profile with the time-dependent shift $\mb X(t)$. Roughly speaking, the generic Riemann solution consisting of a rarefaction wave, a contact discontinuity, and a shock wave is time-asymptotically stable to the one-dimensional Boltzmann equation.
			
		\end{remark}
		\begin{remark}
			The shift $\mb X(t)$ (defined in \eqref{X(t)}) is proved to satisfy the time-asymptotic behavior \eqref{Xdotlimit}, which implies that
			$$\lim\limits_{t\to +\infty}\frac{\mb X(t)}{t}=0.$$
			That is, the shift function $\mb X(t)$ grows sub-linearly with respect to the time $t$ and then the shifted Boltzmann shock profile still keeps the original traveling wave shape time-asymptotically.
		\end{remark}
		
		\
		
		\section{Proof of main result}
		\setcounter{equation}{0}
		\subsection{Change of variables}
		For simplicity, we rewrite \eqref{(1.22)} in non-divergence form, simultaneously performing the change of variables associated to the speed of propagation of the shock $(t,x)\mapsto (t,y=x-\s t)$:
		\begin{equation}\label{NS-1}
			\left\{
			\begin{array}{lll}
				\displaystyle v_t-\s v_y-u_{1y}=0,\\
				\displaystyle
				u_{1t}-\sigma u_{1y}+p_y=\f43\left( \frac{\mu(\t)u_{1y}}{v}\right) _y-\int\x_1^2\Pi_{1y}d\x,\\
				\di u_{it}-\s u_{iy}=\left( \frac{\mu(\t)u_{iy}}{v}\right) _y-\int\x_1\x_i\Pi_{1y}d\x,~~~i=2,3,\\
				\displaystyle
				\t_t-\s \t_y+pu_{1y}=\left( \f{\k(\t)\t_y}{v}\right) _y+\f43\mu(\t) \frac{u_{1y}^2}{v}+\mu(\t) \frac{u_{2y}^2+u_{3y}^2}{v}
				-\int\x_1\left( \f{|\x|^2}{2}-\sum_{i=1}^3u_i\xi_i \right) \Pi_{1y}d\x.
			\end{array}
			\right.
		\end{equation}
		Similarly, we rewrite \eqref{Lag-B} (together with \eqref{data1}), \eqref{Lag-G} and \eqref{(1.21)} as the following, respectively,
		\begin{equation}\label{tran-Lag-B}
			\left\{
			\begin{array}{lll}
				\displaystyle f_t-\s f_y-\f{u_1}{v}f_y+\f{\xi_1}{v}f_y=Q(f,f),\\[3mm]
				\displaystyle
				f(0,y,\xi)=f_0(y,\xi),
			\end{array}
			\right.
		\end{equation}
		\begin{equation}\label{tran-G}
			\mb{G}_t-\s\mb{G}_y-\f{u_1}{v}\mb{G}_y+\f{1}{v}\mb{P}_1(\xi_1\mb{M}_y)+\f{1}{v}\mb{P}_1(\xi_1\mb{G}_y)=\mb{L}_\mb{M}\mb{G}+Q(\mb{G},\mb{G}),
		\end{equation}
		\begin{equation*}
			\Pi_1=\mb{L}_\mb{M}^{-1}\left[ \mb{G}_t-\s\mb{G}_y-\f{u_1}v\mb{G}_y+\f{1}{v}\mb{P}_1(\xi_1\mb{G}_y)-Q(\mb{G},\mb{G})\right].
		\end{equation*}
		
		Now we consider the basic wave patterns under the new variables. It is easy to check that, from \eqref{rarefaction-equ}, the approximate rarefaction wave $(v,u,\t)(t,y)=(v^R, u^R, \theta^R)(t,y+\s t) $ satisfies 
		\begin{equation*}
			\begin{cases}
				\displaystyle v_t -\s  v_y-u_{1y}= 0, \\[1mm]
				\displaystyle u_{1t} -\s  u_{1y}+ p_y = 0,\\[1mm]
				\di u_i=0,~~~~~~~~~~~ i=2,3,\\[1mm]
				\di \theta_t-\s\theta_y+p(v,\t) u_{1y}=0. \\
			\end{cases}
		\end{equation*}
		From \eqref{contact-equ},
		the viscous contact wave  $(v,u,\t)(t,y)=(v^C, u^C, \theta^C)(t,y+\s t)$ satisfies
		\begin{eqnarray}\label{vcex}
			\begin{cases}
				\di v_t-\s v_y-u_{1y} = 0, \cr
				\di u_{1t}-\s u_{1y}+ p_y
				=\frac{4}{3}\left( \frac{\mu(\theta)u_{1y}}{ v} \right)_y+Q^C_1, \cr
				\di u_i=0,\quad i=2,3,\\
				\di \theta_t-\s\theta_y+  p u_{1y}
				=\left( \frac{\kappa(\theta)\theta_{y}}{v} \right) _y+\frac{4}{3}\mu(\theta)\frac{(u_{1y})^2}{v}+Q^C_2.
			\end{cases}
		\end{eqnarray}
		The ansatz we consider is the superposition of the approximate rarefaction wave, the
		viscous contact wave, and the shock profile shifted by $\mb X(t)$ (to be defined in \eqref{X(t)}), which can be expressed in the following form:
		\begin{equation}\label{ansatz}
			\begin{array}{l}
				\left(\begin{array}{cc} \bar v\\[1mm] \bar u_1 \\[1mm] \bar u_2 \\[1mm] \bar u_3 \\[1mm] \bar \t
				\end{array}
				\right)(t,y)= \left(\begin{array}{cc}v^{R}(t,y+\s t)+ v^{C}(y+\s t)+ v^{S}(y-\mb X(t))-v_*-v^*\\ [1mm]u_1^{R}(t,y+\s t)+ u_1^{C}(t,y+\s t)+ u_1^{S}(y-\mb X(t))-u_{1*}-u_1^* \\[1mm] 0 \\[1mm] 0 \\[1mm]
					\t^{R}(t,y+\s t)+ \t^{C}(y+\s t)+ \t^{S}(y-\mb X(t))-\t_*-\t^*
				\end{array}
				\right),
			\end{array}
		\end{equation}
		\begin{equation}\label{barM}
			\bar{\mb M}(t,y,\xi):=\mb M_{[\bar v(t,y),\bar u(t,y),\bar \theta (t,y)]}(\xi),
		\end{equation}
		and
		\begin{equation*}
			\bar f(t,y,\xi):=\bar{\mb M}(t,y,\xi)+\mb G^S(y-\mb X(t),\xi),
		\end{equation*}
		Direct calculation yields that the ansatz $(\bar v,\bar u,\bar \theta)$ satisfies
		\begin{equation}\label{bar-system}
			\left\{
			\begin{array}{ll}
				\di \bar v_t-\s \bar v_y+\dot{\mb X}(t)(v^S)^{-\mb X}_y-\bar u_{1y}=0,\\[3mm]
				\di \bar u_{1t}-\s \bar u_{1y}+\dot{\mb X}(t)(u^S_1)^{-\mb X}_y+\bar p_y=\f43\left( \frac{\mu(\bar\t)\bar u_{1y}}{\bar v}\right) _y+Q_1-\int\x_1^2(\Pi^S_{1})^{-\mb X}_yd\x,\\[3mm]
				\di \bar u_{i}=0,~~~i=2,3,\\
				\di \bar \theta_t-\s\bar \theta_y+\dot{\mb X}(t)(\theta^S)^{-\mb X}_y+\bar p\bar u_{1y}=\left(\frac{\kappa(\bar\t)\bar \theta_y}{\bar v}\right)_y+\f43\mu(\bar\t)\frac{\bar u_{1y}^2}{\bar v}+Q_2\\
				\di\qquad\qquad\qquad\qquad\qquad\qquad-\int\x_1\left( \f12|\x|^2-(u_1^S)^{-\mb X}\xi_1\right) (\Pi^S_{1})^{-\mb X}_yd\x,
			\end{array}
			\right.
		\end{equation}
		where $\bar p=\frac{2\bar\theta}{3\bar v}$ and the error terms are
		\begin{equation*}
			Q_i:=Q^I_i+Q^R_i+Q^C_i, \quad i=1,2,
		\end{equation*}
		with the wave interactions terms
		\begin{equation}\label{QI1}
			\begin{array}{ll}
				Q_1^I:=&\di  \left(\bar p-p^R-p^C-(p^S)^{-\mb X}\right)_y\\[3mm]
				&\di -\f43\left(\frac{\mu(\bar\t)\bar u_{1y}}{\bar v}-\frac{\mu(\t^R)u_{1y}^R}{v^R}-\frac{\mu(\t^C)u^C_{1y}}{v^C}-\frac{\mu((\t^S)^{-\mb X})(u_1^S)^{-\mb X}_y}{(v^S)^{-\mb X}}\right)_y,
			\end{array}
		\end{equation}
		\begin{equation}\label{QI2}
			\begin{array}{ll}
				\di 
				\qquad\ \  Q_2^I:= &\di \left(\bar p\bar u_{1y}-p^Ru^R_{1y}-p^Cu^C_{1y}-(p^S)^{-\mb X}(u_1^S)^{-\mb X}_y\right)\\[3mm]
				&\di -\left(\frac{\kappa(\bar \t)\bar \theta_y}{\bar v}-\frac{\kappa(\t^R)\theta^R_y}{v^R}-\frac{\kappa((\theta^S)^{-\mb X})(\theta^S)^{-\mb X}_y}{(v^S)^{-\mb X}}-\frac{\kappa(\t^C)\theta^C_y}{v^C}\right)_y \\[6mm]
				&\di -\f43\left(\frac{\mu(\bar \t)\bar u^2_{1y}}{\bar v}-\frac{\mu(\t^R)(u^R_{1y})^2}{v^R}-\frac{\mu(\t^C)(u^C_{1y})^2}{v^C}-\frac{\mu((\t^S)^{-\mb X})((u^S)^{-\mb X}_{1y})^2}{(v^S)^{-\mb X}}\right),
			\end{array}
		\end{equation}
		and the error terms due to the inviscid rarefaction wave
		\begin{equation*}
			Q^R_1:=-\f43\left(\frac{\mu(\t^R)u^R_{1y}}{v^R}\right)_y,\qquad
			Q^R_2:=-\left(\frac{\kappa(\t^R)\theta^R_y}{v^R}\right)_y -\f43 \frac{\mu(\t^R)(u^R_{1y})^2}{v^R},
		\end{equation*}
		and the error terms $Q_1^C, Q_2^C$ due to the viscous contact wave given in \eqref{QC12}.

		\subsection{Construction of weight function}
		For the continuation argument, the main tool is the weighted energy estimates (Proposition \ref{main proposition}). Thus we define the weight function $a(y)$ by
		\begin{equation}\label{weight}
			\di a(y):=1+\frac{\lambda}{\d_S}(v^S(y)-v^*),
		\end{equation}
		where $\delta_S:=v_+-v^*$ and $\l$ is a small constant satisfying $\d_S\ll\l\leqslant C\sqrt{\d_S}$. From now on we will fix $\lambda=\d_S^{\f34}$ for simplicity. Notice that 
		\begin{equation}\label{abound}
			1<a(y)<1+\delta_S^{\f34}<2,
		\end{equation}
		and 
		\begin{equation}\label{a-prime}
			a^\prime(y)=\delta_S^{-\f14}v^S_y>0,
		\end{equation}
		and so,
		\begin{equation*}
			|a'|\sim\delta_S^{-\f14} |u^S_{1y}|\sim \delta_S^{-\f14} |\theta^S_y|.
		\end{equation*}

		\subsection{Construction of shift}
		Our construction of shift depends on the weight function  $a:\mathbb R\to\mathbb R$ defined in \eqref{weight}. The definition of the shift  $\mb X(t)$ is a solution to the ODE:
		\begin{equation}\label{X(t)}
			\left\{
			\begin{array}{ll}
				\di \dot{\mb{X}}(t)=-\frac{H}{\delta_S} \int_{-\i}^{+\i}a^{-\mb{X}}\left((v^S)^{-\mb{X}}_{y}\frac{\bar p }{\bar v}(v-\bar v)+(u^S_{1})^{-\mb{X}}_y(u_1-\bar u_1)+(\t^S)^{-\mb{X}}_{y}\frac{\t-\bar\t}{\bar\theta}\right) dy,\\[5mm]
				\di \mb X(0)=0,
			\end{array}
			\right.
		\end{equation}
		where $H$ is the specific constant chosen as $
		\di H:=\frac{20p^*}{3(v^*)^2(\sigma^*)^3}\f{5\mu (\t^*)+3\kappa (\t^*)}{10\mu (\t^*)+3\kappa (\t^*)} $, which will be used in Section 4.2.

		The following lemma ensures that \eqref{X(t)} has a unique absolutely continuous solution defined on any interval in time $[0,T]$.

		\begin{lemma}\label{lem:xex}
			For any $c_1,c_2>0$, there exists a constant $C>0$ such that the following is true.  For any $T>0$, and any function $v,\t\in L^\infty ((0,T)\times \mathbb R)$ verifying
			\begin{equation*}
				c_1 \leqslant v(t,y), \theta(t,y)\leqslant c_2,\qquad \forall (t,y)\in [0,T]\times\mathbb R ,
			\end{equation*}
			the ODE \eqref{X(t)} has a unique absolutely continuous solution $\mb X(t)$ on $[0,T]$. Moreover,
			\begin{equation}\label{roughx}
				|{\mb X}(t)| \leqslant Ct,\quad \forall t\in[0,T].
			\end{equation}
		\end{lemma}
		{\it Proof.} The existence and uniqueness of the absolutely continuous solution $\mb X(t)$ is a direct consequence of the well-known Cauchy-Lipschitz Theorem. Moreover, since
		\begin{equation*}
			\begin{array}{l}
				\di |\dot{\mb X}(t)|\leqslant \frac{C}{\delta_S}  \|\left( v-\bar v,u_1-\bar u_1,\t-\bar\t\right)  \|_{L^\infty} \int_{-\i}^{+\i} |( v^S_y, u^S_{1y},\theta^S_y)| dy\leqslant C,
			\end{array}
		\end{equation*}
		we have \eqref{roughx}.

		\subsection{A priori estimates}
		
		Denote the perturbation around the ansatz \eqref{ansatz} by
		\begin{equation}\label{perdef}
			(\phi,\psi,\zeta)(t,y):=(v-\bar v, u-\bar u, \t-\bar \t)(t,y),
		\end{equation}
		\begin{equation}\label{tildeGdef}
			\wt{\mb G}(t,y,\xi):=\mb G(t,y,\xi)-\mb G^S(y-\mb X(t),\xi),
		\end{equation}		
		\begin{equation}\label{tildefdef}
			\tilde{f}(t,y,\xi):=f(t,y,\xi)-F^S(y-\mb X(t),\xi),
		\end{equation}
		where $u=(u_1,u_2,u_3)$, $\psi=(\psi_1,\psi_2,\psi_3)$.

		Subtracting \eqref{bar-system} from \eqref{NS-1} yields the system for perturbation $(\phi,\psi,\zeta)$:
		\begin{equation}\label{pereq1}
			\left\{
			\begin{array}{ll}
				\di \phi_t-\s \phi_y-\dot{\mb X}(t)(v^S)^{-\mb X}_y-\psi_{1y}=0,\\[3mm]
				\di \psi_{1t}-\s \psi_{1y}-\dot{\mb X}(t)(u_1^S)^{-\mb X}_y+(p-\bar p)_y=\f43\left(\frac{\mu(\t)u_{1y}}{v}-\frac{\mu(\bar\t)\bar u_{1y}}{\bar v}\right)_y-Q_1\\[3mm]
				\di~~~~~~~~~~~~~~~~~~~~~~~~~~~~~-\int \xi_1^2\left( \Pi_1-(\Pi_1^S)^{-\mb X}\right)_y d\xi,\\[3mm]
				\di \psi_{it}-\s \psi_{iy}=\left( \frac{\mu(\t)\psi_{iy}}{v}\right) _y-\int\x_1\x_i\Pi_{1y}d\x,~~~i=2,3,\\[3mm]
				\di \zeta_t-\s\zeta_y -\dot{\mb X}(t)(\theta^S)^{-\mb X}_y+(pu_{1y} -\bar p\bar u_{1y})=\left(\frac{\kappa(\t)\theta_y}{v}-\frac{\kappa(\bar\t)\bar \theta_y}{\bar v}\right)_y+\f43\left(\frac{\mu(\t)u_{1y}^2}{v}-\frac{\mu(\bar\t)\bar u_{1y}^2}{\bar v}\right)\\[3mm]
				\qquad\qquad\qquad\di+\frac{\mu(\t)(\psi_{2y}^2+\psi_{3y}^2)}{v}-Q_2-\int\x_1\frac{|\xi|^2}{2}\left( \Pi_{1}-(\Pi_1^S)^{-\mb X}\right)_y d\x+\sum_{i=2}^3\psi_i\int\x_1\xi_i\Pi_{1y}d\x\\[3mm]
				\qquad\qquad\qquad\di +\left[u_1\int\x_1^2\Pi_{1y}d\x-(u_1^S)^{-\mb X}\int\x_1^2(\Pi_1^S)^{-\mb X}_yd\x \right] .
			\end{array}
			\right.
		\end{equation}
		
		We now derive the equation for the non-fluid component $\wt{\mb G}(t,y,\xi)$. From \eqref{tran-G} and \eqref{tildeGdef}, we have
		\begin{equation}\label{tildeG}
			\begin{array}{ll}
				\di \wt{\mb{G}}_{t}-\mb{L}_\mb{M}\wt{\mb{G}}=&\di \sigma \wt{\mb{G}}_y+\dot{\mb X}(t)(\mb{G}^S)^{-\mb X}_y+\frac{u_1}{v} \wt{\mb{G}}_y-\frac{1}{v}\mb{P}_1(\xi_1\wt{\mb{G}}_{y})+\left(\frac{u_1}{v}-\frac{(u_1^S)^{-\mb X}}{(v^S)^{-\mb X}} \right) (\mb{G}^S)^{-\mb X}_y
				\\[4mm]
				&\di-\left( \frac{1}{v}\mb{P}_1(\xi_1(\mb{G}^S)^{-\mb X}_y)-\frac{1}{(v^S)^{-\mb X}}   (\mb{P}^{S}_1)^{-\mb X} (\x_1(\mb{G}^{S})^{-\mb X}_y)\right) \\[4mm]
				&\di-\f1v\mb{P}_1(\xi_1\mb{M}_y)+\f1{(v^S)^{-\mb X}}(\mb{P}^{S}_1)^{-\mb X}(\xi_1(\mb{M}^S)^{-\mb X}_y)
				+Q(\wt{\mb{G}},\wt{\mb{G}})
				\\[4mm]
				&\di+Q(\wt{\mb{G}},(\mb{G}^S)^{-\mb X})+Q((\mb{G}^S)^{-\mb X},\wt{\mb{G}})+(\mb{L}_{\mb{M}}-\mb{L}_{(\mb{M}^S)^{-\mb X}})(\mb{G}^S)^{-\mb X},\\
			\end{array}
		\end{equation}
		where the operator $(\mb{P}^{S}_1)^{-\mb X}$ is the microscopic projection according to local Maxwellian of Boltzmann shock profile with shift $\mb X(t)$. Let
		\begin{equation*}
			\wt{\mb{G}}_0:=\frac{3}{2v\t}\mb{L}_{\mb{M}}^{-1} \mb{P}_1\left[\xi_1\mb{M}\left( \xi_1(u^R_{1y}+u^C_{1y})+\frac{\left|\xi-u \right|^2 }{2\t}(\t^R_{y}+\t^C_{y})\right) \right] 
		\end{equation*}
		and
		\begin{equation}\label{G1}
			\wt{\mb{G}}_1(t,y,\xi):=\wt{\mb{G}}(t,y,\xi)-\wt{\mb{G}}_0(t,y,\xi),
		\end{equation}
		then $\wt{\mb{G}}_1(t,y,\xi)$ satisfies
		\begin{equation}\label{G1eq}
			\begin{array}{ll}
				\di \wt{\mb{G}}_{1t}-\mb{L}_\mb{M}\wt{\mb{G}}_1=&\di \sigma \wt{\mb{G}}_y+\dot{\mb X}(t)(\mb{G}^S)^{-\mb X}_y+\frac{u_1}{v} \wt{\mb{G}}_y-\frac{1}{v}\mb{P}_1(\xi_1\wt{\mb{G}}_{y})+\left(\frac{u_1}{v}-\frac{(u_1^S)^{-\mb X}}{(v^S)^{-\mb X}} \right) (\mb{G}^S)^{-\mb X}_y\\[4mm]
				&\di-\left( \frac{1}{v}\mb{P}_1(\xi_1(\mb{G}^S)^{-\mb X}_y)-\frac{1}{(v^S)^{-\mb X}}   (\mb{P}^{S}_1)^{-\mb X} (\x_1(\mb{G}^{S})^{-\mb X}_y)\right) \\[4mm]
				&\di-\frac{3}{2v\t}\mb{P}_1\left[\xi_1\mb{M}\left( \xi\cdot\psi_y+\xi_1(u^S_{1})^{-\mb X }_y+\frac{\left|\xi-u \right|^2 }{2\t}(\z_{y}+(\t^S)^{-\mb X}_{y})\right) \right]\\[4mm]
				&\di+\f1{(v^S)^{-\mb X}}(\mb{P}^{S}_1)^{-\mb X}(\xi_1(\mb{M}^S)^{-\mb X}_y))
				+Q(\wt{\mb{G}},\wt{\mb{G}})
				+Q(\wt{\mb{G}},(\mb{G}^S)^{-\mb X})\\[4mm]
				&\di+Q((\mb{G}^S)^{-\mb X},\wt{\mb{G}})+(\mb{L}_{\mb{M}}-\mb{L}_{(\mb{M}^S)^{-\mb X}})(\mb{G}^S)^{-\mb X}-\wt{\mb{G}}_{0t}.\\
			\end{array}
		\end{equation}
		Notice that in \eqref{G1} and \eqref{G1eq},
		$\wt{\mb{G}}_0$ is subtracted from $\wt{\mb{G}}$ when carrying out the energy estimates because $|(u^R_{1y}, \t^R_{y}, \t^C_{y})|^2$ is not integrable globally in time $t$.
		
		Similarly, we have the equation for $\tilde f(t,y,\xi)$ by \eqref{tildefdef}:
		\begin{equation*}
			\begin{array}{ll}
				\di ~~\tilde f_t-\s\tilde f_y+\frac{\xi_1-u_1}{v}\tilde f_y\di -\dot{\mb X}(t)(F^S)^{-\mb X}_y+\left(\frac{\xi_1-u_1}{v}-\frac{\xi_1-(u_1^S)^{-\mb X}}{(v^S)^{-\mb X}} \right) (F^S)^{-\mb X}_y\\[5mm]
				\di =\mb{L}_\mb{M}\wt{\mb{G}}+Q(\wt{\mb{G}},\wt{\mb{G}})+(\mb{L}_{\mb{M}}-\mb{L}_{(\mb{M}^S)^{-\mb X}})(\mb{G}^S)^{-\mb X}+Q(\wt{\mb{G}},(\mb{G}^S)^{-\mb X})+Q((\mb{G}^S)^{-\mb X},\wt{\mb{G}}).
			\end{array}
		\end{equation*}

		Consider the reformulated system \eqref{pereq1} and \eqref{tildeG}. To prove the global existence on the time interval $[0, T]$, we should close the following a priori estimate. Set the a priori assumption
		\begin{equation}\label{prioriassumption}
			\begin{array}{ll}
				\di \mathcal{N}(T)^2&\di:=\sup_{t\in[0,T]}\Bigg\{
				\|(\phi,\psi,\zeta)\|^2_{H^1}+\iint\f{|\wt{\mb{G}}_1|^2}{\mb{M}_{_{\#}}}d\xi dy+\iint\f{|\wt{\mb{G}}_y |^2}{\mb{M}_{_\#}}d\xi dy\\[4mm]
				&\di~~~~~~~~~~
				+\iint\f{|\wt{\mb{G}}_t |^2}{\mb{M}_{_\#}}d\xi dy+\iint\f{|\tilde f_{yy}|^2}{\mb{M}_{_\#}}d\xi dy+\iint\f{|\tilde f_{yt}|^2}{\mb{M}_{_\#}}d\xi dy\Bigg\}\leqslant \ve_1^2,
			\end{array}
		\end{equation}
		where the small positive constant $\ve_1$ and the global Maxwellian $\mb{M}_{_\#}$ are to be chosen later. Note that the a priori assumption (\ref{prioriassumption}) implies that
		\begin{equation*}
			\|(\phi,\psi,\z)\|^2_{L^\i}\leqslant
			C\ve_1^2.
		\end{equation*}
		From \eqref{perdef} and \eqref{prioriassumption}, we have
		\begin{equation*}
			\begin{array}{ll}
				\|(v_y,u_y,\t_y)\|_{L^2}^2\leqslant \|(\phi_y,\psi_y,\z_y)\|_{L^2}^2+\|(\bar v_y,\bar u_y,\bar \t_y)\|_{L^2}^2\leqslant C(\d_1+\ve_1)^2.
			\end{array}
		\end{equation*}
		
		Noting that $(\phi,\psi,\z)$ also satisfies
		\begin{equation}\label{(4.18)}
			\left\{
			\begin{array}{ll}
				\di \phi_t-\s \phi_y-\dot{\mb X}(t)(v^S)^{-\mb X}_y-\psi_{1y}=0,\\[1mm]
				\di \psi_{1t}-\s \psi_{1y}-\dot{\mb X}(t)(u_1^S)^{-\mb X}_y+(p- p^R-p^c-(p^S)^{-\mb X})_y=-(u^C_{1t}-\s u^C_{1y})-\int \xi_1^2\wt{\mb G}_y d\xi,\\[1mm]
				\di \psi_{it}-\s \psi_{iy}=-\int \xi_1\xi_i\wt{\mb G}_y d\xi,~~~i=2,3,\\[2mm]
				\di \zeta_t-\s\zeta_y -\dot{\mb X}(t)(\theta^S)^{-\mb X}_y+(pu_{1y}- p^Ru^R_{1y}- p^Cu^C_{1y}- (p^S)^{-\mb X}(u_1^S)^{-\mb X}_{y} )\\[2mm]
				\qquad\qquad\qquad\di =-\left(\frac{\kappa(\t^C)\theta^C_y}{v^C}\right)_y-\int \xi_1\frac{|\xi|^2}{2}\wt{\mb G}_y d\xi+(u_1-(u_1^S)^{-\mb X})\int \xi_1^2(\wt{\mb G}^S)^{-\mb X}_y d\xi\\
				\qquad\qquad\qquad\qquad\di+u_1\int \xi_1^2\wt{\mb G}_y d\xi+\sum_{i=2}^3\int \xi_1\xi_i\wt{\mb G}_y d\xi,
			\end{array}
			\right.
		\end{equation}
		we have
		\begin{equation*}
			\|(\p_{t},\psi_{t},\z_{t})\|_{L^2}^2\leqslant C(\d_1+\ve_1)^2
		\end{equation*}
		and
		\begin{equation*}
			\begin{array}{ll}
				\|(v_{t},u_{t},\t_{t})\|_{L^2}^2\leqslant \|(\phi_{t},\psi_{t},\z_{t})\|_{L^2}^2+ C \|(\bar v_t,\bar u_t,\bar \t_t)\|_{L^2}^2\leqslant C(\d_1+\ve_1)^2,
			\end{array}
		\end{equation*}
		where we have used
		$$\left(\int \xi_1^2\wt{\mb G}_y d\xi \right)^2\leqslant C\int \frac{|\wt{\mb G}_y|^2}{\mb M_{_\#}}d\xi$$
		and
		\begin{equation}\label{xprop}
			|\dot{\mb{X}}(t)|\leqslant C\|(\phi ,\psi_1,\z)\|_{L^\infty}\leqslant C\ve_1,\qquad t\in [0,T].
		\end{equation}
		For the two-order derivatives,
		\begin{equation*}
			\begin{array}{ll}
				\di\iint\f{|f_{yy}|^2}{\mb{M}_{_\#}}d\xi dy
				&\di \leqslant C\iint\f{| \tilde f_{yy}|^2}{\mb{M}_{_\#}} d\xi dy
				+C\iint\f{|(\mb{M}^{S})_{yy}^{-\mb{X}}|^2+|(\mb{G}^{S})_{yy}^{-\mb{X}}|^2}{\mb{M}_{_\#}}d \xi dy\\[4mm]
				&\di \leqslant C(\d_1+\ve_1)^2.
			\end{array}
		\end{equation*}
		Hence we have
		\begin{equation}\label{vandf}
			\begin{array}{ll}
				\di&\quad\|(v_{yy},u_{yy},\t_{yy})\|_{L^2}^2\\[3mm]
				&\di\leqslant C\left\| \partial_{yy}\left(\rho,\rho
				u,\rho\left(\t+\f{|u|^2}{2}\right)\right)\right\| _{L^2}^2+C\int\left| \partial_{y}\left(\rho,\rho
				u,\rho\left(\t+\f{|u|^2}{2}\right)\right)\right| ^4dy\\[4mm]
				&\di 
				\leqslant C\iint\f{| f_{yy}|^2}{\mb{M}_{_\#}}d\xi dy+C\int \left| (v_{y},u_{y},\t_{y})\right|^4  dy \\[4mm]
				&\di \leqslant C(\delta_1+\ve_1)^2
			\end{array}
		\end{equation}
		and
		\begin{equation*}
			\|(\phi_{yy},\psi_{yy},\z_{yy})\|_{L^2}^2\leqslant \|(v_{yy},u_{yy},\t_{yy})\|_{L^2}^2+\|(\bar v_{yy},\bar u_{yy},\bar\t_{yy})\|_{L^2}^2\leqslant
			C(\d_1+\ve_1)^2,
		\end{equation*}
		which together with the Sobolev inequality yields that
		\begin{equation*}
			\|(\phi_{y},\psi_{y},\z_{y})\|^2_{L^{\i}}\leqslant C(\d_1+\ve_1)^2.
		\end{equation*}
		For the microscopic variable, we have
		\begin{equation*}
			\begin{array}{ll}
				\di\left\|\int\f{|\wt{\mb{G}}_1|^2}{\mb{M}_{_\#}}d \xi\right\|_{L^\i_y} \leqslant
				C\left(\iint\f{|\wt{\mb{G}}_1|^2}{\mb{M}_{_\#}}d \xi
				dy\right)^{\f{1}{2}}\cdot\left(\iint\f{|\wt{\mb{G}}
					_{1y}|^2}{\mb{M}_{_\#}}d \xi dy\right)^{\f{1}{2}}\leqslant
				C\ve_1^2.
			\end{array}
		\end{equation*}
		Furthermore, by noticing the facts that $f=\mb{M}+\mb{G}$ and
		$F^{S}=\mb{M}^{S}+\mb{G}^{S}$, it holds that
		\begin{equation*}
			\begin{array}{ll}
				\di\iint\f{|\wt{\mb{G}}_{yy}|^2}{\mb{M}_{_\#}}d\xi dy&\di\leqslant
				C\iint\f{|\wt f_{yy}|^2}{\mb{M}_{_\#}}d \xi dy+C
				\iint\f{|\partial_{yy}(\mb{M}-(\mb{M}^{S})^{-\mb{X}})|^2}{\mb{M}_{_\#}}d \xi dy\\[3mm]
				&\leqslant C(\d_1+\ve_1)^2,
			\end{array}
		\end{equation*}
		which together with the Sobolev embedding gives
		\begin{equation*}
			\begin{array}{ll}
				\di\left\|\int\f{|
					\wt{\mb{G}}_y|^2}{\mb{M}_{_\#}}d \xi\right\|_{L_y^{\i}}&\di\leqslant
				C\left(\iint\f{| \wt{\mb{G}}_y
					|^2}{\mb{M}_{_\#}}d \xi
				dy\right)^{\f{1}{2}}\cdot\left(\iint\f{|
					\wt{\mb{G}}
					_{yy}|^2}{\mb{M}_{_\#}}d \xi dy\right)^{\f{1}{2}}\\
				&\di\leqslant C(\d_1+\ve_1)^2.
			\end{array}
		\end{equation*}
		Finally, by \eqref{xprop}, we also have similar conclusion when taking the derivative of $t$,
		\begin{equation*}
			\iint\f{|f_{yt}|^2}{\mb{M}_{_\#}}d\xi dy\leqslant
			C(\d_1+\ve_1)^2,
		\end{equation*}
		\begin{equation*}
			\|(\phi_{yt},\psi_{yt},\z_{yt})\|_{L^2}^2\leqslant \|(v_{yt},u_{yt},\t_{yt})\|_{L^2}^2+\|(\bar v_{yt},\bar u_{yt},\bar\t_{yt})\|_{L^2}^2\leqslant
			C(\d_1+\ve_1)^2,
		\end{equation*}
		\begin{equation*}
			\|(\phi_{t},\psi_{t},\z_{t})\|^2_{L^{\i}}\leqslant C(\d_1+\ve_1)^2.
		\end{equation*}
		\begin{equation*}
			\begin{array}{ll}
				\di\iint\f{|\wt{\mb{G}}_{yt}|^2}{\mb{M}_{_\#}}d\xi dy\leqslant C(\d_1+\ve_1)^2,
			\end{array}
		\end{equation*}
		\begin{equation*}
			\begin{array}{ll}
				\di\left\|\int\f{|
					\wt{\mb{G}}_t|^2}{\mb{M}_{_\#}}d \xi\right\|_{L_y^{\i}}\leqslant C(\d_1+\ve_1)^2.
			\end{array}
		\end{equation*}
		
		To close the a priori assumption \eqref{prioriassumption} and to prove Theorem
		\ref{maintheorem}, we need  the following a priori estimates.
		
		\begin{proposition} \label{main proposition} (A priori estimates)
			For each $(v_+,u_+,\theta_+)$ with $v_+,\t_+>0$, let $(v_-,u_-,\t_-)$ satisfies that the Riemann problem of Euler equations \eqref{(1.16)} consists of a 1-rarefaction wave, a 2-contact discontinuity, and a 3-shock wave. Suppose that $f(t,y,\xi)$ is a solution to \eqref{tran-Lag-B} on $t\in[0,T]$. Then there exist positive constants $C_0,\delta_1,\varepsilon_1$ $(\d_1,\ve_1<1)$ and a global Maxwellian $\mb M_{_\#}=\mb M_{[v_{\#},u_{\#},\t_{\#}]}$ $(v_{\#},\t_{\#}>0)$ independent of time $T$ such that if the wave strength $\delta_R+\d_C+\d_S\leqslant\delta_1$ and $\di\mathcal{N}(T)\leqslant\ve_1$, then it holds that
			\begin{equation}\label{summary}
				\begin{array}{ll}
					\di \mathcal{N}(T)^2+\delta_S\int_0^T| \dot{\mb{X}}(t)| ^2dt+\int_0^T\left( \Lambda^R+\Lambda^S
					\right) dt+\sum_{|\b|=1}\int_0^T\|\partial^\b(\p,\psi,\z)\|_{L^2}^2dt\\[3mm]
					\di\quad
					+ \int_0^T\left( \|(\phi_{yy},\psi_{yy},\z_{yy})\|_{L^2}^2+\|(\phi_{yt},\psi_{yt},\z_{yt})\|_{L^2}^2\right) dt\\[3mm]
					\di\quad
					+ \int_0^T\iint\f{(1+|\x|) \left(|\wt{\mb{G}}_{1}|^2+|\wt{\mb{G}}_{y}|^2+|\wt{\mb{G}}_{t}|^2+ |\wt{\mb{G}}_{yy}|^2+|\wt{\mb{G}}_{yt}|^2\right) }{\mb{M}_{_\#}}d \xi dydt\\[4mm]
					\di\leqslant C_0\left( \mathcal{N}(0)^2+\d_1^{\f12}\right),
				\end{array}
			\end{equation}
			where $\partial^\b$ denotes the derivatives with respect to $y$ or $t$, and
			\begin{align*}
				\begin{aligned}
					& \Lambda^R:= \int_{-\i}^{+\i}|v^R_y| |(\phi, \z)|^2 dy ,\\
					&\Lambda^S:=\int_{-\i}^{+\i}|(v^S)^{-\mb X}_y| |(\phi, \psi, \z)|^2 dy.
				\end{aligned}
			\end{align*}
		\end{proposition}
		
		The Proof of Theorem \ref{main proposition} will be given in Section 5 and Section 6.

		\subsection{Local-in-time existence}
		
		In this subsection, we give the local-in-time existence of solution to the Cauchy problem of 1D Boltzmann equation in Lagrangian coordinate \eqref{tran-Lag-B}. To state the local-in-time existence precisely, we don't need the shift $\mb X(t)$ in the ansatz \eqref{ansatz} and \eqref{barM}. In order to highlight this difference, we define the following notation
		\begin{equation}\label{hat-v}
			\begin{array}{l}
				\left(\begin{array}{cc} \hat v\\[1mm] \hat u_1 \\[1mm] \hat u_2 \\[1mm] \hat u_3 \\[1mm] \hat \t
				\end{array}
				\right)(t,y):= \left(\begin{array}{cc}v^{R}(t,y+\s t)+ v^{C}(y+\s t)+ v^{S}(y)-v_*-v^*\\ [1mm] u_1^{R}(t,y+\s t)+ u_1^{C}(t,y+\s t)+ u_1^{S}(y)-u_{1*}-u_1^* \\ 0 \\ 0 \\
					\t^{R}(t,y+\s t)+ \t^{C}(y+\s t)+ \t^{S}(y)-\t_*-\t^*
				\end{array}
				\right),
			\end{array}
		\end{equation}
		and
		\begin{equation}\label{hat-M}
			\hat{\mb M}(t,y,\xi):=\mb M_{[\hat v,\hat u,\hat \t](t,y)}(\xi).
		\end{equation}
		Note that $(\hat v,\hat u,\hat \t,\hat{\mb M})=(\bar v,\bar u,\bar\t,\bar {\mb M})$ when $\mb X(t)=0$. It is easy to check that
		$$\|\partial^\b(\hat v,\hat u,\hat \t)\|_{L^2}^2\leqslant C(\d_R+\d_C+\d_S),~~1\leqslant |\b|\leqslant3,$$
		and
		$$\iint\frac{|\partial^\b\hat{\mb M}(t,y,\xi)|^2}{\mb M_{_\#}(\xi)}d\xi dy\leqslant C(\d_R+\d_C+\d_S),~~1\leqslant |\b|\leqslant3.$$
		The local-in-time existence of solution to Boltzmann equation \eqref{tran-Lag-B} can be stated as follows, whose proof will be given in Appendix A.2.
		\begin{proposition} \label{localexistence}
			(Local-in-time existence) There exist two independent positive constants $\d_2$, $\ve_2$ such that if the wave strength satisfies $\delta_R+\d_C+\d_S<\d_2$ and the initial data satisfies $f_0(y,\xi)\geqslant0$ and
			\begin{equation}\label{epsilon2}
				\begin{array}{ll}
					\di \left\| f_0(y,\xi)-\hat{\mb M}(0,y,\xi)\right\|_{H^2_y\left(L^2_\xi\left(\frac{1}{\sqrt{\mb M_{_\#}}} \right)  \right) }\leqslant \chi<\ve_2,
				\end{array}
			\end{equation}
			for some positive constant $\chi$,
			then there exists a positive constant $T_0=T_0(\chi)$ such that the Cauchy problem \eqref{tran-Lag-B} admits a unique solution $f(t,y,\xi)$ on $[0,T_0]\times \mathbb R\times \mathbb R^3$ satisfying $f(t,y,\xi)\geqslant0$ and
			\begin{equation*}
				\sup_{t\in [0,T_0]} \left\| f(t,y,\xi)-\hat{\mb M}(t,y,\xi)\right\|_{H^2_y\left(L^2_\xi\left(\frac{1}{\sqrt{\mb M_{_\#}}} \right)  \right) }\leqslant 2\chi.
			\end{equation*}
		\end{proposition}

		\subsection{Global-in-time existence}
		
		Before we construct the global-in-time solution of \eqref{tran-Lag-B}, we first claim that there exist positive constants $\ve_3<\ve_1$ and $\delta_3<\min\left\lbrace \d_1,\d_2 \right\rbrace $ independent of time $T$ such that if $\mathcal{N}(T)\leqslant\ve_3$ and $\d_R+\d_C+\d_S\leqslant\d_3$, then the solution of \eqref{tran-Lag-B} can be further extended  by one more small step in time. In fact, 
		$$\|f(T,y,\xi)-\mb M_{[\bar v,\bar u,\bar \theta]}(T,y,\xi)\|_{H^2_y\left(L^2_\xi\left(\frac{1}{\sqrt{\mb M_{_\#}}} \right)  \right) }\leqslant C(\mathcal{N}(T)+\d_R+\d_C+\d_S).$$
		Thus we can extend the solution in time by Proposition \ref{localexistence} as long as $\mathcal{N}(T)\leqslant\ve_3$ and $\d_R+\d_C+\d_S\leqslant\d_3$ are taken suitably small.
		
		\begin{remark}
			In order to get the global-in-time existence, we need to repeatedly apply the local-in-time existence (Proposition \ref{localexistence}). However, Proposition \ref{localexistence} is based on the smallness of both the perturbation and wave strength, therefore the claim is crucial in the extension of the solution. 
		\end{remark}
		
		We will use the following two steps to obtain the global-in-time existence in Theorem \ref{maintheorem}.
		
		{\it \underline{Step 1}. (Local-in-time existence)} We will construct the local-in-time solution based on the initial value \eqref{epsilon0} by Proposition \ref{localexistence}. By the construction of the ansatz, we have
		\begin{equation*}
			\begin{array}{ll}
				\di \sum_\pm\left\| \hat{\mb M}(0,y,\xi)-\mb M_{[v_\pm,u_\pm,\t_\pm]}(\xi) \right\|_{L^2_y\left(\mathbb R_\pm,L^2_\xi\left(\frac{1}{\sqrt{\mb M_{_\#}}} \right)  \right) }\\[6mm]
				\di \quad +\left\| \hat{\mb M}_{y}(0,y,\xi) \right\|_{H^1_y\left(L^2_\xi\left(\frac{1}{\sqrt{\mb M_{_\#}}} \right)  \right) }\leqslant C(\d_R+\d_C+\d_S),
			\end{array}
		\end{equation*}
		which and  \eqref{epsilon0} give \eqref{epsilon2} with $\chi=C(\ve_0+\d_R+\d_C+\d_S)<\ve_2$ provided that $\ve_0,\d_R,\d_C,\d_S$ are chosen suitably small. By Proposition \ref{localexistence}, the Boltzmann equation \eqref{tran-Lag-B} admits a unique solution $f(t,y,\xi)$ on $[0,T_0]\times \mathbb R\times \mathbb R^3$ satisfying
		\begin{equation}\label{(3.54)}
			\sup_{t\in [0,T_0]} \left\| f(t,y,\xi)-\hat{\mb M}(t,y,\xi)\right\|_{H^2_y\left(L^2_\xi\left(\f1{\sqrt{\mb M_{_\#}}} \right) \right)  }\leqslant 2C(\ve_0+\d_R+\d_C+\d_S).
		\end{equation}
		
		Next we need to verify the a priori assumption \eqref{prioriassumption}. It can be easily seen from \eqref{(3.54)} that for any $t\in [0,T_0]$,
		\begin{equation}\label{assumption1}
			\begin{array}{ll}
				\di&\di
				\|(v-\hat v,u-\hat u,\t-\hat\t)\|^2_{H^1}+\iint\f{|\mb{G}|^2+|\mb{G}_y|^2+|f_{yy}|^2}{\mb{M}_{_\#}}d\xi dy\\[4mm]
				\leqslant&\di C\left\| f(t,y,\xi)-\hat{\mb M}(t,y,\xi)\right\|^2_{H^2_y\left(L^2_\xi\left(\f1{\sqrt{\mb M_{_\#}}} \right) \right)  }+C(\d_R+\d_C+\d_S)^2\\[4mm]
				\leqslant&\di C\ve_0^2+C(\d_R+\d_C+\d_S)^2.
			\end{array}
		\end{equation}
		Using the arguments similar to the energy estimates, one can get
		\begin{equation*}
			\begin{array}{ll}
				\di\di
				\iint\f{|\mb{G}_t |^2+|f_{yt}|^2}{\mb{M}_{_\#}}d\xi dy \leqslant\di C\ve_0^2+C(\d_R+\d_C+\d_S)^2,
			\end{array}
		\end{equation*}
		which together with \eqref{roughx} gives
		\begin{equation}\label{assumption2}
			\begin{array}{ll}
				\di\di
				\iint\f{|\wt{\mb{G}}_1 |^2+|\wt{\mb{G}}_y |^2+|\wt{\mb{G}}_t |^2+|\wt f_{yt}|^2+|\wt f_{yy}|^2}{\mb{M}_{_\#}}d\xi dy \leqslant\di C\ve_0^2+C(\d_R+\d_C+\d_S)^2,
			\end{array}
		\end{equation}
		Since the only difference between $(\hat v,\hat u,\hat \t)$ and $(\bar v,\bar u,\bar\t)$ is the construction of shift $\mb X(t)$, we have
		\begin{equation}\label{assumption3}
			\begin{array}{ll}
				\di
				\|(\hat v-\bar v,\hat u-\bar u,\hat\t-\bar\t)\|^2_{H^1}
				&\di=\|(v^S-(v^S)^{-\mb X},u^S-(u^S)^{-\mb X},\t^S-(\t^S)^{-\mb X})\|^2_{H^1}\\[4mm]
				&\leqslant\di C\d_S(1+t),
			\end{array}
		\end{equation}
		where in the last inequality we have used
		\begin{equation*}
			\begin{array}{ll}
				\di&\di
				\|v^S-(v^S)^{-\mb X}\|^2_{H^1}\\[4mm]
				\leqslant&\di C \|v^S-v_+\|^2_{L^2(\mathbb R_+)}+C \|v^S-v^* \|^2_{L^2(\mathbb R_-)}+C \|(v^S)^{-\mb X}-v_+\|^2_{L^2(\mathbb R_+)}+C \|(v^S)^{-\mb X}-v^*\|^2_{L^2(\mathbb R_-)}\\[4mm]
				&~~~~\di+C\|v^S_y\|^2_{L^2}+C\|(v^S)^{-\mb X}_y\|^2_{L^2}\\[4mm]
				\leqslant&\di C\d_S+\int_{-\mb X(t)}^{+\i}(v^S-v_+)^2 dy+\int^{-\mb X(t)}_{-\i}(v^S-v^*)^2 dy\\[4mm]
				\leqslant&\di C\d_S+\int_{-\mb X(t)}^{0}(v^S-v_+)^2 dy+\int_{0}^{+\i} C\d_S^2 e^{-c\d_Sy}dy+\int^{-\mb X(t)}_{0}(v^S-v^*)^2 dy+\int^{0}_{-\i} C\d_S^2 e^{c\d_Sy}dy\\[4mm]
				\leqslant&\di C\d_S(1+t).\\[4mm]
			\end{array}
		\end{equation*}
		Therefore, the combination of \eqref{assumption1}-\eqref{assumption3} gives
		$$\mathcal{N}(t)^2\leqslant C_1\ve_0^2+C_1(\d_R+\d_C+\d_S)(1+t),$$
		for some positive constant $C_1$. Thus the a priori assumption $\mathcal{N}(T_0)^2\leqslant\ve_3^2$ holds provided $\ve_0,\d_R,\d_C,\d_S$ are chosen smaller, that is, $\ve_0<\ve_4$ and $\d_R+\d_C+\d_S<\d_4$ for some $\ve_4\in(0,\ve_3),\d_4\in(0,\d_3)$.

		{\it \underline{Step 2}.} We now consider the maximal existence time
		$$T_{M}:=\sup\left\lbrace t>0|~\mathcal{N}(t)\leqslant\ve_3 \right\rbrace.$$
		If $T_{M}<+\i$, then the continuity argument gives 
		\begin{equation}\label{v31}
			\mathcal{N}(T_{M})=\ve_3.
		\end{equation}
		However, it holds from Proposition \ref{main proposition} that
		\begin{equation}\label{v32}
			\begin{array}{ll}
				\di\mathcal{N}(T_{M})^2&\di\leqslant C_0 (\mathcal{N}(0)^2+\d_0^{\f12})\\[3mm]
				&\leqslant\di C_0(C_1\ve_0^2+C_1\d_0+\d_0^{\f12}) \\[3mm]
				&\leqslant\di \frac{\ve_3^2}{4}+\frac{\ve_3^2}{4} =\di \frac{\ve_3^2}{2}
			\end{array}
		\end{equation}
		if we take $\d_R+\d_C+\d_S<\d_0:=\min\left\lbrace \d_4, \frac{\ve_3^4}{64(C_0^2+1)(C_1+1)}\right\rbrace $ and $\ve_0:=\min \left\lbrace \ve_4,\frac{\ve_3}{2\sqrt{C_0C_1}}\right\rbrace$. Therefore the contradiction between \eqref{v31} and \eqref{v32} yields $T_M=+\i$, which completes the proof of the global-in-time existence.

		\subsection{Time-asymptotic behaviors}
		
		In this subsection, we prove the time-asymptotic behaviors and \eqref{mainlimit} and \eqref{Xdotlimit}. By \eqref{summary}, it holds that
		\begin{equation}\label{bounded}
			\begin{array}{ll}
				\di \iint\frac{|f-\bar f|^2}{\mb{M}_{_\#}}d\xi dy&\di \leqslant 2\iint\frac{|\mb{M}-\bar{\mb M}|^2}{\mb{M}_{_\#}}d\xi dy+2\iint\frac{|\wt{\mb G}|^2}{\mb{M}_{_\#}}d\xi dy\\[4mm]
				&\di \leqslant C\left\|(\phi,\psi,\z) \right\| _{L^2}^2+C\iint\frac{|\wt{\mb G}_0|^2}{\mb{M}_{_\#}}d\xi dy+C\iint\frac{|\wt{\mb G}_1|^2}{\mb{M}_{_\#}}d\xi dy\\[4mm]
				&\di \leqslant  C\left( \mathcal{N}(0)^2+\d_1^{\f12}\right),
			\end{array}
		\end{equation}
		and
		\begin{equation}\label{bounded2}
			\begin{array}{ll}
				\di \int_0^{+\i}\iint\frac{|(f-\bar f)_y|^2}{\mb{M}_{_\#}}d\xi dydt
				&\di \leqslant C\int_0^{+\i}\iint\left( \frac{|(\mb M-\bar{\mb M})_y|^2}{\mb{M}_{_\#}}+\frac{|\wt{\mb G}_y|^2}{\mb{M}_{_\#}}\right)d\xi  dydt\\[4mm]
				&\di \leqslant C\left( \mathcal{N}(0)^2+\d_1^{\f12}\right).
			\end{array}
		\end{equation}
		We also have
		\begin{equation}\label{derivativebounded}
			\begin{array}{ll}
				&\di \int_0^{+\i}\left| \frac{d}{dt}\iint\frac{|(f-\bar f)_y|^2}{\mb{M}_{_\#}}d\xi dy\right| dt\\[4mm]
				\di \leqslant&\di\int_0^{+\i}\iint\frac{|(f-\bar f)_y|^2}{\mb{M}_{_\#}}d\xi  dydt+\int_0^{+\i}\iint \frac{|(f-\bar f)_{yt}|^2}{\mb{M}_{_\#}}d\xi  dydt\\[4mm]
				\di \leqslant&\di \int_0^{+\i}\iint \frac{|(f-\bar f)_y|^2}{\mb{M}_{_\#}}d\xi  dydt+2\int_0^{+\i}\iint\left(  \frac{|(\mb M-\bar{\mb M})_{yt}|^2}{\mb{M}_{_\#}}+\frac{|\wt{\mb G}_{yt}|^2}{\mb{M}_{_\#}}\right) d\xi  dydt\\[4mm]
				\di \leqslant&\di C\left( \mathcal{N}(0)^2+\d_1^{\f12}\right).
			\end{array}
		\end{equation}
		\eqref{bounded2} and \eqref{derivativebounded} yield that
		$$\lim\limits_{t\to +\i}\iint\frac{|(f-\bar f)_y|^2}{\mb{M}_{_\#}}d\xi dy=0,$$
		which as well as \eqref{bounded} and the Sobolev inequality
		$$
		\di \left\|\int\frac{|f-\bar f|^2}{\mb{M}_{_\#}} d \xi\right\|^2_{L^{\i}_y}
		\leqslant C\left(\iint\frac{|f-\bar f|^2}{\mb{M}_{_\#}} d\xi dy\right)\cdot\left(\iint\frac{|(f-\bar f)_y|^2}{\mb{M}_{_\#}} d \xi dy\right)
		$$
		easily leads to
		\begin{equation}\label{fbarf}
			\lim\limits_{t\to +\i}\left\|f-\bar f\right\|_{L^\i_y\left(L^2_\xi\left(\frac{1}{\sqrt{\mb M_{_\#}}} \right)  \right) }=0.
		\end{equation}
		Finally, observing that
		\begin{equation*}
			\begin{array}{ll}
				\di \lim\limits_{t\to +\infty}\bigg\| \bar f(t,y,\xi)-\big[\mb M_{\left[ v^R\left(t,y \right) ,u^R\left(t,y\right),\t^R\left(t,y \right)\right] }(\xi)+ \mb M_{\left[ v^C\left(\frac{y+\s t}{\sqrt{1+t}} \right) ,u^C\left(t,\frac{y+\s t}{\sqrt{1+t}} \right),\t^C\left(\frac{y+\s t}{\sqrt{1+t}} \right)\right] }(\xi) \\[5mm]
				~~~~~~~~~~~+F^S(y-\mb X(t),\xi)-\mb M_{[v_*,u_*,\t_*]}(\xi)-\mb M_{[v^*,u^*,\t^*]}(\xi)\big] \bigg\|_{L^\i_y\left(L^2_\xi\left(\frac{1}{\sqrt{\mb M_{_\#}}} \right)  \right) } =0,
			\end{array}
		\end{equation*}
		we can prove	\eqref{mainlimit} by combining with \eqref{Randr}. In addition, by \eqref{fbarf}, we have
		\begin{equation*}
			\lim\limits_{t\to +\infty}\|(\phi,\psi,\z)(t, \cdot)\|_{L^\i}=0,
		\end{equation*}
		which together with \eqref{xprop} implies  \eqref{Xdotlimit}. Therefore, we complete the proof of Theorem \ref{maintheorem}.
		
		\
		
		\section{Properties for the collision operator}
		\setcounter{equation}{0}
		In this section, we list some basic lemmas based on the celebrated H-Theorem for later use. 
		First we give the properties for the gain part $Q_+(g,h)$ and loss part $Q_-(g,h)$, whose proof can be found in \cite{Golse-Perthame-Sulem,Metivier-Zumbrun}.
		
		\begin{lemma}\label{Lemma 4.1} There exists a positive constant $C$ such that for the loss part $Q_-(g,h)$,
			$$
			\int\f{(1+|\xi|)^{-1}Q_-^2(g,h)}{\wt{\mb{M}}}d\xi\leqslant
			C\int\f{(1+|\xi|)g^2}{\wt{\mb{M}}}d\xi\cdot\int\f{h^2}{\wt{\mb{M}}}d\xi,
			$$
			and for the gain part $Q_+(g,h)$,
			$$
			\int\f{(1+|\xi|)^{-1}Q_+^2(g,h)}{\wt{\mb{M}}}d\xi\leqslant
			C\int\f{g^2}{\wt{\mb{M}}}d\xi\cdot\int\f{(1+|\xi|)h^2}{\wt{\mb{M}}}d\xi.
			$$
			As a result,
			$$
			\int\f{(1+|\xi|)^{-1}Q^2(g,h)}{\wt{\mb{M}}}d\xi\leqslant
			C\left[ \int\f{(1+|\xi|)g^2}{\wt{\mb{M}}}d\xi\cdot\int\f{h^2}{\wt{\mb{M}}}d\xi+
			\int\f{g^2}{\wt{\mb{M}}}d\xi\cdot\int\f{(1+|\xi|)h^2}{\wt{\mb{M}}}d\xi\right] .
			$$
			where $\wt{\mb{M}}$ can be any Maxwellian so that the above
			integrals are well-defined.
		\end{lemma}
		
		Now we recall some basic properties of the linearized collision operator $\mb L_{\mb M}$. For the hard sphere model, $\mb L_{\mb M}$ takes the form, cf. \cite{Grad},
		\begin{equation}\label{LM}
			(\mb L_\mathbf{M}h) (\xi)=-\nu_\mathbf{M}(\xi)h(\xi)+\sqrt{\mathbf{M}(\xi)}\mb K_\mathbf{M}\left(\left(\frac{h}{\sqrt{\mathbf{M}}}\right)(\xi)\right).
		\end{equation}
		Here $\di \nu_\mathbf{M}(\xi)=\int_{{\mathbb R}^3}\!\!\int_{{\mathbb S}_+^2} |(\xi-\xi_*)\cdot \Omega|\mb M(\xi_*)d\Omega d\xi_*$ and  $\mb K_{\mathbf{M}}(\cdot) = -\mb K_{1\mathbf{M}}(\cdot) +\mb K_{2\mathbf{M}}(\cdot)$ is a compact $L^2$-operator. The collision frequency $\nu_{\mb M}$ and $\mb K_{i\mathbf{M}}$ have the following expressions
		\begin{equation}\label{4.2}
			\begin{aligned}\nu_{\mathbf{M}}(\xi)&=\frac{2\rho}{\sqrt{2\pi R\theta}}\left\{\left(\frac{R\theta}{|\xi-u|}+|\xi-u|\right)\int_0^{|\xi-u|}\exp\left(-\frac{y^2}{2R\theta}\right)dy+R\theta\exp\left(-\frac{|\xi-u|^2}{2R\theta}\right)\right\},\\k_{1\mathbf{M}}(\xi,\xi_*)&=\frac{\pi\rho}{\sqrt{(2\pi R\theta)^3}}|\xi-\xi_*|\exp\left(-\frac{|\xi-u|^2}{4R\theta}-\frac{|\xi_*-u|^2}{4R\theta}\right),\\k_{2\mathbf{M}}(\xi,\xi_*)&=\frac{2\rho}{\sqrt{2\pi R\theta}}|\xi-\xi_*|^{-1}\exp\left(-\frac{|\xi-\xi_*|^2}{8R\theta}-\frac{(|\xi|^2-|\xi_*|^2)^2}{8R\theta|\xi-\xi_*|^2}\right),\end{aligned}
		\end{equation}
		where $k_{i\mb M}(\xi,\xi_*)(i=1,2)$ is the kernel of the operator $\mb K_{i\mb M}$, respectively, and $\nu_{\mb M}(\xi)\sim(1+|\xi|)$. Furthermore, we have the following properties for the operator $\mb K_{\mathbf{M}}$ and $\mb L_{\mathbf{M}}$, whose proof can be found in \cite{Liu-Yang-Yu-Zhao,Ukai-Yang-Zhao}.
		
		\begin{lemma}\label{Property of K} If $\di \t/2<\t_{_\#}<\t$, then the operator $\mb K_{\mathbf{M}}$ satisfies that
			$$\int\frac{\left|\sqrt{{\mathbf{M}}}\mb K_{\mathbf{M}}\left(\frac h{\sqrt{\mathbf{M}}}\right)\right|^2}{\mathbf{M}_{_\#}}d\xi\leqslant C \int\frac{h^2}{\mathbf{M}_{_\#}}d\xi.$$
			
		\end{lemma}
		
		\begin{lemma}\label{Lemma 4.3} If $\di \t/2<\t_{_\#}<\t$, then there exist two
			positive constants $\wt\sigma=\wt\sigma(v,u,\t;
			v_{_\#},u_{_\#},\t_{_\#})$ and
			$\eta_0=\eta_0(v,u,\t;v_{_\#},u_{_\#},\t_{_\#})$ such that if
			$|v-v_\#|+|u-u_\#|+|\t-\t_\#|<\eta_0$, we have for
			$g(\xi)\in  \mathfrak{N}^\bot$,
			$$
			-\int\f{g\mb{L}_\mb{M}g}{\mb{M}_{_\#}}d\xi\geqslant
			\wt\sigma\int\f{(1+|\xi|)g^2}{\mb{M}_{_\#}}d\xi,
			$$
			and, by Cauchy inequality,  for each $g(\xi)\in  \mathfrak{N}^\bot$,
			$$\int\f{1+|\xi|}{\mb{M}_{_\#}}|\mb{L}_\mb{M}^{-1}g|^2d\xi\leqslant
			\wt\sigma^{-2}\int\f{(1+|\xi|)^{-1}g^2}{\mb{M}_{_\#}}d\xi.
			$$
		\end{lemma}

		\
		
		\section{Lower order estimates}
		\setcounter{equation}{0}
		The remaining part of the paper is dedicated to the proof of Proposition \ref{main proposition}. The lower order estimates to the system \eqref{pereq1} and \eqref{G1eq}
		are given in the following proposition.
		
		\begin{proposition}\label{lowerorderestimate}
			Under the assumptions of Proposition \ref{main proposition}, there exists a positive
			constant C such that
			\begin{equation*}
				\begin{array}{ll}
					\di \sup_{t\in[0,T]}\left[ \left\|(\phi,\psi,\z,\phi_y) \right\| _{L^2}^2+\iint\f{|\wt{\mb{G}}_1|^2}{\mb{M}_{_\#}}d \xi dy\right] +\delta_S\int_0^T| \dot{\mb{X}}(t)| ^2dt+\int_0^T\left( \Lambda^R+\Lambda^S\right) dt\\[3mm]
					\di\qquad\qquad+\sum_{|\b|=1}\int_0^T\|\partial^\b(\p,\psi,\z)\|_{L^2}^2dt
					+\int_0^T\iint\f{(1+|\xi|)|\wt{\mb{G}}_1|^2}{\mb{M}_{_\#}}d \xi dydt\\[3mm]
					\di\quad\leqslant C\mathcal{N}(0)^2+C\ve_1^2\int_0^T\|\psi_{1yy}\|_{L^2}^2 dt+C\d_C\int_0^T\f{1}{1+t} \int e^{-\frac{c_0 |y+\sigma t|^2}{1+t}} |(\phi,\psi,\z)|^2 dydt\\[3mm]
					\di\qquad\qquad
					+C \int_0^T\iint\f{(1+|\xi|)\left( |\wt{\mb{G}}_y|^2+|\wt{\mb{G}}_t|^2+|\wt{\mb{G}}_{yy}|^2+|\wt{\mb{G}}_{yt}|^2\right) }{\mb{M}_{_\#}}d \xi dydt+C\d_1^{\f12}.
				\end{array}
			\end{equation*}
		\end{proposition}

		{\it Proof of Proposition \ref{lowerorderestimate}.} The proof is divided into the following seven steps.
		
		{\it \underline{Step 1}. Construction of weighted relative entropy method.} Let
		$$\Phi(z):=z-1-\ln z,$$
		so $\Phi^\prime(z):=1-\f1z$. Using the same method as \cite{Huang-Wang-Yang-2010-2}, it holds that
		\begin{equation}\label{(4.2)}
			\begin{array}{ll}
				\di\quad 
				\left[\f23\bar\theta\Phi\left(\frac{v}{\bar v}\right)+\bar\t\Phi\left(\frac{\theta}{\bar \theta}\right) + \sum_{i=1}^3 \f12 \psi_i^2 \right]_t-\s\left[\f23\bar\theta\Phi\left(\frac{v}{\bar v}\right)+\bar\t\Phi\left(\frac{\theta}{\bar \theta}\right) + \sum_{i=1}^3 \f12 \psi_i^2 \right]_y \\[5mm]
				\di
				=( \bar\theta_t -\s\bar\theta_y) \left[  \f23\Phi\left(\frac{v}{\bar v}\right) +\Phi\left(\frac{\theta}{\bar \theta}\right) \right]+ \dot{\mb X}(t) \left[ (v^S)^{-\mb{X}}_{y}\frac{\bar p }{\bar v}\phi+(u^S_{1})^{-\mb{X}}_y\psi_1+(\t^S)^{-\mb{X}}_{y}\frac{\z}{\bar\theta} \right]  \\[4mm]
				\di\qquad
				+\left[\f\z\t\left(\kappa(\t)\frac{\t_y}{v}-\kappa(\bar\t)\frac{\bar\t_y}{\bar v}\right) -(p-\bar p)\psi_1+\f43\psi_1\left(\mu(\t)\frac{u_{1y}}{v}-\mu(\bar\t)\frac{\bar u_{1y}}{\bar v}  \right)+\f{\mu(\t)}{v} \sum_{i=2}^3 \psi_i\psi_{iy} \right]_y \\[4mm]
				\di\qquad
				+\left[ -\frac{\bar p\bar u_{1y}}{v\bar v}\phi^2 -\bar u_{1y}\frac{\z(p-\bar p)}{\theta}+\bar p\bar u_{1y}\frac{\z^2}{\t\bar\t}\right]  -\bigg(\f{\z}{\t} \bigg)_y \frac{\t_y}{v} (\k(\t)-\k(\bar \t))-\f43\frac{u_{1y}\psi_{1y}}{v}(\mu(\t)-\mu(\bar \t))\\[4mm]
				\di\qquad
				+\frac{\kappa(\bar\t)}{v\bar v \t}\bar \t_y \z_y \phi+\frac{\kappa(\bar\t)}{\t^2}\z \t_y \left( \frac{\t_y}{v}-\frac{\bar\t_y}{\bar v}\right)+\f\z\t\left[\f43\left(\mu( \t)\frac{u_{1y}^2}{v}-\mu(\bar \t)\frac{\bar u_{1y}^2}{\bar v} \right) +\f{\mu(\t)}{v} \sum_{i=2}^3 \psi_{iy}^2 -Q_2\right]   \\[4mm]
				\di\qquad
				+\frac{4}{3}\frac{\mu(\bar\t)}{v\bar v} \bar u_{1y}\psi_{1y}\phi -Q_1\psi_1-\f{\z^2}{\t\bar \t}\left[\left(\kappa(\bar\t)\frac{\bar\t_y}{\bar v} \right)_y +\f43 \mu(\bar\t)\frac{\bar u_{1y}^2}{\bar v}+Q_2 \right]  \\[4mm]
				\di\qquad
				-\frac{\k(\bar\t)}{v\t}\z_y^2-\f43\frac{\mu(\bar \t)}{v} \psi_{1y}^2-\frac{\mu(\t)}{v}\sum_{i=2}^3\psi_{iy}^2 \\[4mm]
				\di\qquad
				-\psi_1\int \xi_1^2\left(\Pi_1-(\Pi_1^S)^{-\mb X} \right)_y d\xi-\sum_{i=2}^3\psi_i\int \xi_1\xi_i\Pi_{1y} d\xi-\f\z\t\int \xi_1\f{|\xi|^2}2\left(\Pi_1-(\Pi_1^S)^{-\mb X} \right)_y d\xi \\[4mm]
				\di\qquad
				+\f\z\t\sum_{i=2}^3\psi_i\int \xi_1\xi_i\Pi_{1y} d\xi +\f\z\t\left( u_1\int\xi_1^2\Pi_{1y}d\xi-(u_1^S)^{-\mb X}\int\xi_1^2(\Pi_1^S)^{-\mb X}_y d\xi\right) \\[4mm]
				\di\qquad
				+\f{\z^2}{\t\bar \t}\int\xi_1\left(\f12|\xi|^2-(u_1^S)^{-\mb X}\xi_1 \right)(\Pi_1^S)^{-\mb X}_y d\xi.
			\end{array}
		\end{equation}
		Here we will use weighted energy method with the weight function $a^{-\mb X}:=a(y-\mb X(t))$. Simple calculation shows that
		\begin{equation}\label{(4.3)}
			\begin{array}{ll}
				\di\quad
				\left[ a^{-\mb X}\left(\f23\bar\theta\Phi\left(\frac{v}{\bar v}\right)+\bar\t\Phi\left(\frac{\theta}{\bar \theta}\right) + \sum_{i=1}^3 \f12 \psi_i^2 \right)\right]_t-\s\left[ a^{-\mb X}\left(\f23\bar\theta\Phi\left(\frac{v}{\bar v}\right)+\bar\t\Phi\left(\frac{\theta}{\bar \theta}\right) + \sum_{i=1}^3 \f12 \psi_i^2 \right)\right]_y \\[3mm]
				\di
				=a^{-\mb X}\left\lbrace \left[\f23\bar\theta\Phi\left(\frac{v}{\bar v}\right)+\bar\t\Phi\left(\frac{\theta}{\bar \theta}\right) + \sum_{i=1}^3 \f12 \psi_i^2 \right]_t-\s\left[\f23\bar\theta\Phi\left(\frac{v}{\bar v}\right)+\bar\t\Phi\left(\frac{\theta}{\bar \theta}\right) + \sum_{i=1}^3 \f12 \psi_i^2 \right]_y\right\rbrace\\[3mm]
				\di \qquad
				-\dot{\mb X}(t)a^{-\mb X}_y\left(\f23\bar\theta\Phi\left(\frac{v}{\bar v}\right)+\bar\t\Phi\left(\frac{\theta}{\bar \theta}\right) + \sum_{i=1}^3 \f12 \psi_i^2 \right)-\s a^{-\mb X}_y\left(\f23\bar\theta\Phi\left(\frac{v}{\bar v}\right)+\bar\t\Phi\left(\frac{\theta}{\bar \theta}\right) + \sum_{i=1}^3 \f12 \psi_i^2 \right).
			\end{array}
		\end{equation}
		We are going to simplify the first and fourth term of the right-hand side of \eqref{(4.2)} to a more convenient form. Since
		\begin{equation}\label{(4.4)}
			\bar\theta_t -\s\bar\theta_y = (\theta_t^R-\sigma \theta_y^R)+(\theta_t^C-\sigma \theta_y^C) - \dot{\mb X}(t) (\theta^S)^{-\mb X}_y-\sigma (\theta^S)^{-\mb X}_y
		\end{equation}
		and $\bar u_{1y} = u_{1y}^R +u_{1y}^C +(u_{1}^S)^{-\mb X}_y$, we have
		\begin{equation*}
			\begin{array}{ll}
				\quad\di\f23\Big((\theta_t^R-\sigma \theta_y^R)+(\theta_t^C-\sigma \theta_y^C)-\sigma (\theta^S)^{-\mb X}_y \Big)\Phi(\frac{v}{\bar v})-\frac{\bar p\bar u_{1y}}{v\bar v}\phi^2\\[3mm]
				\di=\left[ -u^R_{1y}\left( \f23p^R\Phi(\frac{v}{\bar v})+\frac{\bar p\phi^2}{v\bar v}\right) \right] +\left[ \f23(\theta_t^C-\sigma \theta_{y}^C)\Phi(\frac{v}{\bar v})-\frac{\bar p u^C_{1y}}{v\bar v}\phi^2\right]\\[4mm]
				\di\qquad +\left[ -\f23\sigma (\theta^S)^{-\mb X}_y \Phi(\frac{v}{\bar v})-\frac{(u_1^S)^{-\mb X}_y \bar p }{v\bar v}\phi^2\right] \\[3mm]
				\di =:A_1+A_2+A_3.
			\end{array}
		\end{equation*}
		Using the fact $|\phi|\leqslant C\ve_0$, and by $\Phi(1)=\Phi'(1)=0, \Phi''(1)=1$,
		\begin{equation*}
			\Phi\left( \frac{v}{\bar v}\right)  =  \frac{\phi^2}{2\bar v^2} + O(|\phi|^3),
		\end{equation*}
		and 
		\begin{align*}
			\begin{aligned}
				|\bar p - p^R| & \leqslant C( |\bar v - v^R| +|\bar \theta - \theta^R| ) \\
				&\leqslant C(|v^S - v^*| +|v^C - v_*| +|\theta^S - \theta^*| +|\theta^C - \theta_*| )  \leqslant C(\delta_S +\delta_C), 
			\end{aligned}
		\end{align*}
		we have
		\begin{align*}
			A_{1} & \leqslant -\frac{4\bar p}{3\bar v^2}u^R_{1y} \phi^2 + C(\delta_1 +\ve_1 ) |u^R_{1y}| \phi^2.
		\end{align*}
		Using $\eqref{vcex}_4$, we have
		\[
		A_{2} \leqslant C( |u^C_{1y}|+ |\theta^C_{yy}| + |\theta^C_y||v^C_y| + |Q^C_2| ) \phi^2,
		\]
		which together with \eqref{vc-pe} and \eqref{QC12} yields
		\[
		A_{2} \leqslant C \d_C (1+t)^{-1}e^{-\frac{c_0 |y+\sigma t|^2}{1+t}}\phi^2.
		\]
		Using \eqref{shock-vu}, \eqref{theta-s}, \eqref{sm1} and 
		\begin{align}\label{perror}
			\begin{aligned}
				|\bar p - p^*| & \leqslant C( |\bar v - v^*| +|\bar \theta - \theta^*| ) \\
				&\leqslant C(|v^R - v_*| + |v^S - v^*| +|v^C - v^*| +|\theta^R - \theta_*|+|\theta^S - \theta^*| +|\theta^C - \theta^*| ) \\
				& \leqslant C(\delta_R+\delta_S +\delta_C) \leqslant C\delta_1, 
			\end{aligned}
		\end{align}
		we have
		\[
		A_{3} \leqslant \frac{4p^*\sigma^*}{3(v^*)^2}(v^S)^{-\mb X}_y \phi^2+ C(\delta_1 +\ve_1 ) (v^S)^{-\mb X}_y \phi^2.
		\]
		Similarly, by $p-\bar p =\frac{2}{3v} \z -\frac{\bar p}{v} \phi $ and $(u_1^S)^{-\mb X}_{y} <0$, it holds that
		\begin{align*}
			&\quad\left( (\theta_t^R-\sigma \theta_y^R)+(\theta_t^C-\sigma \theta_y^C)-\sigma (\theta^S)^{-\mb X}_y \right)  \Phi\left(\frac{\theta}{\bar \theta}\right) - \frac{\bar u_{1y}}{\theta}\z(p-\bar p) +\bar p \bar u_{1y} \frac{\z^2}{\theta\bar\theta}\\
			& = u^R_{1y}\left[  - p^R \Phi\left(\frac{\theta}{\bar \theta}\right) - \frac{ \z}{\theta} \left( \frac{2}{3v} \z -\frac{\bar p}{v}\phi \right)  +\bar p \frac{\z^2}{\theta\bar\theta}  \right]  \\
			&\quad -\left[ -(\theta_t^C-\sigma \theta_y^C) \Phi(\frac{\bar\theta}{\theta}) +  u_{1y}^C \frac{ \z}{\theta} \left( \frac{2}{3v} \z -\frac{\bar p}{v} \phi \right)  -\bar p u^C_{1y} \frac{\z^2}{\theta\bar\theta} \right]  \\
			&\quad +\left[ -\s(\theta^S)^{-\mb X}_y \Phi(\frac{\bar\theta}{\theta}) -(u^S_1)^{-\mb X}_y \frac{\z}{\theta} \left( \frac{2}{3v} \z -\frac{\bar p}{v} \phi \right)   \right]  +\bar p (u_1^S)^{-\mb X}_y \frac{\z^2}{\theta\bar\theta}  \\
			&\leqslant u^R_{1y} \left(  -\frac{1}{3\bar v\bar\theta}\z + \frac{\bar p}{\bar v \bar \theta} \phi\right)  \z +C \d_C(1+t)^{-1}e^{-\frac{c_0 |y+\sigma t|^2}{1+t}}|(\phi,\z)|^2\\
			&\quad + \frac{\sigma^* (v^S)^{-\mb X}_y}{2v^*\theta^*}\Big( \f23\z- 2p^* \phi\Big)\z  +C(\delta_1 +\ve_1) (|u^R_{1y} |+ |(v^S)^{-\mb X}_y | )  |(\phi,\z)|^2.
		\end{align*}
		Therefore, by combining the  estimates above together with $\bar p = \frac{R\bar\theta}{\bar v}$,
		\begin{equation}\label{(4.5)}
			\begin{array}{ll}
				\di \quad
				\int a^{-\mb X}  \left[ \f23\Big((\theta_t^R-\sigma \theta_y^R)+(\theta_t^C-\sigma \theta_y^C)-\sigma (\theta^S)^{-\mb X}_y \Big)\Phi(\frac{v}{\bar v})-\frac{\bar p\bar u_{1y}}{v\bar v}\phi^2 \right]  dy\\[3mm]
				\di\qquad
				+\int a^{-\mb X}  \left[ \Big((\theta_t^R-\sigma \theta_y^R)+(\theta_t^C-\sigma \theta_y^C)-\sigma (\theta^S)^{-\mb X}_y \Big) \Phi\left(\frac{\theta}{\bar \theta}\right) - \frac{\bar u_{1y}}{\theta}(p-\bar p)\z +\bar p \bar u_{1y} \frac{\z^2}{\theta\bar\theta}  \right]  dy \\[3mm]
				\di
				\leqslant\int a^{-\mb X} (v^S)_y^{-\mb X} \left[  \frac{4p^*\sigma^*}{3(v^*)^2} \phi^2+\frac{\sigma^*}{3v^*\theta^*}\z^2-\frac{p^*\sigma^*}{v^*\theta^*}\phi\z\right] dy-C\Lambda^R\\[3mm]
				\di \qquad
				+ C \d_C\f1{1+t} \int a^{-\mb X}  e^{-\frac{c_0 |y+\sigma t|^2}{1+t}} |(\phi,\z)|^2 dy+C(\delta_1 +\ve_1 ) \int a^{-\mb X}  (|u^R_{1y} |+ |(v^S)^{-\mb X}_y | )  |(\phi,\z)|^2 dy.
			\end{array}
		\end{equation}
		Plugging \eqref{(4.2)}, \eqref{(4.4)}, \eqref{(4.5)} into \eqref{(4.3)}, and integrating over $\mathbb R\times [0,T]$ with respect to $y$ and $t$, we have
		\begin{equation}\label{weightedre}
			\begin{array}{ll}
				\di \quad
				\int a^{-\mb X}\left(\f23\bar\theta\Phi\left(\frac{v}{\bar v}\right)+\bar\t\Phi\left(\frac{\theta}{\bar \theta}\right) + \sum_{i=1}^3 \f12 \psi_i^2 \right)dy\Bigg|_{t=0}^{t=T}+\frac{\delta_S}{H}\int_0^T|\dot{\mb X}(t)|^2dt+C\int_0^T\Lambda^Rdt\\[4mm]
				\di\qquad\qquad
				+\int_0^T\int a^{-\mb X}  \left(  \f43\frac{\mu(\bar \t)}{v} \psi_{1y}^2+\frac{\mu(\t)}{v}\sum_{i=2}^3\psi_{iy}^2+\frac{\k(\bar\t)}{v\t}\z_y^2 \right)   dydt \\[4mm]
				\di
				\leqslant\int_0^T\int a^{-\mb X}_y \left[ (p-\bar p)\psi_1-\s\left(\f23\bar\theta\Phi\left(\frac{v}{\bar v}\right)+\bar\t\Phi\left(\frac{\theta}{\bar \theta}\right) + \sum_{i=1}^3 \f12 \psi_i^2 \right) \right] dydt\\[4mm]
				\di \qquad
				+ \int_0^T\int a^{-\mb X} (v^S)_y^{-\mb X} \left[  \frac{4p^*\sigma^*}{3(v^*)^2} \phi^2+\frac{\sigma^*}{3v^*\theta^*}\z^2-\frac{p^*\sigma^*}{v^*\theta^*}\phi\z\right] dydt+\sum_{i=1}^{10} B_i+\sum_{i=1}^5 K_i,
			\end{array}
		\end{equation}
		where
		$$ B_1:=- \int_0^T\int a^{-\mb X}_y \left[  \f\z\t\left(\kappa(\t)\frac{\t_y}{v}-\kappa(\bar\t)\frac{\bar\t_y}{\bar v}\right) +\f43\psi_1\left(\mu(\t)\frac{u_{1y}}{v}-\mu(\bar\t)\frac{\bar u_{1y}}{\bar v}  \right)+\f{\mu(\t)}{v} \sum_{i=2}^3 \psi_i\psi_{iy}\right] dydt,$$
		\begin{equation*}
			\begin{array}{ll}
				\di B_2:=\int_0^T\int a^{-\mb X} \left[  \frac{\kappa(\bar\t)}{v\bar v \t}\bar \t_y \z_y \phi+\frac{\kappa(\bar\t)}{\t^2}\z \t_y \left( \frac{\t_y}{v}-\frac{\bar\t_y}{\bar v}\right)+\f43\f\z\t\left(\mu( \t)\frac{u_{1y}^2}{v}-\mu(\bar \t)\frac{\bar u_{1y}^2}{\bar v} \right)\right.\\[4mm]
				\di\qquad\qquad\qquad\qquad\qquad\qquad\qquad\qquad
				\left.+\f\z\t\f{\mu(\t)}{v} \sum_{i=2}^3 \psi_{iy}^2+\frac{4}{3}\frac{\mu(\bar\t)}{v\bar v} \bar u_{1y}\psi_{1y}\phi\right]dydt,
			\end{array}
		\end{equation*}
		$$B_3:=- \int_0^T\int a^{-\mb X}  \left( Q_1\psi_1+Q_2\f\z\t +Q_2\frac{\z^2}{\t\bar\t}\right)dy dt,$$
		$$B_4:=-\int_0^T\int a^{-\mb X} \f{\z^2}{\t\bar \t}\left[\left(\kappa(\bar\t)\frac{\bar\t_y}{\bar v} \right)_y +\f43 \mu(\bar\t)\frac{\bar u_{1y}^2}{\bar v} \right] dydt,$$
		$$B_5:=- \int_0^T\int a^{-\mb X}  \left[ \bigg(\f{\z}{\t} \bigg)_y \frac{\t_y}{v} (\k(\t)-\k(\bar \t))+\f43\frac{u_{1y}\psi_{1y}}{v}(\mu(\t)-\mu(\bar \t))\right] dydt,$$
		$$B_6:=C\d_C\int_0^T\f{1}{1+t} \int  e^{-\frac{c_0 |y+\sigma t|^2}{1+t}} |(\phi,\z)|^2 dydt,$$
		$$B_7:=C(\delta_1 +\ve_1 ) \int_0^T\int a^{-\mb X}  (|u^R_{1y} |+ |(v^S)^{-\mb X}_y | )  |(\phi,\z)|^2 dydt,$$
		$$B_8:=-\int_0^T\dot{\mb X}(t) \int a^{-\mb X}(\theta^S)^{-\mb X}_y  \left(\f23\bar\theta\Phi\left(\frac{v}{\bar v}\right)+\bar\t\Phi\left(\frac{\theta}{\bar \theta}\right) \right) dydt,$$
		$$B_9:=-\int_0^T\dot{\mb X}(t)\int a^{-\mb X}_y\left(\f23\bar\theta\Phi\left(\frac{v}{\bar v}\right)+\bar\t\Phi\left(\frac{\theta}{\bar \theta}\right) + \sum_{i=1}^3 \f12 \psi_i^2 \right)dydt,$$
		$$B_{10}:=\int_0^T\int a^{-\mb X}\f{\z^2}{\t\bar \t}\int\xi_1\left(\f12|\xi|^2-(u_1^S)^{-\mb X}\xi_1 \right)(\Pi_1^S)^{-\mb X}_y d\xi dydt,$$
		$$K_1:=-\int_0^T\int a^{-\mb X}\psi_1\int \xi_1^2\left(\Pi_1-(\Pi_1^S)^{-\mb X} \right)_y d\xi dydt,$$
		$$K_2:=-\int_0^T\int a^{-\mb X} \sum_{i=2}^3\psi_i\int \xi_1\xi_i\Pi_{1y} d\xi dydt,$$
		$$K_3:=-\int_0^T\int a^{-\mb X}\f\z\t\int \xi_1\f{|\xi|^2}2\left(\Pi_1-(\Pi_1^S)^{-\mb X} \right)_y d\xi dydt,$$
		$$K_4:=\sum_{i=2}^3\int_0^T\int a^{-\mb X} \f\z\t\psi_i\int \xi_1\xi_i\Pi_{1y} d\xi dydt,$$
		$$K_5:=\int_0^T\int a^{-\mb X}\f\z\t\left( u_1\int\xi_1^2\Pi_{1y}d\xi-(u_1^S)^{-\mb X}\int\xi_1^2(\Pi_1^S)^{-\mb X}_y d\xi\right)dydt.$$

		{\it \underline{Step 2}. Estimation on main bad terms.} In this step, we will investigate the terms on the right-hand side of \eqref{weightedre}. First we introduce a new variable $z$:
		\begin{equation*}
			z:=\frac{(v^S)^{-\mb X}-v^*}{\delta_S}.
		\end{equation*}
		Since $\mb X(t)$ is bounded on $[0,T]$ by \eqref{xprop}, it follows from $\delta_S:=v_+-v^*>0$ and $v^S_y>0$ that for any fixed $t$, the change of variable $y\in\mathbb R\mapsto z\in (0,1)$ is well-defined, together with
		\begin{equation}\label{dery}
			\frac{dz}{dy} =\frac{(v^S)^{-\mb{X}}_y}{\delta_S}>0.
		\end{equation}
		
		We need a lemma to estimate the wave interaction terms in the right-hand side of \eqref{weightedre}.
		\begin{lemma}\label{waveinteraction}
			Let $\mb X$ be the shift defined by \eqref{X(t)} and $Q_i^I(~i=1,2)$ is defined in \eqref{QI1} and \eqref{QI2}. Under the same hypotheses as in Proposition \ref{main proposition}, the following holds: for $i=1,2,$ and $t\in[0,T]$, there exist constants $C,c>0$ independent of $T, \delta_R, \delta_S$ and $\delta_C$ such that
			\begin{equation*}
				\begin{array}{ll}
					\qquad\qquad\qquad\qquad\qquad\qquad\|Q_i^I\|_{L^2} \leqslant C\delta_S (\delta_R+\delta_C) e^{-c \delta_S t}+C\delta_R\delta_Ce^{-ct},\\[3mm]
					
					\| |(v^S)^{-\mb X}_y| |\big(v^R -v_*,\theta^R -\theta_*\big) |\|_{L^2}+\| |(v^S)^{-\mb X}_y| |\big(v^C -v^*,\theta^C -\theta^*\big) |\|_{L^2}\leqslant C\delta_S^{\f32} (\delta_R+\delta_C) e^{-c \delta_S t}.
				\end{array}
			\end{equation*}
		\end{lemma}
		\begin{remark}
			The proof is almost the same as \cite[Lemma 4.2]{Kang-Vasseur-Wang-2024}, and we omit the details.
		\end{remark}

		To begin with, we concentrate on the first term of the right-hand side of \eqref{weightedre}. Direct calculation yields
		\begin{equation*}
			\begin{array}{ll}
				\di\quad
				(p-\bar p)\psi_1-\s\left(\f23\bar\theta\Phi\left(\frac{v}{\bar v}\right)+\bar\t\Phi\left(\frac{\theta}{\bar \theta}\right) + \sum_{i=1}^3 \f12 \psi_i^2 \right)\\[3mm]
				\di
				= \psi_1\left( \frac{2}{3v} \z -\frac{\bar p}{v} \phi  \right)  -\sigma \left(\f23\bar\theta\Phi\left(\frac{v}{\bar v}\right)+\bar\t\Phi\left(\frac{\theta}{\bar \theta}\right) + \sum_{i=1}^3 \f12 \psi_i^2 \right)\\[3mm]
				\di
				\leqslant\psi_1\left( \frac{2}{3v^*} \z -\frac{p^*}{v^*} \phi  \right)-\frac{\sigma^*\theta^*}{3(v^*)^2}\phi^2 - \frac{\sigma^*}{2\theta^*}\z^2 -\frac{\s^*}{2}\psi_1^2-\s\sum_{i=2}^3 \f12 \psi_i^2\\[5mm]
				\di\qquad\qquad
				+C\left( \d_S+|\phi| + |\bar v - v^*| + |\bar\theta - \theta^*| \right)  |(\phi, \psi_1,\z)|^2\\[3mm]
				\di
				=  -\frac{\sigma^*\theta^*}{3(v^*)^2} \left( \phi+\frac{1}{\sigma^*}\psi_1  \right) ^2 - \frac{\sigma^*}{2\theta^*} \left( \z-\frac{2\theta^*}{3v^*\s^*}\psi_1\right) ^2-\s\sum_{i=2}^3 \f12 \psi_i^2\\[5mm]
				\di\qquad\qquad
				+C\left( \delta_S +|\phi|+ + |(v^R - v_* ,\theta^R - \theta_*)| +|(v^C - v^*,\theta^C - \theta^*)|  \right)  |(\phi, \psi_1,\z)|^2.
			\end{array}
		\end{equation*}
		Thus
		\begin{equation}\label{B2-G}
			\begin{array}{ll}
				\di\quad
				\int_0^T\int a^{-\mb X}_y \left[ (p-\bar p)\psi_1-\s\left(\f23\bar\theta\Phi\left(\frac{v}{\bar v}\right)+\bar\t\Phi\left(\frac{\theta}{\bar \theta}\right) + \sum_{i=1}^3 \f12 \psi_i^2 \right) \right] dydt\\[4mm]
				\di
				\leqslant-\int_0^T (\Lambda_1+\Lambda_2+\Lambda_3) dt+C\delta_S \int_0^T\int a_y^{-\mb X}  |(\phi, \psi_1,\z)|^2 dydt+ C \int_0^T \int a_y^{-\mb X}  |(\phi, \psi_1,\z)|^3 dydt\\[4mm]
				\di\qquad\qquad
				+ C  \int_0^T\int a_y^{-\mb X} \big( |(v^R - v_* ,\theta^R - \theta_*)| +|(v^C - v^*,\theta^C - \theta^*)|  \big) |(\phi, \psi_1,\z)|^2 dydt,
			\end{array}
		\end{equation}
		where the good terms $\Lambda_1$, $\Lambda_2$, $\Lambda_3$ are defined by
		\begin{equation*}
			\begin{array}{ll}
				\di\Lambda_1:= \frac{\sigma^*\theta^*}{3(v^*)^2} \int a_y^{-\mb X}  \left( \phi+\frac{ 1}{\sigma^*}\psi_1  \right) ^2 dy,\\[4mm]
				\di\Lambda_2:=\frac{\sigma^*}{2\theta^*} \int a_y^{-\mb X}  \left( \z-\frac{2\theta^*}{3v^*\s^*}\psi_1\right) ^2 dy,\\[4mm]
				\di\Lambda_3:= \f{\s}2\sum_{i=2}^3 \int a_y^{-\mb X} \psi_i^2dy.
			\end{array}
		\end{equation*}
		Here we estimate the last three terms in \eqref{B2-G} in the following way. First, using \eqref{a-prime},
		$$
		\delta_S \int_0^T\int a_y^{-\mb X}  |(\phi, \psi_1,\z)|^2 dydt \leqslant C\d_S^{\f34} \Lambda^S.
		$$
		Using \eqref{a-prime} and the interpolation inequality, we have
		\begin{equation*}
			\begin{array}{ll}
				\di\quad
				\int_0^T\int a_y^{-\mb X}  |(\phi, \psi_1,\z)|^3 dydt\\[3mm]
				\di
				\leqslant C \int_0^T\int  |a_y^{-\mb X} |\Big|\phi+\frac{1}{\sigma^*}\psi_1\Big| ^3dydt+C\int_0^T\int |a_y^{-\mb X}||\psi_1|^3dydt\\[3mm]
				\di\quad\quad+C \int_0^T\int |a_y^{-\mb X} |\Big| \z-\frac{2\theta^*}{3v^*\s^*}\psi_1 \Big|^3 dydt\\[3mm]
				\di
				\leqslant C\ve_1\int_0^T (\Lambda_1+\Lambda_2)dt+C \delta_S^{-\f14} \int_0^T\int (v^S)_y^{-\mb X}   \|\psi_1\|_{L^\infty}^2 |\psi_1| dydt\\[4mm]
				\di
				\leqslant C \ve_1\int_0^T(\Lambda_1+\Lambda_2)dt+C\delta_S^{-\f14} \int_0^T \|\psi_{1y}\|_{L^2} \|\psi_1\|_{L^2} \int (v^S)_y^{-\mb X}   |\psi_1| dydt\\[4mm]
				\di
				\leqslant C \ve_1\int_0^T(\Lambda_1+\Lambda_2+\|\psi_{1y}\|_{L^2}^2)dt+C\ve_1 \delta_S^{-\f12}  \int_0^T\left( \int (v^S)_y^{-\mb X}   |\psi_1| dy\right) ^2dt\\[4mm]
				\di
				\leqslant C \ve_1\int_0^T(\Lambda_1+\Lambda_2 + \|\psi_{1y}\|_{L^2}^2 +\Lambda^S)dt.
			\end{array}
		\end{equation*}
		Similarly, using the interpolation inequality and Lemma \ref{waveinteraction}, we have
		\begin{equation*}
			\begin{array}{ll}
				\di\quad
				\int_0^T\int a_y^{-\mb X} \big( |(v^R - v_* ,\theta^R - \theta_*)| +|(v^C - v^*,\theta^C - \theta^*)|  \big) |(\phi, \psi_1,\z)|^2 dydt\\[3mm]
				\di
				\leqslant C\delta_1  \int_0^T(\Lambda_1+\Lambda_2)dt+C \d_S^{-\f14} \int_0^T\int  (v^S)_y^{-\mb X}  \big( |(v^R - v_* ,\theta^R - \theta_*)| +|(v^C - v^*,\theta^C - \theta^*)|  \big) \psi_1^2 dydt\\[3mm]
				\di
				\leqslant C\delta_1  \int_0^T(\Lambda_1+\Lambda_2)dt\\[3mm]
				\di\quad
				+C\int_0^T\delta_S^{-\f14} \|\psi_1\|_{L^4}^2 \Big[\| |(v^S)_y^{-\mb X}| | (v^R - v_* ,\theta^R - \theta_*)|\|_{L^2}+\| |(v^S)_y^{-\mb X} | | (v^C - v^* ,\theta^C - \theta^*)|\|_{L^2}  \Big]dt \\[4mm]
				\di
				\leqslant C\delta_1  \int_0^T(\Lambda_1+\Lambda_2)dt+C\int_0^T \|\psi_{1y}\|_{L^2}^{\f12}  \|\psi_1\|_{L^2}^{\f32}  \delta_S(\delta_R+\delta_C)e^{-c\delta_S t}dt\\[4mm]
				\di
				\leqslant  C\delta_1  \int_0^T(\Lambda_1+\Lambda_2)dt+C\ve_1 \int_0^T\|\psi_{1y}\|_{L^2}^2dt +C\int_0^T \delta_S^{\f43} ( \delta_R^{\f43}+ \delta_C^{\f43}) e^{-c\delta_S t}dt\\[4mm]
				\di
				\leqslant  C\delta_1  \int_0^T(\Lambda_1+\Lambda_2)dt+C\ve_1 \int_0^T\|\psi_{1y}\|_{L^2}^2dt+C\delta_S^{\f13} ( \delta_R^{\f43}+ \delta_C^{\f43}).
			\end{array}
		\end{equation*}
		Combining all the estimates above and taking $\d_1$, $\ve_1$ suitably small, it holds that
		\begin{equation}\label{goodfromweight}
			\begin{array}{ll}
				\di\quad
				\int_0^T\int a^{-\mb X}_y \left[ (p-\bar p)\psi_1-\s\left(\f23\bar\theta\Phi\left(\frac{v}{\bar v}\right)+\bar\t\Phi\left(\frac{\theta}{\bar \theta}\right) + \sum_{i=1}^3 \f12 \psi_i^2 \right) \right] dydt\\[4mm]
				\di\leqslant-C\int_0^T (\Lambda_1+\Lambda_2+\Lambda_3) dt+C(\delta_S^{\f34}+\ve_1)\int_0^T\left( \Lambda^S+\|\psi_{1y}\|_{L^2}^2\right)dt+C\delta_S^{\f13} ( \delta_R^{\f43}+ \delta_C^{\f43}).
			\end{array}
		\end{equation}
		To estimate the second term on the right-hand side of \eqref{weightedre}, we find by Cauchy inequality that
		\begin{equation*}
			\begin{array}{ll}
				\di\quad
				\int a^{-\mb X} (v^S)_y^{-\mb X}\frac{4p^*\sigma^*}{3(v^*)^2} \phi^2 dy\\[3mm]
				\di
				=\frac{4p^*\sigma^*}{3(v^*)^2}  \int a^{-\mb X} (v^S)_y^{-\mb X} \Big|\Big(\phi+\frac{1}{\sigma^*}\psi_1\Big) - \frac{1}{\sigma^*}\psi_1 \Big|^2 dy\\[5mm]
				\di
				\leqslant\frac{4p^*\sigma^*}{3(v^*)^2}\left(1+\delta_S^{\f34}+\delta_S^{\f18}\right) \int  (v^S)_y^{-\mb X} \Big| \frac{1}{\sigma^*}\psi_1 \Big|^2 dy+C\delta_S^{-\f18} \int  (v^S)_y^{-\mb X} \Big|\phi+\frac{1}{\sigma^*}\psi_1 \Big|^2 dy,
			\end{array}
		\end{equation*}
		Using this method, we can conclude that the second term on the right-hand side of \eqref{weightedre} can be estimated by
		\begin{equation}\label{mainbadterm}
			\begin{array}{ll}
				\di\quad
				\int_0^T\int a^{-\mb X} (v^S)_y^{-\mb X} \left[  \frac{4p^*\sigma^*}{3(v^*)^2} \phi^2+\frac{\sigma^*}{3v^*\theta^*}\z^2-\frac{p^*\sigma^*}{v^*\theta^*}\phi\z\right] dydt\\[4mm]
				\di\leqslant\alpha^* \left(1+C\delta_S^{\f18}\right) \delta_S \int_0^T\int_0^1 \psi_1^2 dzdt+ C\delta_S^{\f18}\int_0^T (\Lambda_1 + \Lambda_2)dt,
			\end{array}
		\end{equation}
		where $\di\alpha^*:=\frac{20p^*}{9(v^*)^2\s^*}$ is a positive constant. 
		
		The approach to addressing the primary challenges, $\int_0^T\int_0^1 \psi_1^2 dzdt$, is to use Poincar\'e inequality. In other words, we aim for it to be controlled by $\int_0^T\int \psi_{1y}^2 dydt$, which is in the left-hand side of \eqref{weightedre}. The Poincar\'e inequality is as follows.
		\begin{lemma}\label{lem-poin}
			\cite{Kang-Vasseur-2021} For any $f:[0,1]\to \mathbb R$ satisfying $\di\int_0^1 z(1-z)|f'|^2 dz<\infty$, 
			\begin{equation*}
				\int_0^1\left| f-\int_0^1 f dz \right| ^2 dz\leqslant \frac{1}{2}\int_0^1 z(1-z)|f'|^2 dz.
			\end{equation*}
		\end{lemma}
		
		Consider $\di\int a^{-\mb X}  \f43\frac{\mu(\bar \t)}{v} \psi_{1y}^2dy$. By \eqref{dery}, we have
		\begin{equation}\label{D1}
			\di\int a^{-\mb X}  \f43\frac{\mu(\bar \t)}{v} \psi_{1y}^2dy\geqslant\int\f43\frac{\mu(\bar \t)}{v} \psi_{1y}^2dy=\int_0^1\f43\frac{\mu(\bar \t)}{v}\psi_{1z}^2\frac{dz}{dy} dz=\int_0^1\f43\frac{\mu(\bar \t)}{v}\psi_{1z}^2\frac{(v^S)^{-\mb{X}}_y}{\delta_S} dz.
		\end{equation}
		From \eqref{VS2}$_1$, together with \eqref{sigmadef} and the facts that
		\begin{equation}\label{(5.17)}
			\frac{1}{\delta_S^2}=\frac{z(1-z)}{(v^S-v^*)(v_+-v^S)},
		\end{equation}
		we have
		\begin{equation}\label{(4.14)}
			\begin{array}{ll}
				\di\quad
				\int_0^1\f43\frac{\mu(\bar \t)}{v}\psi_{1z}^2\frac{v^S_y}{\delta_S} dz\\[4mm]
				\di= \int_0^1\f43\frac{\mu(\bar \t)}{v}\psi_{1z}^2\frac{1}{v_+-v^*}\left(-\f43\frac{\mu(\t^S)\s}{v^S} \right) ^{-1}\left((p^S-p^*)+\sigma^2(v^S-v^*)+\int\xi_1^2\Pi_1^Sd\x \right) dz\\[4mm]
				\di=-\int_0^1\frac{\mu(\bar \t)}{\mu(\t^S)}\frac{v^S}{v}\psi_{1z}^2\frac{1}{\s\d_S^2}\left((v_+-v^*)(p^S-p^*)+(p^*-p_+)(v^S-v^*) +(v_+-v^*)\int\xi_1^2\Pi_1^Sd\x \right) dz\\[4mm]
				\di=-\int_0^1\frac{\mu(\bar \t)}{\mu(\t^S)}\frac{v^S}{v}\psi_{1z}^2\frac{1}{\s\d_S^2}\left((p^S-p_+) ( v^S -v^*) + (v_+- v^S)(p^S-p^*) +(v_+-v^*)\int\xi_1^2\Pi_1^Sd\x \right) dz\\[4mm]
				\di=\int_0^1\frac{\mu(\bar \t)}{\mu((\t^S)^{-\mb X})}\frac{(v^S)^{-\mb X}}{v}z(1-z)\psi_{1z}^2\frac{1}{\s}\left(\frac{p^S-p_+}{v^S-v_+}-\frac{p^S-p^*}{v^S-v^*}-\frac{\di \int\xi_1^2\Pi_1^Sd\x}{v^S-v^*}-\frac{\di \int\xi_1^2\Pi_1^Sd\x}{v_+-v^S}\right) dz.
			\end{array}
		\end{equation}
		Note that we drop the mark "$-\mb X$" from \eqref{(5.17)} to \eqref{(4.20)} for simplicity. By \eqref{2ndorder}, it holds that
		\begin{equation}\label{pestimate}
			\begin{array}{ll}
				\di \frac{p^S-p_+}{v^S-v_+}-\frac{p^S-p^*}{v^S-v^*}= \frac{ p^*}{9(v^*)^2}\frac{200\mu(\t^*)}{10\mu(\t^*)+3\k(\t^*)}  \d_S+O(\delta_S^2).
			\end{array}
		\end{equation}
		For$\di\frac{\di \int\xi_1^2\Pi_1^Sd\x}{v^S-v^*}$, by \eqref{shock-macro} and \eqref{shock-micro}, we find
		\begin{equation}\label{(4.20)}
			\frac{\di \int\xi_1^2\Pi_1^Sd\x}{v^S-v^*}=\frac{\di \int\xi_1^2\Pi_1^Sd\x}{v^S_y}\frac{v^S_y}{v^S-v^*}=O(\d_S^2),
		\end{equation}
		and same  conclusion holds for $\di\frac{\di \int\xi_1^2\Pi_1^Sd\x}{v_+-v^S}$. Thus, substituting \eqref{(4.14)}, \eqref{pestimate} and \eqref{(4.20)} into \eqref{D1}, by Poincar\'e inequality, it holds that
		\begin{equation}\label{(4.22)}
			\begin{array}{ll}
				\di\int a^{-\mb X}  \f43\frac{\mu(\bar \t)}{v} \psi_{1y}^2dy&\di\geqslant\alpha^*\frac{10\mu(\t^*)}{10\mu(\t^*)+3\k(\t^*)}(1-C(\d_1+\ve_1))\d_S\int_0^1 z(1-z)\psi_{1z}^2dz\\[5mm]
				&\di \geqslant2\alpha^*\frac{10\mu(\t^*)}{10\mu(\t^*)+3\k(\t^*)}(1-C(\d_1+\ve_1))\d_S\left[\int_0^1 \psi_1^2 dz-\left( \int_0^1\psi_1dz \right)  ^2 \right].
			\end{array}
		\end{equation}
		
		Comparing \eqref{(4.22)} to \eqref{mainbadterm}, we find that the coefficient before $\int_0^1 \psi_1^2 dz$ may not be large enough to control the bad term in \eqref{mainbadterm}. However, another good term $\int \z_y^2   dy$ may help.
		
		Similar to \eqref{(4.22)}, we have
		\begin{equation}\label{(4.23)}
			\begin{array}{ll}
				\quad\di\int a^{-\mb X}   \frac{\k(\bar\t)}{v\t}\z_y^2   dy\\[3mm]
				\di\geqslant\int \frac{\k(\bar\t)}{v\t}\z_y^2   dy\\[3mm]
				\di=\int_0^1 \frac{\k(\bar\t)}{v\t}\z_z^2 \frac{v^S_y}{\delta_S} dz\\[3mm]
				\di=\int_0^1 \frac{\k(\bar\t)}{v\t}\z_z^2\left( \f43\frac{\mu(\bar \t)}{v}\right)^{-1}\f43\frac{\mu(\bar \t)}{v} \frac{v^S_y}{\delta_S}dz\\[3mm]
				\di\geqslant \f34 \frac{\k(\t^*)}{\t\mu(\t^*)}\alpha^*\frac{10\mu(\t^*)}{10\mu(\t^*)+3\k(\t^*)}(1-C(\d_1+\ve_1))\d_S\int_0^1z(1-z)\z_z^2dz\\[3mm]
				\di\geqslant \f34 \frac{\k(\t^*)}{\t\mu(\t^*)}\alpha^*\frac{10\mu(\t^*)}{10\mu(\t^*)+3\k(\t^*)}(1-C(\d_1+\ve_1))\d_S\left[\int_0^1 \z^2 dz-\left( \int_0^1\z dz \right)  ^2 \, \right].
			\end{array}
		\end{equation}
		Here we again drop the mark "$-\mb X$" for simplicity. Now it remains to relate $\di \left[\int_0^1 \z^2 dz-\left( \int_0^1\z dz \right)  ^2 \,\right]$ to $\di \left[\int_0^1 \psi_1^2 dz-\left( \int_0^1\psi_1dz \right)  ^2 \,\right]$. Cauchy inequality and H{\"o}lder inequality yield that
		\begin{equation}\label{(4.24)}
			\begin{array}{ll}
				\quad\di\int_0^1 \z^2 dz-\left( \int_0^1\z dz \right)  ^2 \\[3mm]
				\di=\int_0^1 \left( \z-\frac{2\theta^*}{3v^*\s^*}\psi_1+\frac{2\theta^*}{3v^*\s^*}\psi_1\right)^2  dz-\left( \int_0^1\left( \z-\frac{2\theta^*}{3v^*\s^*}\psi_1+\frac{2\theta^*}{3v^*\s^*}\psi_1 \right)  dz \right)^2\\[3mm]
				\di\geqslant\left( \frac{4(\t^*)^2}{9(v^*\s^*)^2}-\delta_S^{\f18}\right) \int_0^1 \psi_1^2 dz-C\delta_S^{-\f18}\int_0^1 \left( \z-\frac{2\theta^*}{3v^*\s^*}\psi_1\right)^2  dz\\[3mm]
				\di\qquad\qquad\qquad\qquad-2\left( \int_0^1\left( \z-\frac{2\theta^*}{3v^*\s^*}\psi_1 \right)  dz \right)^2-\frac{8(\t^*)^2}{9(v^*\s^*)^2}\left( \int_0^1\psi_1 dz \right)  ^2\\[3mm]
				\di
				\geqslant \left( \frac{2\t^*}{5}-\delta_S^{\f18}\right) \int_0^1 \psi_1^2 dz-C\delta_S^{-\f18}\int \d_S^{-1} (v^S)^{-\mb X}_y \left( \z-\frac{2\theta^*}{3v^*\s^*}\psi_1\right)^2  dy\\[3mm]
				\di\qquad\qquad\qquad\qquad
				-2\int_0^1\left( \z-\frac{2\theta^*}{3v^*\s^*}\psi_1 \right)^2  dz -\frac{4\t^*}{5}\left( \int_0^1\psi_1 dz \right)  ^2\\[3mm]
				\di\geqslant\left( \frac{2\t^*}{5}-\delta_S^{\f18}\right) \int_0^1 \psi_1^2 dz-C\d_S^{-\f78}\Lambda_2-\frac{4\t^*}{5}\left( \int_0^1\psi_1 dz \right)  ^2.
			\end{array}
		\end{equation}
		Therefore, combining \eqref{(4.22)}, \eqref{(4.23)}, \eqref{(4.24)} yields that
		\begin{equation}\label{(4.25)}
			\begin{array}{ll}
				\quad\di\int a^{-\mb X} \left(  \f43\frac{\mu(\bar \t)}{v} \psi_{1y}^2dy+\frac{\k(\bar\t)}{v\t}\z_y^2\right)dy \\[3mm]
				\di\geqslant2\alpha^*(1-C(\d_1+\ve_1)-C\d_S^{\f18})\d_S\int_0^1 \psi_1^2 dz-C\d_S^{\f18}\Lambda_2-4\alpha^*\d_S\frac{5\mu(\t^*)+3\k(\t^*)}{10\mu(\t^*)+3\k(\t^*)}\left( \int_0^1\psi_1 dz \right)  ^2.
			\end{array}
		\end{equation}
		The first term on the right-hand side of \eqref{(4.25)} is aimed to control the the bad term in \eqref{mainbadterm}, while the last term needs to be absorbed by $\frac{\delta_S}{H}\int_0^T|\dot{\mb X}(t)|^2dt$. Note that
		$$\int a^{-\mb X}(v^S)^{-\mb{X}}_y \frac{ \bar p}{\bar v}\phi dy=-\int a^{-\mb X} \frac{(v^S)^{-\mb{X}}_y \bar p }{\bar v \s^*}\psi_1 dy +\int a^{-\mb X} \frac{(v^S)^{-\mb{X}}_y \bar p}{\bar v}\left( \phi+\frac{1}{\sigma^*}\psi_1\right)dy. $$
		Then using \eqref{abound}, \eqref{a-prime}, and \eqref{perror}, we have
		$$\left|\int a^{-\mb X}(v^S)^{-\mb{X}}_y \frac{ \bar p}{\bar v}\phi dy + \frac{p^*\delta_S}{v^*\sigma^*}\int_0^1 \psi_1 dz\right| \leqslant C\delta_S (\d_S^{\f34}+\delta_1)\int_0^1|\psi_1| dz+ C\d_S^{\f14}\int a^{-\mb X}_y \left|\phi+\frac{1}{\sigma^*}\psi_1\right| dy.$$
		Similar to the above, together with $\s^* =\sqrt\frac{5 p^*}{3v^*}$, we can obtain
		\begin{equation*}
			\begin{array}{ll}
				\di\left| \dot{\mb X} -2\sigma^*H\int_0^1 \psi_1 dz\right| &\di \leqslant C (\d_S^{\f34}+\delta_1)\int_0^1|\psi_1| dz+ C\d_S^{-\f34}\int a^{-\mb{X}}_y \left|\phi+\frac{1}{\sigma^*}\psi_1\right| dy\\[4mm]
				&\di\quad + C\d_S^{-\f34}\int a^{-\mb{X}}_y \left| \z-\frac{2\theta^*}{3v^*\s^*}\psi_1 \right| dy,
			\end{array}
		\end{equation*}
		which yields
		\begin{equation*}
			\begin{array}{ll}
				\quad\di2(\sigma^*)^2H^2\left(\int_0^1 \psi_1 dz\right)^2 - |\dot{\mb X}|^2 \\[3mm]
				\di\leqslant\left( \left|2\sigma^*H \int_0^1 \psi_1 dz\right| - |\dot{\mb X} | \right)^2\\[3mm]
				\di\leqslant C(\d_S^{\f34}+\delta_1)^2 \int_0^1|\psi_1|^2 dy + C\d_S^{-\f32}  ( \Lambda_1 + \Lambda_2)\int a^{-\mb{X}}_y dy\\[3mm]
				\di\leqslant C(\d_S^{\f34}+\delta_1)^2 \int_0^1|\psi_1|^2 dy +C\d_S^{-\f34}  ( \Lambda_1 + \Lambda_2).
			\end{array}
		\end{equation*}
		Thus we can conclude
		\begin{equation}\label{Xgood}
			\frac{\delta_S}{2H} |\dot{\mb X}|^2 \geqslant (\sigma^*)^2H\delta_S\left(\int_0^1 \psi_1 dz\right)^2 -C\d_S(\d_S^{\f34}+\delta_1)^2 \int_0^1\psi_1^2 dz
			- C\d_S^{\f14}  ( \Lambda_1 + \Lambda_2). 
		\end{equation}
		Therefore, by choosing $H=\frac{3\alpha^*}{(\s^*)^2}\f{5\mu (\t^*)+3\kappa (\t^*)}{10\mu (\t^*)+3\kappa (\t^*)}$,  combining \eqref{goodfromweight}, \eqref{mainbadterm}, \eqref{(4.25)}, and \eqref{Xgood},  and using the smallness of $\d_1$, $\d_S$, $\ve_1$, we have
		\begin{equation}\label{(4.33)}
			\begin{array}{ll}
				\di \quad\frac{\delta_S}{2H}\int_0^T|\dot{\mb X}(t)|^2dt+\f34\int_0^T\int a^{-\mb X}  \left(  \f43\frac{\mu(\bar \t)}{v} \psi_{1y}^2+\frac{\k(\bar\t)}{v\t}\z_y^2 \right)   dydt\\[4mm]
				\di\quad\qquad
				-\int_0^T\int a^{-\mb X}_y \left[ (p-\bar p)\psi_1-\s\left(\f23\bar\theta\Phi\left(\frac{v}{\bar v}\right)+\bar\t\Phi\left(\frac{\theta}{\bar \theta}\right) + \sum_{i=1}^3 \f12 \psi_i^2 \right) \right] dydt\\[4mm]
				\di\quad\qquad
				-\int_0^T\int a^{-\mb X} (v^S)_y^{-\mb X} \left[  \frac{4p^*\sigma^*}{3(v^*)^2} \phi^2+\frac{\sigma^*}{3v^*\theta^*}\z^2-\frac{p^*\sigma^*}{v^*\theta^*}\phi\z\right] dydt\\[4mm]
				\di
				\geqslant C\int_0^T (\Lambda_1+\Lambda_2+\Lambda_3) dt+C\d_S\int_0^T\int_0^1 \psi_1^2 dzdt-C(\delta_S^{\f34}+\ve_1)\int_0^T \Lambda^Sdt -C\delta_S^{\f13} ( \delta_R^{\f43}+ \delta_C^{\f43}).
			\end{array}
		\end{equation}
		Finally, due to
		$$\delta_S \int_0^T\int_0^1 \psi_1^2 dzdt = \int_0^T\int  (v^S)^{-\mb X}_y \psi_1^2 dydt,$$
		we can take
		\begin{equation*}
			\begin{array}{ll}
				\di\int_0^T \Lambda_1 dt&\di= C\d_S^{-\f14}\int_0^T \int  (v^S)^{-\mb X}_y\left( \phi+\frac{ 1}{\sigma^*}\psi_1  \right) ^2 dydt\\[3mm]
				&\di\geqslant C\int_0^T \int  (v^S)^{-\mb X}_y \phi^2dydt-\tilde\nu\int_0^T\int  (v^S)^{-\mb X}_y \psi_1^2 dydt,
			\end{array}
		\end{equation*}
		with $\tilde\nu$ suitably small. Likewise, 
		$$\int_0^T \Lambda_2 dt\geqslant C\int_0^T \int  (v^S)^{-\mb X}_y \z^2dydt-\tilde\nu\int_0^T\int  (v^S)^{-\mb X}_y \psi_1^2 dydt,$$
		$$\int_0^T \Lambda_3 dt\geqslant C\int_0^T \int  (v^S)^{-\mb X}_y \left( \psi_2^2+\psi_3^2\right) dydt.$$
		As a consequence, we can conclude from \eqref{(4.33)} that
		\begin{equation}\label{4.2conclusion}
			\begin{array}{ll}
				\di \quad\frac{\delta_S}{2H}\int_0^T|\dot{\mb X}(t)|^2dt+\f34\int_0^T\int a^{-\mb X}  \left(  \f43\frac{\mu(\bar \t)}{v} \psi_{1y}^2+\frac{\k(\bar\t)}{v\t}\z_y^2 \right)   dydt\\[4mm]
				\di\quad\qquad
				-\int_0^T\int a^{-\mb X}_y \left[ (p-\bar p)\psi_1-\s\left(\f23\bar\theta\Phi\left(\frac{v}{\bar v}\right)+\bar\t\Phi\left(\frac{\theta}{\bar \theta}\right) + \sum_{i=1}^3 \f12 \psi_i^2 \right) \right] dydt\\[4mm]
				\di\quad\qquad
				-\int_0^T\int a^{-\mb X} (v^S)_y^{-\mb X} \left(   \frac{4p^*\sigma^*}{3(v^*)^2} \phi^2+\frac{\sigma^*}{3v^*\theta^*}\z^2-\frac{p^*\sigma^*}{v^*\theta^*}\phi\z\right)  dydt\\[4mm]
				\di
				\geqslant C\int_0^T  \Lambda^Sdt -C\delta_S^{\f13} ( \delta_R^{\f43}+ \delta_C^{\f43}).
			\end{array}
		\end{equation}
		
		\
		
		{\it \underline{Step 3}. Estimation on $B_i$$(i=1,\cdots, 10)$} With the priori groundwork laid, the subsequent energy estimation method becomes relatively simple. The key point is to make full use of Cauchy inequality. We will estimate $B_i$$(i=1,\cdots, 10)$ one by one.
		
		Since Lemma \ref{rarefaction}, Lemma \ref{contact} and Lemma \ref{Lemma-shock} yield
		\begin{equation*}
			\begin{array} {ll}
				\di\quad\int_0^T\int  |\bar \theta_y|^2  |(\phi,\z)|^2 dy dt\\[3mm]
				\di\leqslant C\int_0^T\int \left(  |\theta^R_y|^2+|\theta^C_y|^2+|(\theta^S)^{-\mb X}_y|^2\right)  |(\phi,\z)|^2 dy dt\\[3mm]
				\di\leqslant C\delta_1 \int_0^T(\Lambda^R+\Lambda^S)dt+ C\delta_C^2 \int_0^T(1+t)^{-1}\int  e^{-\frac{2c_0|y+\sigma t|^2}{1+t}}|(\phi,\z)|^2 dydt ,
			\end{array}
		\end{equation*}
		we have for the first term in $B_1$,
		\begin{equation*}
			\begin{array}{ll}
				\di \quad
				\left| - \int_0^T\int a^{-\mb X}_y \left[  \f\z\t\left(\kappa(\t)\frac{\t_y}{v}-\kappa(\bar\t)\frac{\bar\t_y}{\bar v}\right) \right] dydt\right| \\[4mm]
				\di
				=\left| \int_0^T\int \d_S^{-\f14} (v^S)^{-\mb X}_y\f\z\t\left[ (\k(\t)-\k(\bar \t))\frac{\t_y}{v}+\k(\bar \t) \left(\frac{\z_y}{v}+\bar \t_y\left(\f1v-\f1{\bar v} \right)  \right) \right] dydt\right| \\[4mm]
				\di
				\leqslant C\int_0^T \int \d_S^{-\f14} (v^S)^{-\mb X}_y|\z\z_y|dydt+C\int_0^T \int \d_S^{-\f14}(v^S)^{-\mb X}_y| \bar \t_y||(\phi,\z)|^2dydt\\[4mm]
				\di
				\leqslant \f1{200}\int_0^T \int\frac{\k(\bar \t)}{v\t}\z_y^2dydt+C\delta_1\int_0^T(\Lambda^R+\Lambda^S)dt+C\delta_C^2 \int_0^T(1+t)^{-1}\int  e^{-\frac{2c_0|y+\sigma t|^2}{1+t}}|(\phi,\z)|^2 dydt.
			\end{array}
		\end{equation*}
		Using this method, it is easy to check that
		\begin{equation*}
			\begin{array}{ll}
				\di 
				B_1 \leqslant \f1{200}\int_0^T\int  \left(  \f43\frac{\mu(\bar \t)}{v} \psi_{1y}^2+\frac{\mu(\t)}{v}\sum_{i=2}^3\psi_{iy}^2+\frac{\k(\bar\t)}{v\t}\z_y^2 \right)   dydt+C\delta_1\int_0^T(\Lambda^R+\Lambda^S)dt\\[4mm]
				\di\qquad\qquad\qquad\qquad\qquad\qquad\qquad
				+C\delta_C^2 \int_0^T(1+t)^{-1}\int  e^{-\frac{2c_0|y+\sigma t|^2}{1+t}}|(\phi,\z)|^2 dydt.
			\end{array}
		\end{equation*}
		Similarly, by the priori assumption $\left\| (\phi,\psi,\z)\right\| _{L^\i}\leqslant C\ve_1$,
		\begin{equation*}
			\begin{array}{ll}
				\di 
				B_2\leqslant \f1{200}\int_0^T\int  \left(  \f43\frac{\mu(\bar \t)}{v} \psi_{1y}^2+\frac{\mu(\t)}{v}\sum_{i=2}^3\psi_{iy}^2+\frac{\k(\bar\t)}{v\t}\z_y^2 \right)   dydt+C\delta_1\int_0^T(\Lambda^R+\Lambda^S)dt\\[4mm]
				\di\qquad\qquad\qquad\qquad\qquad\qquad\qquad
				+C\delta_C^2 \int_0^T(1+t)^{-1}\int  e^{-\frac{2c_0|y+\sigma t|^2}{1+t}}|(\phi,\z)|^2 dydt.
			\end{array}
		\end{equation*}
		For $B_3$, it holds that
		\begin{equation*}
			\begin{array}{ll}
				\di 
				B_3\leqslant C\int_0^T( \|Q_1^I\|_{L^2}+\|Q_1^C\|_{L^2}+\|Q_2^I\|_{L^2}+\|Q_2^C\|_{L^2} ) \|(\psi_1,\z)\|_{L^2}  dt\\[4mm]
				\di\qquad\qquad\qquad\qquad\qquad\qquad\qquad
				+C\int_0^T\|(Q_1^R, Q_2^R)\|_{L^1}\|(\psi_{1}, \z)\|_{L^\i} dt.
			\end{array}
		\end{equation*}
		Then, by Lemma \ref{waveinteraction} and the facts that
		$$\int_0^{+\i}\|Q_i^R\|_{L^1}^{\f43}dt \leqslant C\delta_R^{\f13},\quad i=1,2,$$
		$$\|Q_1^C\|_{L^2}\leqslant C \delta_C (1+t)^{-\frac 54},\qquad \|Q_2^C\|_{L^2}\leqslant C \delta_C (1+t)^{-\frac 74},$$
		it can readily be seen that
		\begin{equation*}
			\begin{array}{ll} 
				B_3&\di\leqslant C\int_0^T\left(\delta_S (\delta_R+\delta_C) e^{-C \delta_S t}+\delta_R\delta_Ce^{-Ct}+\delta_C (1+t)^{-\frac 54}\right)\|(\psi_1,\z)\|_{L^2} dt  \\[4mm]
				&\di\qquad+C\int_0^T\|(Q_1^R, Q_2^R)\|_{L^1}\|(\psi_{1}, \z)\|_{L^2}^{\f12}\|(\psi_{1y}, \z_y)\|_{L^2}^{\f12} dt\\[4mm]
				&\di\leqslant C\d_1+\f1{200}\int_0^T\int  \left(  \f43\frac{\mu(\bar \t)}{v} \psi_{1y}^2+\frac{\k(\bar\t)}{v\t}\z_y^2 \right)   dydt+C\int_0^T\|(Q_1^R, Q_2^R)\|_{L^1}^{\frac 43}\|(\psi_1, \z)\|_{L^2}^{\frac 23}dt\\[4mm]
				&\di\leqslant\f1{200}\int_0^T\int  \left(  \f43\frac{\mu(\bar \t)}{v} \psi_{1y}^2+\frac{\k(\bar\t)}{v\t}\z_y^2 \right)   dydt+C\d_1^{\f13}\sup_{t\in[0,T]}\|(\psi_1, \z)\|_{L^2}^{\frac 23}+C\d_1.
			\end{array}
		\end{equation*}
		The estimates for $B_4$, $B_5$, $B_7$ are similar to $B_1$, $B_2$. For $B_8$,
		\begin{equation*}
			\begin{array}{ll}
				\di 
				B_8&\di\leqslant \frac{\d_S}{16H}\int_0^T|\dot{\mb X}(t)|^2dt +\frac{C}{\d_S}\int_0^T\left(\int|(\theta^S)^{-\mb X}_y||(\phi,\z)|^2 dy\right) ^2 dt\\[4mm]
				&\di
				\leqslant \frac{\d_S}{16H}\int_0^T|\dot{\mb X}(t)|^2dt +C\d_S\ve_1^2\int_0^T\Lambda^Sdt,
			\end{array}
		\end{equation*}
		and similarly for $B_9$,
		\begin{equation*}
			\begin{array}{ll}
				\di 
				B_9
				\leqslant \frac{\d_S}{16H}\int_0^T|\dot{\mb X}(t)|^2dt +C\ve_1\int_0^T\Lambda^Sdt.
			\end{array}
		\end{equation*}
		At last, by Lemma \ref{Lemma-shock}, we can deduce
		$$B_{10}\leqslant C\d_S^2\int_0^T (v^S)^{-\mb X}_y\z^2 dydt\leqslant C\d_S\int_0^T\Lambda^Sdt.$$
		
		Therefore, combining all the estimates above, it holds that
		\begin{equation}\label{(4.36)}
			\begin{array}{ll}
				\di 
				\sum_{i=1}^{10} B_i\leqslant \frac{\d_S}{8H}\int_0^T|\dot{\mb X}(t)|^2dt +C(\d_1+\ve_1) \int_0^T(\Lambda^S+\Lambda^R)dt\\[4mm]
				\di\qquad\qquad
				+C\delta_C \int_0^T(1+t)^{-1}\int  e^{-\frac{c_0|y+\sigma t|^2}{1+t}}|(\phi,\z)|^2 dydt\\[4mm]
				\di\qquad\qquad
				+\f1{20}\int_0^T\int  \left(  \f43\frac{\mu(\bar \t)}{v} \psi_{1y}^2+\frac{\mu(\t)}{v}\sum_{i=2}^3\psi_{iy}^2+\frac{\k(\bar\t)}{v\t}\z_y^2 \right)   dydt\\[4mm]
				\di\qquad\qquad
				+C\d_1^{\f13}\sup_{t\in[0,T]}\|(\psi_1, \z)\|_{L^2}^{\frac 23}+C\d_1.
			\end{array}
		\end{equation}
		
		\

		{\it \underline{Step 4}. Estimation on $K_i$$(i=1,\cdots, 5)$.} In the following, we turn to estimate $K_i$$(i=1,\cdots, 5)$. Here, we only estimate $K_1$, as $K_i$$(i=2,\cdots, 5)$ can be estimated similarly. Let $\mb{M}_{_\#}$ be a global Maxwellian with the state $(v_{_\#},u_{_\#},\t_{_\#})$ satisfying $\f12\t<\t_\#<\t$ and
		$|v-v_{_\#}|+|u-u_{_\#}|+|\t-\t_{_\#}|\leqslant \eta_0$ so that Lemma \ref{Lemma 4.3} holds, and Lemmas \ref{Lemma-shock}  holds with $\mb{M}_0$ being replaced by $\mb{M}_{_\#}$.  By the definition of $\Pi_1$ and $\Pi_1^{S}$, we have
		\begin{equation*}
			\begin{array}{ll}
				\di \Pi_1-(\Pi_1^{S})^{-\mb X} \di
				=&\di {\mb{L}}_\mb{M}^{-1}\left[ \wt{\mb{G}}_t-\sigma \wt{\mb{G}}_y -\f{u_1}{v}\wt{\mb{G}}_{y}+\f{1}{v}\mb{P}_1(\x_1\wt{\mb{G}}_{y})-Q(\wt{\mb{G}},\wt{\mb{G}})\right]+J \\[3mm]
				&\di
				-{\mb{L}}_\mb{M}^{-1}\left[Q(\wt{\mb{G}},(\mb{G}^{S})^{-\mb X})+ Q(({\mb{G}}^{S})^{-\mb X},\wt{\mb{G}})\right] -\dot{\mb X}(t){\mb{L}}_\mb{M}^{-1}(\mb{G}^{S})^{-\mb X}_y,
			\end{array}
		\end{equation*}
		where
		\begin{equation}\label{J3}
			\begin{array}{ll}
				J=&\di
				\left( {\mb{L}}_\mb{M}^{-1}-{\mb{L}}_\mb{M^{S}}^{-1}\right) \left[ -\s(\mb{G}^{S})^{-\mb X}_y-Q((\mb{G}^{S})^{-\mb X},(\mb{G}^{S})^{-\mb X})\right] 
				\\[3mm]
				&\di ~~
				-\left( \f{u_1}{v}{\mb{L}}_\mb{M}^{-1}-\f{(u_1^{S})^{-\mb X}}{(v^{S})^{-\mb X}}{\mb{L}}_\mb{M^{S}}^{-1}\right) (\mb{G}^{S})^{-\mb X}_y\\[3mm]
				&\di ~~
				+\left( \f1v{\mb{L}}_\mb{M}^{-1}\mb{P}_1-\f{1}{(v^{S})^{-\mb X}}{\mb{L}}_\mb{M^S}^{-1}\mb{P}^{S}_1\right) 
				\left( \x_1(\mb{G}^{S})^{-\mb X}_y\right) .
			\end{array}
		\end{equation}
		Hereafter for brevity, we denote ${\mb{L}}_\mb{(\mb{M}^{S})^{-\mb X}}$, ${\mb{L}}_\mb{(\mb{M}^{S})^{-\mb X}}^{-1}$,  and $(\mb{P}^{S}_1)^{-\mb X}$ as ${\mb{L}}_\mb{M^{S}}$,  ${\mb{L}}_\mb{M^{S}}^{-1}$, and $\mb{P}^{S}_1$, respectively.
		
		First, using integration by parts, it holds that
		\begin{equation*}
			\begin{array}{ll}
				\di 
				K_1&\di= \int_0^T\int\left(  a^{-\mb X}\psi_1\right) _y\int \xi_1^2\left(\Pi_1-(\Pi_1^S)^{-\mb X} \right) d\xi dydt\\[4mm]
				&\di=\int_0^T\int \left(  a^{-\mb X}\psi_{1y}+\d_S^{-\f14}(v^S)_y^{-\mb X}\psi_1\right) \int \xi_1^2\left(\Pi_1-(\Pi_1^S)^{-\mb X} \right) d\xi dydt\\[4mm]
				&\di\leqslant\f1{200}\int_0^T\int\f43\frac{\mu(\bar \t)}{v} \psi_{1y}^2dydt+C\d_S\int_0^T\Lambda^Sdt\\[4mm]
				&\di~~~~~~~~~~~~+C\int_0^T\int\left(\int \xi_1^2 (\Pi_1-(\Pi_1^{S})^{-\mb X})d\xi\right) ^2dydt\\[4mm]
				&\di=:\f1{200}\int_0^T\int\f43\frac{\mu(\bar \t)}{v} \psi_{1y}^2dydt+C\d_S\int_0^T\Lambda^Sdt+C\sum_{i=1}^4K_{1}^i.
			\end{array}
		\end{equation*}
		Here,
		$$K_1^1=\int_0^T\int\left(\int \xi_1^2 {\mb{L}}_\mb{M}^{-1}\left[ \wt{\mb{G}}_t-\sigma \wt{\mb{G}}_y -\f{u_1}{v}\wt{\mb{G}}_{y}+\f{1}{v}\mb{P}_1(\x_1\wt{\mb{G}}_{y})-Q(\wt{\mb{G}},\wt{\mb{G}})\right]d\xi\right) ^2dydt,$$
		$$K_1^2=\int_0^T\int\left(\int \xi_1^2 {\mb{L}}_\mb{M}^{-1}\left[ Q(\wt{\mb{G}},(\mb{G}^{S})^{-\mb X})+ Q(({\mb{G}}^{S})^{-\mb X},\wt{\mb{G}})\right]d\xi\right) ^2dydt,$$
		$$K_1^3=\int_0^T\int\left(\int -\xi_1^2 \dot{\mb X}(t){\mb{L}}_\mb{M}^{-1}(\mb{G}^{S})^{-\mb X}_yd\xi\right) ^2dydt,~~~~~~~K_1^4=\int_0^T\int\left(\int \xi_1^2 J d\xi\right) ^2dydt.$$
		
		It follows from H{\"o}lder inequality and Lemma \ref{Lemma 4.3} that
		\begin{equation*}
			\begin{array}{ll}
				\di\int_0^T\int\left(\int \xi_1^2 {\mb{L}}_\mb{M}^{-1} \wt{\mb{G}}_td\xi\right) ^2dydt&\di\leqslant C\int_0^T\iint \frac{| {\mb{L}}_\mb{M}^{-1} \wt{\mb{G}}_t|^2}{\mb M_{_\#}}d\xi dydt\\[3mm]
				&\di\leqslant C\int_0^T\iint \frac{ (1+|\xi|)^{-1}|\wt{\mb{G}}_t|^2}{\mb M_{_\#}}d\xi dydt.
			\end{array}
		\end{equation*}
		By Lemma \ref{Lemma 4.1},
		\begin{equation}\label{Qestimate}
			\begin{array}{ll}
				\quad\di\int_0^T\int\left(\int \xi_1^2 {\mb{L}}_\mb{M}^{-1} Q(\wt{\mb{G}},\wt{\mb{G}})d\xi\right) ^2dydt\\[3mm]
				\di\leqslant C\int_0^T\iint \frac{ (1+|\xi|)^{-1}Q^2(\wt{\mb{G}},\wt{\mb{G}})}{\mb M_{_\#}}d\xi dydt\\[3mm]
				\di\leqslant C\int_0^T\iint \frac{ (1+|\xi|)|\wt{\mb{G}}|^2}{\mb M_{_\#}}d\xi\cdot\int\frac{|\wt{\mb{G}}|^2}{\mb M_{_\#}} d\xi dydt\\[3mm]
				\di\leqslant C\int_0^T\iint \frac{ (1+|\xi|)(|\wt{\mb{G}}_0|^2+|\wt{\mb{G}}_1|^2)}{\mb M_{_\#}}d\xi\cdot\int\frac{|\wt{\mb{G}}_0|^2+|\wt{\mb{G}}_1|^2}{\mb M_{_\#}} d\xi dydt\\[3mm]
				\di\leqslant C(\d_1+\ve_1)\int_0^T\iint \frac{ (1+|\xi|)|\wt{\mb{G}}_1|^2}{\mb M_{_\#}}d\xi dydt+C\d_1.
			\end{array}
		\end{equation}
		The other terms in $K_1^1$ can be estimated similarly, and we conclude
		\begin{equation*}
			\di K_1^1\di\leqslant C(\d_1+\ve_1)\int_0^T\iint \frac{ (1+|\xi|)|\wt{\mb{G}}_1|^2}{\mb M_{_\#}}d\xi dydt+C\int_0^T\iint \frac{ (1+|\xi|)(|\wt{\mb{G}}_t|^2+|\wt{\mb{G}}_y|^2)}{\mb M_{_\#}}d\xi dydt+C\d_1.
		\end{equation*}
		For $K_1^2$ and $K_1^3$, it holds that
		\begin{equation*}
			\begin{array}{ll}
				\di 
				K_1^2&\di\leqslant C\int_0^T\iint \frac{ (1+|\xi|)^{-1}\left[ Q^2(\wt{\mb{G}},(\mb{G}^{S})^{-\mb X})+Q^2((\mb{G}^{S})^{-\mb X},\wt{\mb{G}})\right] }{\mb M_{_\#}}d\xi dydt\\[4mm]
				&\di\leqslant C\int_0^T\iint \frac{ (1+|\xi|)\left( |\wt{\mb{G}}_0|^2+|\wt{\mb{G}}_1|^2\right)} {\mb M_{_\#}}\cdot\int \frac{ (1+|\xi|){|(\wt{\mb{G}}^{S})^{-\mb X}|^2}}{\mb M_{_\#}}d\xi dydt\\[4mm]
				&\di\leqslant C(\d_1+\ve_1)\int_0^T\iint \frac{ (1+|\xi|)|\wt{\mb{G}}_1|^2}{\mb M_{_\#}}d\xi dydt+C\d_1.
			\end{array}
		\end{equation*}
		$$K_1^3\leqslant C\int_0^T |\dot{\mb X}(t)|^2\iint \frac{ (1+|\xi|)|(\mb{G}^{S})^{-\mb X}_y|^2}{\mb M_{_\#}}d\xi dydt\leqslant C\d_S^5\int_0^T |\dot{\mb X}(t)|^2dt.$$
		By \eqref{J3}, we have
		\begin{equation*}
			\begin{array}{ll}
				\di 
				K_1^4&\di\leqslant C\int_0^T\int\left(\int \xi_1^2 \left( {\mb{L}}_\mb{M}^{-1}-{\mb{L}}_\mb{M^{S}}^{-1}\right) (\mb{G}^{S})^{-\mb X}_y d\xi\right) ^2dydt\\[4mm]
				&\di\quad+C\int_0^T\int\left(\int \xi_1^2 \left( {\mb{L}}_\mb{M}^{-1}-{\mb{L}}_\mb{M^{S}}^{-1}\right) Q((\mb{G}^{S})^{-\mb X},(\mb{G}^{S})^{-\mb X}) d\xi\right) ^2dydt\\[4mm]
				&\di\quad+C\int_0^T\int\left(\int \xi_1^2 \left( \f{u_1}{v}{\mb{L}}_\mb{M}^{-1}-\f{(u_1^{S})^{-\mb X}}{(v^{S})^{-\mb X}}{\mb{L}}_\mb{M^{S}}^{-1}\right) (\mb{G}^{S})^{-\mb X}_y d\xi\right) ^2dydt\\[4mm]
				&\di\quad+C\int_0^T\int\left(\int \xi_1^2 \left( \f1v{\mb{L}}_\mb{M}^{-1}\mb{P}_1-\f{1}{(v^{S})^{-\mb X}}{\mb{L}}_\mb{M^S}^{-1}\mb{P}^{S}_1\right) 
				\left( \x_1(\mb{G}^{S})^{-\mb X}_y\right)  d\xi\right) ^2dydt\\[4mm]
				&\di :=\sum_{i=1}^4 K_1^{4i}.
			\end{array}
		\end{equation*}
		For $K_1^{41}$, it holds that
		\begin{equation*}
			\begin{array}{ll}
				\quad\di\int_0^T\int\left(\int \xi_1^2 \left( {\mb{L}}_\mb{M}^{-1}-{\mb{L}}_\mb{M^{S}}^{-1}\right) (\mb{G}^{S})^{-\mb X}_y d\xi\right) ^2dydt\\[3mm]
				\di= \int_0^T\int\left(\int\xi_1^2{\mb{L}}_\mb{M}^{-1}\bigg[Q\big((\mb{M}^{S})^{-\mb X}-\mb{M},{\mb{L}}_\mb{M^{S}}^{-1}
				(\mb{G}^{S})^{-\mb X}_y\big)\right.\\[3mm]
				\di~~~~~~~~~~~~~~~~~~~~~\left.+Q\big({\mb{L}}_\mb{M^{S}}^{-1}
				(\mb{G}^{S})^{-\mb X}_y,(\mb{M}^{S})^{-\mb X}-\mb{M}\big)\bigg]d\xi\right) ^2dydt\\[3mm]
				\di\leqslant C\int_0^T\iint \frac{ (1+|\xi|)|(\mb{M}^{S})^{-\mb X}-\mb{M}|^2}{\mb M_{_\#}}d\xi\cdot\int\frac{(1+|\xi|)|{\mb{L}}_\mb{M^{S}}^{-1}
					(\mb{G}^{S})^{-\mb X}_y|^2}{\mb M_{_\#}} d\xi dydt\\[3mm]
				\di\leqslant C\delta_S\int_0^T\int |(v^{S})^{-\mb X}_y|^2| \left(v-(v^{S})^{-\mb X}, u-(u^{S})^{-\mb X},\t-(\t^{S})^{-\mb X} \right)|^2  dydt\\[3mm]
				\di\leqslant C\delta_S\int_0^T \Lambda^S dt+C\d_1.
			\end{array}
		\end{equation*}
		Similar estimates hold for the remaining
		terms $K_1^{4i}(i=2,3,4)$ and $K_j(j=2,3,4,5)$, and we conclude that
		\begin{equation}\label{(4.42)}
			\begin{array}{ll}
				\di 
				\sum_{i=1}^{5} K_i\leqslant C\d_S^5\int_0^T|\dot{\mb X}(t)|^2dt+C\d_1  \int_0^T(\Lambda^S+\Lambda^R)dt+C\d_1\\[4mm]
				\di\qquad\qquad
				+C\delta_C^2 \int_0^T(1+t)^{-1}\int  e^{-\frac{2c_0|y+\sigma t|^2}{1+t}}|(\phi,\z)|^2 dydt\\[4mm]
				\di\qquad\qquad
				+\f1{20}\int_0^T\int  \left(  \f43\frac{\mu(\bar \t)}{v} \psi_{1y}^2+\frac{\mu(\t)}{v}\sum_{i=2}^3\psi_{iy}^2+\frac{\k(\bar\t)}{v\t}\z_y^2 \right)   dydt\\[4mm]
				\di\qquad\qquad
				+C(\d_1+\ve_1)\int_0^T\iint \frac{ (1+|\xi|)|\wt{\mb{G}}_1|^2}{\mb M_{_\#}}d\xi dydt\\[4mm]
				\di\qquad\qquad
				+C\int_0^T\iint \frac{ (1+|\xi|)(|\wt{\mb{G}}_t|^2+|\wt{\mb{G}}_y|^2)}{\mb M_{_\#}}d\xi dydt.
			\end{array}
		\end{equation}

		{\it \underline{Step 5}. Estimation on the microscopic component $\wt{\mb G}_1$.} The microscopic component $\wt{\mb{G}}_1$ can be estimated by using \eqref{G1eq}. Multiplying \eqref{G1eq} by $\frac{\wt{\mb{G}}_1}{\mb{M}_{_\#}}$ gives
		\begin{equation}\label{M.1}
			\begin{array}{ll}
				&\di\left( \frac{|\wt{\mb{G}}_1|^2}{2\mb{M}_{_\#}}\right) _{t}-\f{\wt{\mb{G}}_1}{\mb{M}_{_\#}}\mb{L}_\mb{M}\wt{\mb{G}}_1\\[3mm]
				=&\di \bigg\{\sigma \wt{\mb{G}}_y+\dot{\mb X}(t)(\mb{G}^S)^{-\mb X}_y+\frac{u_1}{v} \wt{\mb{G}}_y+\left(\frac{u_1}{v}-\frac{(u_1^S)^{-\mb X}}{(v^S)^{-\mb X}} \right) (\mb{G}^S)^{-\mb X}_y
				-\frac{1}{v}\mb{P}_1(\xi_1\wt{\mb{G}}_{y})\\[3mm]
				&\di-\left(\frac{1}{v}-\frac{1}{(v^S)^{-\mb X}} \right)\mb{P}_1(\xi_1(\mb{G}^S)^{-\mb X}_y)-\frac{1}{(v^S)^{-\mb X}}  \left( \mb{P}_1-(\mb{P}^{S}_1)^{-\mb X}\right) (\x_1(\mb{G}^{S})^{-\mb X}_y)\\[3mm]
				&\di-\frac{3}{2v\t}\mb{P}_1\left[\xi_1\mb{M}\left( \xi\cdot\psi_y+\xi_1(u^S_{1})^{-\mb X }_y+\frac{\left|\xi-u \right|^2 }{2\t}(\z_{y}+(\t^S)^{-\mb X}_{y})\right) \right]\\[3mm]
				&\di+\f1{(v^S)^{-\mb X}}(\mb{P}^{S}_1)^{-\mb X}(\xi_1(\mb{M}^S)^{-\mb X}_y)
				+Q(\wt{\mb{G}},(\mb{G}^S)^{-\mb X})+Q(({\mb{G}}^S)^{-\mb X},\wt{\mb{G}})\\[3mm]
				&\di+Q(\wt{\mb{G}},\wt{\mb{G}})+(\mb{L}_{\mb{M}}-\mb{L}_{\mb{M}^S})(\mb{G}^S)^{-\mb X}-\wt{\mb{G}}_{0t}\bigg\}\f{\wt{\mb{G}}_1}{\mb{M}_{_\#}}.\\
			\end{array}
		\end{equation}
		The Cauchy inequality implies
		\begin{equation*}
			\begin{array}{ll}
				\quad\di\int_0^T\iint
				\dot{\mb X}(t)(\mb{G}^S)^{-\mb X}_y\f{\wt{\mb{G}}_1}{\mb{M}_{_\#}}d\xi
				dydt\\[3mm]
				\di\leqslant \f{\wt\s}{32}\int_{0}^T\iint\f{(1+|\xi|)|\wt{\mb{G}}_1|^2}{\mb{M}_{_\#}}d\xi
				dydt+C\int_{0}^T|\dot{\mb X}(t)|^2\iint\f{(1+|\xi|)|(\mb{G}^S)^{-\mb X}_y|^2}{\mb{M}_{_\#}}d\xi
				dydt\\[3mm]
				\di\leqslant
				\f{\wt\s}{32}\int_{0}^T\iint\f{(1+|\xi|)|\wt{\mb{G}}_1|^2}{\mb{M}_{_\#}}d\xi
				dydt+C\d_S^5\int_{0}^T|\dot{\mb X}(t)|^2dt
			\end{array}
		\end{equation*}
		and
		\begin{equation*}
			\begin{array}{ll}
				\quad\di\int_0^T\iint
				\left(\frac{u_1}{v}-\frac{(u_1^S)^{-\mb X}}{(v^S)^{-\mb X}} \right)(\mb{G}^S)^{-\mb X}_y\f{\wt{\mb{G}}_1}{\mb{M}_{_\#}}d\xi
				dydt\\[3mm]
				\di\leqslant \f{\wt\s}{32}\int_{0}^T\iint\f{(1+|\xi|)|\wt{\mb{G}}_1|^2}{\mb{M}_{_\#}}d\xi
				dydt+C\d_S\int_{0}^T\int |(v^S)^{-\mb X}_y|^2\left| \left( u_1-(u_1^S)^{-\mb X},v-(v^S)^{-\mb X}\right)\right| ^2
				dydt\\[3mm]
				\di\leqslant
				\f{\wt\s}{32}\int_{0}^T\iint\f{(1+|\xi|)|\wt{\mb{G}}_1|^2}{\mb{M}_{_\#}}d\xi
				dydt+C\d_S\int_{0}^T\Lambda^Sdt+C\d_1.
			\end{array}
		\end{equation*}
		Meanwhile, notice that
		$\di\mb{P}_1(\xi_1\wt{\mb{G}}_y)=\xi_1\wt{\mb{G}}_y-\sum_{j=0}^{4}\langle\xi_1\wt{\mb{G}}_y,\chi_j\rangle\chi_j$, and we have
		\begin{equation*}
			\begin{array}{ll}
				\quad\di\int_0^T\iint
				-\f1v\mb{P}_1(\xi_1\wt{\mb{G}}_y)\f{\wt{\mb{G}}_1}{\mb{M}_{_\#}}d\xi
				dydt\\[3mm]
				\di= \int_0^T\iint -\f1v\xi_1\wt{\mb{G}}_{y}\f{\wt{\mb{G}}_1}{\mb{M}_{_\#}}d \xi
				dydt
				+\sum_{j=0}^{4}\int_0^T\iint\f1v\langle \xi_1\wt{\mb{G}}_y,\chi_j\rangle\chi_j\f{\wt{\mb{G}}_1}{\mb{M}_{_\#}}d \xi dydt\\[3mm]
				\di\leqslant \f{\wt\s}{32}\int_{0}^T\iint\f{(1+|\xi|)|\wt{\mb{G}}_1|^2}{\mb{M}_{_\#}}d\xi
				dydt+C\int_{0}^T\iint\f{(1+|\xi|)|\wt{\mb{G}}_y|^2}{\mb{M}_{_\#}}d\xi
				dydt.
			\end{array}
		\end{equation*}
		We also have
		\begin{equation*}
			\begin{array}{ll}
				\quad\di\int_0^T\iint
				\left\lbrace -\frac{3}{2v\t}\mb{P}_1\left[\xi_1\mb{M}\left( \xi\cdot\psi_y+\xi_1(u^S_{1})^{-\mb X }_y+\frac{\left|\xi-u \right|^2 }{2\t}(\z_{y}+(\t^S)^{-\mb X}_{y})\right) \right]\right.\\[3mm]
				\di~~~~~~~~~~~~~~~~~~~~~~~~~~+\f1{(v^S)^{-\mb X}}(\mb{P}^{S}_1)^{-\mb X}(\xi_1(\mb{M}^S)^{-\mb X}_y) \Bigg\} \f{\wt{\mb{G}}_1}{\mb{M}_{_\#}}d\xi
				dydt\\[3mm]
				\di\leqslant
				\f{\wt\s}{32}\int_{0}^T\iint\f{(1+|\xi|)|\wt{\mb{G}}_1|^2}{\mb{M}_{_\#}}d\xi
				dydt+C\int_{0}^T\|(\psi_y,\z_y)\|_{L^2}^2dt+C\d_S\int_{0}^T\Lambda^Sdt+C\d_1.
			\end{array}
		\end{equation*}
		The other terms can be estimated similarly. Thus, integrating \eqref{M.1} with respect to $\xi$, $y$ and $t$ and using Lemma \ref{Lemma 4.1}, Lemma \ref{Lemma 4.3} and Cauchy inequality yield that
		\begin{equation}\label{M.10}
			\begin{array}{ll}
				\di \quad\iint\f{|\wt{\mb{G}}_1|^2}{\mb{M}_{_\#}}d\xi
				dy+\int_0^T\iint\f{(1+|\xi|)|\wt{\mb{G}} _1|^2}{\mb{M}_{_\#}}d\xi dydt\\[4mm]
				\di \leqslant C\mathcal{N}(0)^2+C\d_1\int_0^T\Lambda^Sdt+C\d_S^5\int_0^T|\dot{\mb X}(t)|^2dt+C\int_0^T\|(\psi_y,\z_y)\|_{L^2}^2dt\\[4mm]
				\di\qquad
				+C\d_1\int_0^T\|(\phi_t,\psi_t,\z_t)\|_{L^2}^2dt+C\int_0^T\iint\f{(1+|\xi|)|\wt{\mb{G}}_y|^2}{\mb{M}_{_\#}}d\xi dydt+C\d_1.
			\end{array}
		\end{equation}

		{\it \underline{Step 6}. Estimation on $\|\phi_y\|_{L^2}^2$.} In this step, we would like to recover the term $\di \|\p_y\|_{L^2}^2$ in the dissipation rate. Differentiating $\eqref{pereq1}_1$ with respect to $y$ and multiplying the result by $\f43\mu(\bar \t)\phi_y$, we have
		\begin{equation}\label{psiphire3}
			\f43\mu(\bar \t)\left( \frac{\phi_y^2}{2}\right) _t -\f43\s\mu(\bar \t)\left( \frac{\phi_y^2}{2}\right) _y - \dot{\mb X}(t) (v^S)^{-\mb X}_{yy}\f43\mu(\bar \t)\phi_y = \f43\mu(\bar \t)\phi_y\psi_{1yy}.
		\end{equation}
		Multiplying $\eqref{pereq1}_2$  by $-v\phi_{y}$ yields that
		\begin{equation}\label{psiphire4}
			\begin{array}{ll}
				\di -v\phi_{y}\psi_{1t}+\s v\phi_{y} \psi_{1y} + v\phi_{y} \dot{\mb{X}}(t)(u_1^S)^{-\mb{X}}_y
				\\[3mm]
				\di \quad =v\phi_{y} (p-\bar p)_y- v\phi_{y} \f43\left(\frac{\mu(\t)u_{1y}}{v} - \frac{\mu(\bar \t)\bar u_{1y}}{\bar v}\right)_y +v\phi_{y} Q_1+v\p_y\int \xi_1^2\left( \Pi_1-(\Pi_1^S)^{-\mb X}\right)_y d\xi.
			\end{array}
		\end{equation}
		Then, adding \eqref{psiphire3}  to \eqref{psiphire4} yields that
		\begin{equation}\label{(5.4)}
			\begin{array}{ll}
				\di 
				\left( \f23\mu(\bar \t) \phi_y^2-v\psi_1\phi_y\right) _t+ \left( -\f23\s\mu(\bar\t) \phi_y^2+\s v\psi_1\phi_y+v\psi_1(\phi_t-\s\phi_y)\right) _y+\frac{2\theta\phi_y^2}{3v}+(v_t-\s v_y)\psi_1\phi_y\\[3mm]
				\di\quad - \big(v_y\psi_1+v\psi_{1y}) \big(\phi_t-\s\phi_y\big)+\dot{\mb{X}}(t)\phi_y\left( v(u_1^S)^{-\mb{X}}_y-\f43\mu(\bar\t)(v^S)^{-\mb{X}}_{yy}\right) \\[3mm]
				\di  =\f23\z_y\phi_y+\f23v\phi_y\bar\theta_y\left(\frac 1v-\frac 1{\bar v}\right)+\f23v\phi_y\bar v_y \left(\frac{\bar\theta}{\bar v^2}-\frac{\theta}{v^2}\right) -\f43\mu(\bar\t) v\phi_y\bar u_{1yy}\left(\frac 1v-\frac 1{\bar v}\right)\\[3mm]
				\di\quad  -\f43\mu(\bar\t) v\phi_y \bar u_{1y} \left(\frac{\bar v_y}{\bar v^2} - \frac{v_y}{v^2} \right) + \f43\frac{\mu(\bar\t)}{v}\phi_y \psi_{1y}v_y +v\phi_y Q_1+v\phi_y\int \xi_1^2\left( \Pi_1-(\Pi_1^S)^{-\mb X}\right)_y d\xi\\[3mm]
				\di\quad
				+\f23\mu'(\bar\t)\bar\t_t\phi_y^2-\f23\s\mu'(\bar\t)\bar\t_y\phi_y^2-\f43u_{1y}\phi_y\left(\mu(\t)-\mu(\bar\t) \right) _y-\f43v\phi_y\left(\mu(\t)-\mu(\bar\t) \right)\left(\frac{u_{1y}}{v} \right)_y\\[3mm]
				\di\quad
				-\f43\phi_y\mu'(\bar\t)\bar\t_y\left( \psi_{1y}-\frac{\bar u_{1y}}{\bar v}\p\right). 
			\end{array}
		\end{equation}
		Similar to \eqref{(4.42)}, by using Lemma \ref{Lemma 4.1} and Lemma \ref{Lemma 4.3}, it holds that
		\begin{equation}\label{Ax}
			\begin{array}{ll}
				\di \int_0^T\int \left|\int \xi_1^2(\Pi-(\Pi^{S}_1)^{-\mb X})_yd \x\right|^2 dydt\\[3mm]
				\di\leqslant 
				C\d_1\int_0^T\Lambda^S dt+C\d_1\int_0^T\|(\phi_y,\psi_y,\z_y)\|_{L^2}^2 dt+C\d_S^2\int_0^T|\dot{\mb X}(t) |^2 dt+C\d_1\\[3mm]
				\di\quad +C(\d_1+\ve_1)\int_0^T\iint\f{(1+|\xi|)(|\wt{\mb{G}}_1|^2+|\wt{\mb{G}}_y|^2+|\wt{\mb{G}}_t|^2)}{\mb{M}_{_\#}}
				d \x dydt\\[3mm]
				\di\quad
				+C\int_0^T \iint\f{(1+|\xi|)(|\wt{\mb{G}}_{yy}|^2+|\wt{\mb{G}}_{yt}|^2)}{\mb{M}_{_\#}}d \xi dydt.
			\end{array}
		\end{equation}
		Then, integrating \eqref{(5.4)} with respect to $y$, $t$, together with \eqref{Ax}, we have
		\begin{equation}\label{1orderphi}
			\begin{array}{ll}
				\di \sup_{t\in[0,T]}\|\p_y\|_{L^2}^2+\int_{0}^T\|\p_y\|_{L^2}^2dt\leqslant C\mathcal{N}(0)^2+C\sup_{t\in[0,T]}\|\psi_1\|_{L^2}^2+C\delta_S^2\int_0^T|\dot{\mb X}(t)|^2 dt\\[3mm]
				\di\quad
				+C\d_1\int_0^T(\Lambda^S+\Lambda^R) dt+C\delta_C^2 \int_0^T(1+t)^{-1}\int  e^{-\frac{2c_0|y+\sigma t|^2}{1+t}}|(\phi,\psi,\z)|^2 dydt+C\delta_1^{\f12}\\[3mm]
				\di\quad +C\int_0^T\|(\psi_y,\z_y)\|_{L^2}^2 dt+C\ve_1^2\int_0^T\|\psi_{1yy}\|_{L^2}^2dt+C\int_0^T \iint\f{(1+|\xi|)(|\wt{\mb{G}}_{yy}|^2+|\wt{\mb{G}}_{yt}|^2)}{\mb{M}_{_\#}}d \xi dydt\\[3mm]
				\di\quad +C(\d_1+\ve_1)\int_0^T\iint\f{(1+|\xi|)(|\wt{\mb{G}}_1|^2+|\wt{\mb{G}}_y|^2+|\wt{\mb{G}}_t|^2)}{\mb{M}_{_\#}}
				d \x dydt.
			\end{array}
		\end{equation}
		
		We would like to use this occasion to perform the estimation for the time derivatives $\|(\p_{t},\psi_{t},\z_{t})\|_{L^2}^2$, as it is similar to the estimates for $\|(\p_{y},\psi_{y},\z_{y})\|_{L^2}^2$. By multiplying $(\ref{(4.18)})_1$ by $\p_{t}$,
		$(\ref{(4.18)})_2$ by $\psi_{1t}$, $(\ref{(4.18)})_3$ by
		$\psi_{it}~(i=2,3)$, and $(\ref{(4.18)})_4$ by $\z_{t}$
		respectively, and adding them together, after integrating with
		respect to $y$ and $t$, we have
		\begin{equation}\label{1ordert}
			\begin{array}{l}
				\di \int_0^T\|(\p_{t},\psi_{t},\z_{t})\|_{L^2}^2dt\leqslant C\delta_S^2\int_0^T|\dot{\mb X}(t)|^2 dt+
				C\d_1\int_0^T(\Lambda^S+\Lambda^R) dt\\[3mm]
				\di\quad
				+C\delta_C^2 \int_0^T(1+t)^{-1}\int  e^{-\frac{2c_0|y+\sigma t|^2}{1+t}}|(\phi,\psi,\z)|^2 dydt+C\d_1\\[3mm]
				\di\quad
				+C\int_0^T\|(\p_y,\psi_y,\z_y)\|_{L^2}^2dt+C\int_0^T\iint\f{(1+|\xi|)|\wt{\mb{G}}_y|^2}{\mb{M}_{_\#}}d\xi
				dydt.
			\end{array}
		\end{equation}

		{\it \underline{Step 7}. Conclusion.} To conclude, combining \eqref{weightedre}, \eqref{4.2conclusion}, \eqref{(4.36)}, \eqref{(4.42)}, and \eqref{M.10}, with the smallness of $\d_1$ and $\ve_1$ and the fact that
		$$\Phi(z)\sim (z-1)^2,~~|z-1|\ll 1,$$
		we have
		\begin{equation}\label{con1}
			\begin{array}{ll}
				\di \sup_{t\in[0,T]}\left[ \left\|(\phi,\psi,\z) \right\| _{L^2}^2+\iint\f{|\wt{\mb{G}}_1|^2}{\mb{M}_{_\#}}d \xi dy\right] +\delta_S\int_0^T| \dot{\mb{X}}(t)| ^2dt+\int_0^T\left( \Lambda^R+\Lambda^S\right) dt\\
				\di\qquad\qquad+\int_0^T\|( \psi_y,\z_y)\|_{L^2}^2 dt
				+\int_0^T\iint\f{(1+|\xi|)|\wt{\mb{G}}_1|^2}{\mb{M}_{_\#}}d \xi dydt\\[3mm]
				\di\quad\leqslant C\mathcal{N}(0)^2+C\d_1\int_0^T\|(\phi_t,\psi_t,\z_t)\|_{L^2}^2dt+C\int_0^T\iint\f{(1+|\xi|)(|\wt{\mb{G}}_y|^2+|\wt{\mb{G}}_t|^2)}{\mb{M}_{_\#}}
				d \x dydt\\[3mm]
				\di\qquad\qquad
				+C\d_C\int_0^T\f{1}{1+t} \int e^{-\frac{c_0 |y+\sigma t|^2}{1+t}} |(\phi,\z)|^2 dydt+C\d_1^{\f12}.
			\end{array}
		\end{equation}
		Furthermore, \eqref{1orderphi}  and \eqref{1ordert} yield that
		\begin{equation*}
			\begin{array}{ll}
				\di \sup_{t\in[0,T]}\|\p_y\|_{L^2}^2+\int_0^T\|(\phi_y,\p_{t},\psi_{t},\z_{t})\|_{L^2}^2dt\leqslant C\mathcal{N}(0)^2+C\sup_{t\in[0,T]}\|\psi_1\|_{L^2}^2+C\delta_S^2\int_0^T| \dot{\mb X}(t)| ^2 dt\\[3mm]
				\di\quad
				+C\d_1\int_0^T(\Lambda^S+\Lambda^R) dt+C\delta_C^2 \int_0^T(1+t)^{-1}\int  e^{-\frac{2c_0|y+\sigma t|^2}{1+t}}|(\phi,\psi,\z)|^2 dydt+C\delta_1^{\f12}\\[3mm]
				\di\quad +C\int_0^T\|(\psi_y,\z_y)\|_{L^2}^2 dt+C\ve_1^2\int_0^T\|\psi_{1yy}\|_{L^2}^2dt+C(\d_1+\ve_1)\int_0^T \iint\f{(1+|\xi|)|\wt{\mb{G}}_{1}|^2}{\mb{M}_{_\#}}d \xi dydt\\[3mm]
				\di\quad +C\int_0^T\iint\f{(1+|\xi|)(|\wt{\mb{G}}_y|^2+|\wt{\mb{G}}_t|^2+|\wt{\mb{G}}_{y
						y}|^2+|\wt{\mb{G}}_{yt}|^2)}{\mb{M}_{_\#}}
				d \x dydt.
			\end{array}
		\end{equation*}
		Thus, the combination of the above estimate and \eqref{con1}  completes the proof of Proposition \ref{lowerorderestimate}.

		\
		
		\section{Higher order estimates}
		\setcounter{equation}{0}
		
		This section is dedicated to the higher order estimates to \eqref{pereq1} and \eqref{tildeG}.
		\begin{proposition}\label{higherorderestimate}
			Under the assumptions of Proposition \ref{main proposition}, there exists a positive
			constant $C$ such that
			\begin{equation*}
				\begin{array}{ll}
					\di \sup_{t\in[0,T]}\left[ \left\|(\phi_y,\psi_y,\z_y) \right\| _{L^2}^2+\iint\f{\left( |\wt{\mb{G}}_y|^2+|\wt{\mb{G}}_t|^2+|\wt{f}_{yy}|^2+|\wt{f}_{yt}|^2\right) }{\mb{M}_{_\#}}d\xi dy\right] \\[2mm]
					\di \qquad
					+\int_0^T \|(\phi_{yy},\psi_{yy},\z_{yy}, \phi_{yt},\psi_{yt},\z_{yt})\|_{L^2}^2dt
					+\int_0^T\iint \f{(1+|\xi|)\big( | \wt{\mb{G}}_{yy}|^2+|\wt{\mb{G}}_{yt}
						|^2\big)}{\mb{M}_{_\#}}d\xi dydt\\[3mm]
					\di\quad\leqslant C\mathcal{N}(0)^2+C\delta_S^2\int_0^T| \dot{\mb{X}}(t)| ^2dt+C\d_1\int_0^T(\Lambda^S+\Lambda^R)dt\\[2mm]
					\di\qquad
					+C\delta_C^2 \int_0^T(1+t)^{-1}\int  e^{-\frac{2c_0|y+\sigma t|^2}{1+t}}|(\phi,\psi,\z)|^2 dydt
					+C(\d_1+\ve_1)\sum_{|\b|=1}\int_0^T\|\partial^\b(\p,\psi,\z)\|_{L^2}^2dt\\[2mm]
					\di\qquad
					+C(\d_1+\ve_1)\int_0^T\iint
					\f{(1+|\xi|)|\wt{\mb{G}}_1|^2}{\mb{M}_{_\#}}d\x
					dydt+C\d_1^{\f12}.
				\end{array}
			\end{equation*}
		\end{proposition}
		{\it Proof.} The proof is divided into the following three steps.

		{\it \underline{Step 1}. Estimation on	$\|(\phi_y,\psi_y,\z_y)\|_{L^2}^2$.} First, multiplying $(\ref{pereq1})_1$ by $-\phi_{yy}\frac{\bar p \t}{v}$, $(\ref{pereq1})_2$ by $-\psi_{1yy}\frac{3\bar p v}{2}$, $(\ref{pereq1})_3$ by $-\psi_{iyy}$, $(\ref{pereq1})_4$ by $-\z_{yy}$, and adding them together, we have
		\begin{equation}\label{(5.8)}
			\begin{array}{ll}
				\quad\di\left(\frac{\bar p \t}{2v}\phi_y^2+\frac{3\bar p v}{4}\psi_{1y}^2+\f12\sum_{i=2}^3\psi_{iy}^2+\f12 \z_y^2 \right)_t+2\bar p\mu(\bar \t)\psi_{1yy}^2+\frac{\mu(\t)}{v}\sum_{i=2}^3\psi_{iyy}^2+\frac{\k(\bar\t)}{v}\z_{yy}\\[3mm]
				\di=-\dot{\mb X}(t)\left[(v^S)^{-\mb X}_y \frac{\bar p \t}{v} \phi_{yy}+(u_1^S)^{-\mb X}_y \frac{3\bar p v}{2} \psi_{1yy}+(\t^S)^{-\mb X}_y\z_{yy} \right] \\[3mm]
				\di\quad
				+\left[\left( \frac{\bar p \t}{v}\right)_t+\s\left( \frac{\bar p \t}{v}\right)_y  \right] \frac{\phi_y^2}{2}+\left[\left( \frac{3\bar p v}{2}\right)_t+\s\left( \frac{3\bar p v}{2}\right)_y  \right] \frac{\psi_{1y}^2}{2}\\[4mm]
				\di\quad
				-\left(\frac{\bar p \t}{v} \right)_y\phi_t\phi_y+\left(\frac{\bar p \t}{v} \right)_y\psi_{1y}\phi_y-\left(\frac{3\bar p v}{2} \right)_y\psi_{1t}\psi_{1y}-\bar p_y\psi_{1y}\z_y  \\[4mm]
				\di\quad
				+\psi_{1yy}\bar p v \bar\t_y\left(\f1v-\f1{\bar v} \right)-\psi_{1yy}\bar p v \bar v_y\left(\frac{\t}{v^2}-\frac{\bar \t}{\bar v^2} \right) +(p-\bar p)u_{1y}\z_{yy}+\psi_{1yy}\frac{3\bar p v}{2}Q_1+\z_{yy}Q_2 \\[4mm]
				\di\quad
				-2\bar pv\psi_{1yy}\left((\mu(\t)-\mu(\bar \t))\frac{u_{1y}}{v} \right)_y-\z_{yy}\left((\k(\t)-\k(\bar \t))\frac{\t_{y}}{v} \right)_y \\[4mm]
				\di\quad
				-2\bar pv\psi_{1yy}\left[(\mu(\bar\t))_y\left(\frac{u_{1y}}{v}-\frac{\bar u_{1y}}{\bar v} \right)-\mu(\bar\t)\psi_{1y}\frac{v_y}{v^2}+ \mu(\bar\t)\left(\bar u_{1y}\left(\f1v-\f1{\bar v} \right) \right)_y  \right] \\[4mm]
				\di\quad
				-\z_{yy}\left[(\k(\bar\t))_y\left(\frac{\t_{y}}{v}-\frac{\bar\t_{y}}{\bar v} \right)-\k(\bar\t)\z_{y}\frac{v_y}{v^2}+ \k(\bar\t)\left(\bar \t_{y}\left(\f1v-\f1{\bar v} \right) \right)_y  \right] \\[4mm]
				\di\quad
				-\f43\z_{yy}\left( \mu(\t)\frac{u_{1y}^2}{v}-\mu(\bar\t) \frac{\bar u_{1y}^2}{\bar v} \right) -\z_{yy} \frac{\mu(\t)}{v}\sum_{i=2}^3\psi_{iy}^2\\[4mm]
				\di\quad
				+\psi_{1yy}\frac{3\bar p v}{2}\int \xi_1^2\left( \Pi_1-(\Pi_1^S)^{-\mb X}\right)_y d\xi+\psi_{iyy}\int \xi_1\xi_i\left( \Pi_1-(\Pi_1^S)^{-\mb X}\right)_y d\xi\\[4mm]
				\di\quad
				+\z_{yy}\int \xi_1\frac{|\xi|^2}{2}\left( \Pi_1-(\Pi_1^S)^{-\mb X}\right)_y d\xi-\z_{yy}\sum_{i=2}^3 \psi_i\int \xi_1\xi_i\Pi_{1y} d\xi\\[4mm]
				\di\quad
				-\z_{yy}\left[ u_1\int \xi_1^2\Pi_{1y} d\xi-(u_1^S)^{-\mb X}\int \xi_1^2(\Pi_1^S)^{-\mb X}_y d\xi\right] +\left(\cdots\right)_y,
			\end{array}
		\end{equation}
		where $(\cdots)_y$ represents the conservative terms which vanishes after integrating with respect to $y$ over $\mathbb R$. Integrating \eqref{(5.8)} with respect to $y$ and $t$, together with \eqref{Ax}, we have
		\begin{equation}\label{2orderyy}
			\begin{array}{ll}
				\di \sup_{t\in[0,T]}\|(\phi_y,\psi_{y},\z_y)\|_{L^2}^2+\int_{0}^T\|(\psi_{yy},\z_{yy})\|_{L^2}^2dt\leqslant C\mathcal{N}(0)^2+C\delta_S^2\int_0^T|\dot{\mb X}(t)|^2 dt\\[3mm]
				\di\quad
				+C\d_1\int_0^T(\Lambda^S+\Lambda^R) dt+C\delta_C^2 \int_0^T(1+t)^{-1}\int  e^{-\frac{2c_0|y+\sigma t|^2}{1+t}}|(\phi,\z)|^2 dydt\\[3mm]
				\di\quad +C(\d_1+\ve_1)\sum_{|\beta|=1}\int_0^T\|\partial^{\b}(\phi,\psi,\z)\|_{L^2}^2 dt+C\d_1\int_0^T\|\phi_{yy}\|_{L^2}^2dt+C\delta_1^{\f12}\\[3mm]
				\di\quad
				+C(\d_1+\ve_1)\int_0^T\iint\f{(1+|\xi|)(|\wt{\mb{G}}_1|^2+|\wt{\mb{G}}_y|^2+|\wt{\mb{G}}_t|^2)}{\mb{M}_{_\#}}
				d \x dydt\\[3mm]
				\di\quad
				+C\int_0^T \iint\f{(1+|\xi|)\left( |\wt{\mb{G}}_{yy} |^2+|\wt{\mb{G}}_{yt} |^2\right) }{\mb{M}_{_\#}}d \xi dydt.
			\end{array}
		\end{equation}
		
		Again, to recover  $\|\p_{yy}\|_{L^2}^2$ in the dissipation rate, applying
		$\partial_y$ to $(\ref{(4.18)})_2$, we get

		\begin{equation}\label{(4.50)}
			\psi_{1yt}-\s\psi_{1yy}-\dot{\mb X}(t)(u_1^S)^{-\mb X}_{yy}+(p-\bar p)_{yy}+(\bar p-p^R-p^C-(p^S)^{-\mb X})_{yy}
			=-(u^C_{1ty}-\s u^C_{1yy})-\int\xi_1^2\wt{\mb{G}}_{yy}d\xi. 
		\end{equation}
		Note that
		\begin{equation}\label{p-barp}
			(p-\bar p)_{yy}=-\f{p}{v}\p_{yy}+\f{2}{3v}\z_{yy}-\f1v(p-\bar p)\bar v_{yy}
			-\f{\p}{v}\bar p_{yy}-\f{2v_y}{v}(p-\bar p)_y-\f{2\bar p_y}{v}\p_y. 
		\end{equation}
		Multiplying (\ref{(4.50)}) by $-\p_{yy}$ and using \eqref{p-barp}
		imply
		\begin{equation}\label{(4.52)}
			\begin{array}{l}
				\di -\int\psi_{1y}\p_{yy}(T,y)dy+\int_0^T\int\f{p}{2v}\p_{yy}^2dydt\leqslant C\mathcal{N}(0)^2+C\d_1\int_0^T(\Lambda^S+\Lambda^R) dt+C\d_1\\[3mm]
				\di\qquad+C\delta_C^2 \int_0^T(1+t)^{-1}\int  e^{-\frac{2c_0|y+\sigma t|^2}{1+t}}|(\phi,\z)|^2 dydt+C(\d_1+\ve_1)\int_0
				^T\|(\p_y,\psi_y,\z_y)\|_{L^2}^2dt
				\\[3mm]
				\qquad\di+C\int_0^T\|(\psi_{1yy},\z_{yy})\|_{L^2}^2dt+C
				\int_0^T\iint\f{(1+|\xi|)|
					\wt{\mb{G}}_{yy} |^2}{\mb{M}_{_\#}}d\xi dydt+C\delta_S^2\int_0^T|\dot{\mb X}(t)|^2 dt.
			\end{array}
		\end{equation}
		
		To estimate $\|(\p_{yt},\psi_{yt},\z_{yt})\|_{L^2}^2$, we use the
		system (\ref{(4.18)}) again. Differentiating 
		(\ref{(4.18)}) with respect to $y$ and multiplying the four equations of (\ref{(4.18)})
		by $\p_{yt}$, $\psi_{1yt}$, $\psi_{iyt}$ $(i= 2,3)$,
		$\z_{yt}$ respectively, then adding them together and integrating
		with respect to  and $y$ and $t$, we have
		\begin{equation}\label{(4.53)}
			\begin{array}{l}
				\di \int_0^T\|(\p_{yt},\psi_{yt},\z_{yt})\|_{L^2}^2dt\leqslant
				C\int_0^T\|(\p_{yy},\psi_{yy},\z_{yy})\|_{L^2}^2dt+C\delta_S^2\int_0^T|\dot{\mb X}(t)|^2 dt\\
				\quad\di +C\d_1\int_0^T(\Lambda^S+\Lambda^R) dt+C\delta_C^2 \int_0^T(1+t)^{-1}\int  e^{-\frac{2c_0|y+\sigma t|^2}{1+t}}|(\phi,\z)|^2 dydt\\[3mm]
				\di \quad +C\delta_C^2 \int_0^T(1+t)^{-1}\int  e^{-\frac{2c_0|y+\sigma t|^2}{1+t}}|(\phi,\psi,\z)|^2 dydt\\[3mm]
				\di\quad +C(\d_1+\ve_1)\int_0
				^T\|(\p_y,\psi_y,\z_y)\|_{L^2}^2dt+C(\d_1+\ve_1)\int_0^T\iint\f{(1+|\xi|)|\wt{\mb{G}}_y|^2}{\mb{M}_{_\#}}
				d \x dydt\\[3mm]
				\quad\di+C\int_0^T\iint\f{(1+|\xi|)|\wt{\mb{G}}_{yy}|^2}{\mb{M}_{_\#}}
				d\xi dydt+C\d_1.
			\end{array}
		\end{equation}

		A suitable linear combination of \eqref{2orderyy}, \eqref{(4.52)}, and \eqref{(4.53)} gives
		\begin{equation}\label{Step2final}
			\begin{array}{ll}
				\di \sup_{t\in[0,T]}\|(\phi_y,\psi_{y},\z_y)\|^2
				+\int_0^T\left( \|(\phi_{yy},\psi_{yy},\z_{yy})\|_{L^2}^2+\|(\phi_{yt},\psi_{yt},\z_{yt})\|_{L^2}^2\right) dt\\[3mm]
				\di\leqslant C\mathcal{N}(0)^2+C\delta_S^2\int_0^T|\dot{\mb X}(t)|^2 dt+C(\d_1+\ve_1)\sum_{|\beta|=1}\int_0^T\|\partial^{\b}(\phi,\psi,\z)\|_{L^2}^2 dt\\[3mm]
				\di\quad
				+C\d_1\int_0^T(\Lambda^S+\Lambda^R) dt+C\delta_C^2 \int_0^T(1+t)^{-1}\int  e^{-\frac{2c_0|y+\sigma t|^2}{1+t}}|(\phi,\z)|^2 dydt \\[3mm]
				\di\quad +C(\d_1+\ve_1)\int_0^T\iint\f{(1+|\xi|)(|\wt{\mb{G}}_1|^2+|\wt{\mb{G}}_y|^2+|\wt{\mb{G}}_t|^2)}{\mb{M}_{_\#}}
				d \x dydt\\[3mm]
				\di\quad +C\int_0^T \iint\f{(1+|\xi|)\left( |\wt{\mb{G}}_{yy} |^2+|\wt{\mb{G}}_{yt} |^2\right) }{\mb{M}_{_\#}}d \xi dydt+C\delta_1^{\f12}+C\sup_{t\in[0,T]}\left|\int\psi_{1y}\phi_{yy}
				dy \right|.
			\end{array}
		\end{equation}

		{\it \underline{Step 2}. Estimation on the microscopic component.} To close the a priori estimate, we also need to estimate the
		derivatives on the non-fluid component $\wt{\mb{G}} $, i.e.,
		$\partial^\a \wt{\mb{G}}, (|\a|=1,2)$. For this, applying $\partial_y$ on \eqref{tildeG} and multiplying the result by $\f{\wt{\mb{G}}_y}{\mb{M}_{_\#}}$, we have
		\begin{equation}\label{(5.18)}
			\begin{array}{ll}
				\di\quad\left( \frac{|\wt{\mb{G}}_{y}|^2}{2\mb{M}_{_\#}}\right)_t -\f{\wt{\mb{G}}_y}{\mb{M}_{_\#}}(\mb{L}_\mb{M}\wt{\mb{G}})_y\\[3mm]
				=\di \bigg[\sigma \wt{\mb{G}}_y+\dot{\mb X}(t)(\mb{G}^S)^{-\mb X}_y+\frac{u_1}{v} \wt{\mb{G}}_y-\frac{1}{v}\mb{P}_1(\xi_1\wt{\mb{G}}_{y})+\left(\frac{u_1}{v}-\frac{(u_1^S)^{-\mb X}}{(v^S)^{-\mb X}} \right) (\mb{G}^S)^{-\mb X}_y
				\\[3mm]
				\di\quad-\left( \frac{1}{v}\mb{P}_1(\xi_1(\mb{G}^S)^{-\mb X}_y)-\frac{1}{(v^S)^{-\mb X}}   (\mb{P}^{S}_1)^{-\mb X} (\x_1(\mb{G}^{S})^{-\mb X}_y)\right) \\[3mm]
				\di\quad-\f1v\mb{P}_1(\xi_1\mb{M}_y)+\f1{(v^S)^{-\mb X}}(\mb{P}^{S}_1)^{-\mb X}(\xi_1(\mb{M}^S)^{-\mb X}_y)
				+Q(\wt{\mb{G}},\wt{\mb{G}})
				\\[3mm]
				\di\quad+Q(\wt{\mb{G}},(\mb{G}^S)^{-\mb X})+Q((\mb{G}^S)^{-\mb X},\wt{\mb{G}})+(\mb{L}_{\mb{M}}-\mb{L}_{\mb{M}^S})(\mb{G}^S)^{-\mb X}\bigg]_y\f{\wt{\mb{G}}_y}{\mb{M}_{_\#}}.\\
			\end{array}
		\end{equation}
		Integrating \eqref{(5.18)}, together with Lemma \ref{Lemma 4.1}, Lemma \ref{Lemma 4.3} and the notation that
		\begin{equation*}
			\begin{array}{l}
				(\mb{L}_\mb{M}g)_{\varsigma}=\mb{L}_\mb{M}(g_{\varsigma})+Q(g,\mb{M}_{\varsigma})+Q(\mb{M}_{\varsigma},g),
				\quad{\rm for}~~ \varsigma=t,y,
			\end{array}
		\end{equation*}
		we deduce
		\begin{equation}\label{Gy}
			\begin{array}{ll}
				\di \sup_{t\in[0,T]}\iint\f{|\wt{\mb{G}} _y|^2}{\mb{M}_{_\#}}d\xi
				dy+\int_0^T\iint\f{(1+|\xi|)|\wt{\mb{G}} _y|^2}{\mb{M}_{_\#}}d\xi dydt\\[3mm]
				\di \leqslant C\mathcal{N}(0)^2+C\delta_S^2\int_0^T|\dot{\mb X}(t)|^2 dt+C\d_1\int_0^T\Lambda^S dt+C\d_1\\[3mm]
				\di\quad +C(\d_1+\ve_1)\int_0^T\|(\p_y,\psi_y,\z_y)\|_{L^2}^2dt+C\int_0^T\|(\psi_{yy},\z_{yy})\|_{L^2}^2dt\\[3mm]
				\di\quad
				+C(\d_1+\ve_1)\int_0^T\iint\f{(1+|\xi|)|\wt{\mb{G}}_1|^2}{\mb{M}_{_\#}}
				d \x dydt+C\int_0^T \iint\f{(1+|\xi|)|
					\wt{\mb{G}}_{yy} |^2}{\mb{M}_{_\#}}d \xi dydt.
			\end{array}
		\end{equation}
		
		Before taking the derivative of $t$, we should notice from \eqref{X(t)} and H{\"o}lder inequality that
		\begin{equation}\label{Xdot2}
			\begin{array}{ll}
				\di \quad\int_0^T |\ddot{\mb{X}}(t)|^2 \iint \frac{|(\mb{G}^S)^{-\mb X}_y|^2}{\mb M_{_\#}} d\xi dydt\\[3mm]
				\di\leqslant C\d_S^2\int_0^T |\ddot{\mb{X}}(t)|^2 \int |(v^S)^{-\mb X}_y|^2 dydt \\[3mm]
				\di\leqslant C\d_S^5\int_0^T |\ddot{\mb{X}}(t)|^2dt \\[3mm]
				\di\leqslant C\d_S^3\int_0^T \left[\int\frac{\partial}{\partial t}\left(  a^{-\mb{X}}\left((v^S)^{-\mb{X}}_{y}\frac{\bar p }{\bar v}(v-\bar v)+(u^S_{1})^{-\mb{X}}_y(u_1-\bar u_1)+\frac{(\t^S)^{-\mb{X}}_{y}}{\bar\theta}(\t-\bar\t)\right) \right) dy\right] ^2dt\\[3mm]
				\di\leqslant C\d_S\int_0^T \|(\phi_t,\psi_{1t},\z_t)\|_{L^2}^2dt+C\d_S\int_0^T\Lambda^Sdt.\\[3mm]
			\end{array}
		\end{equation}
		Therefore, applying $\partial_t$ on \eqref{tildeG} and multiplying the result by $\f{\wt{\mb{G}}_t}{\mb{M}_{_\#}}$, one can get
		\begin{equation}\label{Gt}
			\begin{array}{ll}
				\di \sup_{t\in[0,T]}\iint\f{|\wt{\mb{G}} _t|^2}{\mb{M}_{_\#}}d\xi
				dy+\int_0^T\iint\f{(1+|\xi|)|\wt{\mb{G}} _t|^2}{\mb{M}_{_\#}}d\xi dydt\\[3mm]
				\di \leqslant C\mathcal{N}(0)^2+C\delta_S^2\int_0^T|\dot{\mb X}(t)|^2 dt+C\d_1\int_0^T\Lambda^S dt+C\d_1\\[3mm]
				\di\quad +C(\d_1+\ve_1)\sum_{|\beta|=1}\int_0^T\|\partial^{\b}(\p,\psi,\z)\|_{L^2}^2dt+C\int_0^T\|(\psi_{yt},\z_{yt})\|_{L^2}^2dt\\[3mm]
				\di\quad
				+C(\d_1+\ve_1)\int_0^T\iint\f{(1+|\xi|)\left( |\wt{\mb{G}}_1|^2+|\wt{\mb{G}}_y|^2\right) }{\mb{M}_{_\#}}
				d \x dydt\\[3mm]
				\di\quad+C\int_0^T \iint\f{(1+|\xi|)|
					\wt{\mb{G}}_{yt} |^2}{\mb{M}_{_\#}}d \xi dydt.
			\end{array}
		\end{equation}

		{\it \underline{Step 3}: Highest order estimates.} Finally, we estimate the highest
		order derivatives, that is, $\di \int\psi_{1y}\p_{yy}dy$ and $\di
		\int_0^T\iint \f{(1+|\xi|)\left( | \wt{\mb{G}}_{yy}|^2+|\wt{\mb{G}}_{yt}
			|^2\right)}{\mb{M}_{_\#}}d\xi dydt$. To
		do so, it is sufficient to study  $\di \iint \f{|
			\wt{f}_{yy}|^2+|
			\wt{f}_{yt}|^2}{\mb{M}_{_\#}}d\xi dy$ in view of \eqref{vandf}. Using the same idea in \cite{Huang-Wang-Yang-2010-1} and noticing \eqref{Xdot2},
		we obtain the estimation for the highest order derivative terms,
		i.e.,
		\begin{equation}\label{M.45}
			\begin{array}{ll}
				\di \iint \f{|
					\wt{f}_{yy}|^2+|
					\wt{f}_{yt}|^2}{\mb{M}_{_\#}}d\xi dy+\int_0^T\iint \f{(1+|\xi|)\left( | \wt{\mb{G}}_{yy}|^2+|\wt{\mb{G}}_{yt}
					|^2\right)}{\mb{M}_{_\#}}d\xi dydt\\[3mm]
				\di \leqslant C\mathcal{N}(0)^2+C\delta_S^2\int_0^T|\dot{\mb X}(t)|^2 dt+C(\eta_0+\d_1+\ve_1)\int_0^T \Lambda^Sdt+C\d_1\\[3mm]
				\di\quad
				+C(\eta_0+\d_1+\ve_1)\sum_{|\b|=1}\int_0^T\|\partial^{\b}(\p,\psi,\z)\|_{L^2}^2dt\\[3mm]
				\di\quad
				+C(\eta_0+\d_1+\ve_1)\int_0^T\left( \|(\p_{yy},\psi_{yy},\z_{yy})\|_{L^2}^2+\|(\p_{yt},\psi_{yt},\z_{yt})\|_{L^2}^2\right) dt\\[3mm]
				\di\quad
				+C(\eta_0+\d_1+\ve_1)\int_0^T\iint
				\f{(1+|\xi|)\left( |\wt{\mb{G}}_1|^2+|\wt{\mb{G}}_y|^2+|\wt{\mb{G}}_t|^2\right) }{\mb{M}_{_\#}}d\x
				dydt,
			\end{array}
		\end{equation}
		where $\eta_0$ defined in Lemma \ref{Lemma 4.3} is chosen to be suitably small , which is crucial to close the a priori assumptions. Actually we can choose $\eta_0=O(1)(\d_1+\ve_1)$.

		Noting that
		$$
		\begin{array}{ll}
			\di \left| \int\psi_{1y}\phi_{yy} dy\right| &\di \leqslant \tilde\nu \|\psi_{1y}\|^2
			+C\|\phi_{yy}\|_{L^2}^2\\
			&\di \leqslant \tilde\nu \|\psi_{1y}\|_{L^2}^2 +C\iint
			\f{|\wt{f}_{yy}|^2}{\mb{M}_{_\#}}d\xi
			dy+C\ve_1^2\|\phi_{y}\|_{L^2}^2+C\d_1,
		\end{array}
		$$
		and combining the estimates \eqref{Step2final}, \eqref{Gy},
		\eqref{Gt}, and \eqref{M.45},
		we complete the proof of Proposition \ref{higherorderestimate}. 
		
		Finally, to prove Proposition \ref{main proposition}, it suffices to estimate 
		$$\di C\delta_C \int_0^T(1+t)^{-1}\int  e^{-\frac{c_0|y+\sigma t|^2}{1+t}}|(\phi,\psi,\z)|^2 dydt,$$
		which is given in the following proposition.
		\begin{proposition}\label{Proposition 6.2}
			Under the assumptions of Proposition \ref{main proposition}, there exists a positive
			constant C such that
			\begin{equation*}
				\begin{array}{ll}
					\di \quad \int_0^T(1+t)^{-1}\int  e^{-\frac{c_0|y+\sigma t|^2}{1+y}}|(\phi,\psi,\z)|^2 dy dt \\[3mm]
					\di \leqslant  C\sup_{t\in[0,T]}\|(\phi,\psi,\z)\|^2_{L^2}  +C\delta_S\int_0^T|\dot{\mb X}(t)|^2 dt+C\int_0^T\left( \Lambda^R+\Lambda^S\right) dt  + C\int_0^T\|(\phi_y,\psi_y,\z_y)\|_{L^2}^2 dt\\[3mm]
					\di\quad
					+C(\d_1+\ve_1)\int_0^T\iint \frac{ (1+|\xi|)|\wt{\mb{G}}_1|^2}{\mb M_{_\#}}d\xi dydt+C\int_0^T\iint \frac{ (1+|\xi|)(|\wt{\mb{G}}_t|^2+|\wt{\mb{G}}_y|^2)}{\mb M_{_\#}}d\xi dydt+C\delta_1^{\frac13}.
				\end{array}
			\end{equation*}
		\end{proposition}
		
		The proof of Proposition \ref{Proposition 6.2} uses the same argument as in \cite{ Huang-Li-Matsumura, Kang-Vasseur-Wang-2024}, with additional treatment for the microscopic part. Since the proof is lengthy and mainly follows from  \cite{ Huang-Li-Matsumura, Kang-Vasseur-Wang-2024}, we omit it for simplicity. Finally,  the proof of Proposition \ref{main proposition} can be completed directly from Proposition \ref{lowerorderestimate}, Proposition \ref{higherorderestimate}, and Proposition \ref{Proposition 6.2}.

	\
	
	\begin{appendix}
		\setcounter{equation}{0}
		\section{Appendix}
		\renewcommand{\thetheorem}{A.\arabic{theorem}}
		\renewcommand{\thelemma}{A.\arabic{lemma}}
		\renewcommand{\theproposition}{A.\arabic{proposition}}
		\renewcommand{\theequation}{A.\arabic{equation}}
		\subsection*{A.1\quad Proof of Lemma \ref{Lemma-shock}}
		In this subsection we prove Lemma \ref{Lemma-shock}. For Boltzmann shock profile $F^S(y,\xi)$ with $y=x-\s t$, it can be parameterized by $\eta(y)$, which satisfies the following Burgers-like equation:
		\begin{eqnarray}\label{eta}
			\begin{cases}
				\di\frac{d\eta}{dy}=\Phi_1(\eta)(\eta-\eta_-)(\eta-\eta_+), \\[3mm]
				\di\eta(-\infty)=\eta_{-},\ \  \eta(+\infty)=\eta_{+},
			\end{cases}
		\end{eqnarray}
		where $\eta_->0>\eta_+$ are constants with $|\eta_--\eta_+|=O(\d_S)$, and $\Phi_1(z)$ is a smooth function on $[\eta_+,\eta_-]$ satisfying $\Phi_1(z)\geqslant c>0$ for some constant $c$. The following lemma gives some properties of $\eta$, which can be obtained by straightforward calculations.

		\begin{lemma}\label{etaproperty}
			The solution $\eta(y)$ to the Burgers-like equation \eqref{eta} satisfies
			\begin{align*}
				\begin{aligned}
					& \frac{d\eta}{dy} <0, ~~~~ \forall y\in\mathbb{R},\\
					& | \eta-\eta_-|\leqslant C\delta_S\ e^{-C\delta_S |y|}, \quad y<0,\\[1mm]
					& | \eta-\eta_+|\leqslant C\delta_S\ e^{-C\delta_S |y|}, \quad y>0,\\[1mm]
					& \left| \frac{d\eta}{dy}\right| \leqslant C\d_S^2\ e^{-C\delta_S |y|},\\[1mm]
					& \di \left| \frac{d^k\eta}{dy^k}\right| \leqslant
					C\d_{S}^{k-1}\left| \frac{d\eta}{dy}\right| , \quad k\geqslant 2.
				\end{aligned}
			\end{align*}		
			More precisely, it holds that
			\begin{equation*}
				\lim\limits_{y\to -\i}\frac{\eta_{y}}{\eta-\eta_-}=O(\d_S),~~~~~~\lim\limits_{y\to -\i}\frac{\eta_{yy}}{\eta_y}=O(\d_S),
			\end{equation*}
			and
			\begin{equation}\label{etasharp}
				\frac{\eta_{yyy}\eta_y-(\eta_{yy})^2}{(\eta_y)^3}=O(1).
			\end{equation}
		\end{lemma}
		
		\
		Under the above preparation, we are now at the stage to prove Lemma \ref{Lemma-shock}. In fact, the macroscopic part has the property
		\begin{equation*}
			\di v^{S}_y\sim u^{S}_{1y}\sim
			\t^{S}_y\sim \eta_y.
		\end{equation*}
		The microscopic part $\mb G^S(y,\xi)$ can be further expressed as
		\begin{equation}\label{GSprecise}
			\mb G^S(y,\xi)=\sqrt{\mb M_0} \Gamma(\omega(\eta),\xi),
		\end{equation}
		where $\mb{M}_0$ is a global Maxwellian  and $\Gamma(\omega(\eta),\xi)$ is a smooth function of $(\omega(\eta),\xi)$. For $\o(\eta)$, it holds that
		\begin{equation}\label{omega}
			\o(\eta)=\Psi(\eta,\o(\eta))(\eta-\eta_-)(\eta-\eta_+),
		\end{equation}
		where $\Psi(z_1,z_2)$ is a smooth function on $[\eta_+,\eta_-]\times\mathbb R$ satisfying $0<c\leqslant\Psi(z_1,z_2)\leqslant C<+\i$ for constants $c$ and $C$. By \eqref{GSprecise} and  \cite{Liu-Yu-2013}, we have
		\begin{equation}\label{Gestimate}
			\di\left( \int
			\f{(1+|\xi|)|\mb{G}^{S}|^2}{\mb{M}_0}d\x\right) ^{\f12}=\left(\int (1+|\xi|)\Gamma^2(\omega(\eta),\xi)d\x\right)^{\f12}=| \o |\Phi_2(\o),
		\end{equation}
		where $\Phi_2(z)$ is a smooth bounded function satisfying $\Phi_2(z)\geqslant c>0$. Therefore, using Lemma \ref{etaproperty} and \eqref{Gestimate}, together with the fact that
		\begin{equation*}
			\di v^{S}_y\sim u^{S}_{1y}\sim
			\t^{S}_y\sim \eta_y\sim\o,
		\end{equation*}
		we have \eqref{shock-macro}-\eqref{equivalence}.
		
		\eqref{shock-vu} is a direct consequence of \eqref{sm1} and $\eqref{shock-profile}_1$. To prove \eqref{theta-s}, we first note that since $v^S_y>0$ for all $y$, the function $v^S(y)$ admits an inverse function: $(v^S)^{-1}:[v^*,v_+]\to \mathbb{R}$. Thus we can view $\t^S, u^S_1$ and $\Pi^S_1$ as a smooth function of $v^S$. Therefore, $\di\frac{d\t^S}{dv^S}=\frac{\t^S_y}{v^S_y}$. Precisely, by \eqref{VS2}, we have
		\begin{equation}\label{first-der0}
			\begin{array}{ll}
				\di \frac{d\t^S}{dv^S}&\di =\frac{4\s^2\mu(\t^S)}{3\kappa(\t^S)}\cdot\frac{(\theta^S-\theta^*)+p^*(v^S -v^*)-\frac12\sigma^2(v^S-v^*)^2-\di \frac{1}{\s}\int\xi_1\left( \frac{1}{2}|\xi|^2-u^S_1 \xi_1\right) \Pi^S_1d\xi}{(p^S-p^*)+\sigma^2(v^S-v^*)+\di \int\xi_1^2\Pi_1^Sd\x}\\
				&\di :=\mathcal F(v^S).
			\end{array}
		\end{equation}
		Direct calculations yield
		\begin{equation*}
			l^*:=\mathcal F(v^*)= \lim_{y\to-\infty} \frac{\theta^S_y}{v^S_y} = \lim_{v^S\to v^*} \frac{d\theta^S}{dv^S} =  \lim_{v^S\to v^*} \frac{\theta^S-\theta^*}{v^S-v^*}<0.
		\end{equation*}
		Then
		\begin{equation*}
			\lim_{v^S\rightarrow v^*}\frac{p^S-p^*}{v^S-v^*}=\lim_{v^S\rightarrow v^*}\Big(\frac{2}{3v^S}\frac{\theta^S-\theta^*}{v^S-v^*}-\frac{p^*}{v^S}\Big)
			=\frac{\f23l^*-p^*}{v^*}.
		\end{equation*}
		By \eqref{Pi1S}, \eqref{GSprecise}, and \eqref{Gestimate}, one can get
		\begin{equation}\label{Pi1Sprecise} 
			\int\xi_1^2\Pi_1^Sd\x=g^1(\o)\o_y+g^2(\o)\o^2,
		\end{equation}
		where $g^1(\o)$ and $g^2(\o)$ are smooth bounded functions. By \eqref{omega}, one has 
		$$\lim\limits_{y\to -\i}\frac{\o_y}{\eta-\eta_-}=O(\d_S^2),$$
		which together with \eqref{Pi1Sprecise}  yields that
		\begin{equation*}
			S_1:=\lim_{v^S\rightarrow v^*}\frac{\di \int\xi_1^2\Pi_1^Sd\x}{v^S-v^*}=\lim_{v^S\rightarrow v^*}\frac{g^1(\o)\o_y+g^2(\o)\o^2}{v^S-v^*}\sim \lim\limits_{y\to -\i}\frac{g^1(\o)\o_y+g^2(\o)\o^2}{\eta-\eta_-} =O(\d_S^2),
		\end{equation*}
		and similarly,
		\begin{equation*}
			S_2:=\di\lim_{v^S\rightarrow v^*}\f1\s\frac{\di \int\xi_1\left( \frac{1}{2}|\xi|^2-u^S_1 \xi_1\right) \Pi^S_1d\xi}{v^S-v^*}=O(\d_S^2).
		\end{equation*}
		Now dividing both the numerator and denominator by $v^S-v^*$ on the right-hand side of \eqref{first-der0}  and then passing the limit $v^S\to v^*$, we have
		$$
		l^* =\frac{4\s^2\mu(\t^*)}{3\kappa(\t^*)}\frac{ l^* +p^*-S_2}{\frac{\f23l^*-p^*}{v^*} +\sigma^2+S_1}.
		$$
		That is, we have the following equation for $l^*$,
		\begin{equation}\label{l*}
			(l^*)^2+\left[ -\f32(p^*-\s^2v^*)-2\frac{\mu(\t^*)}{\kappa(\t^*)}\s^2 v^*+\frac{3v^*}{2}S_1\right] l^*-\frac{4\mu(\t^*)\s^2\theta^*}{3\kappa(\t^*)}+2\s^2\frac{\mu(\t^*)}{\kappa(\t^*)}v^*S_2=0.
		\end{equation}
		Note that $-p^*$ satisfies the identity
		\begin{equation}\label{-p*}
			(-p^*)^2+\left(-\f32(p^*-(\s^*)^2v^*)-2\frac{\mu(\t^*)}{\kappa(\t^*)}(\s^*)^2 v^*\right)(-p^*)-\frac{4\mu(\s^*)^2\theta^*}{3\kappa(\t^*)}=0.
		\end{equation}
		Subtracting \eqref{-p*} from \eqref{l*} 
		and using \eqref{sm1}, we have
		\begin{equation}\label{l*andl*}
			|l^*-(-p^*)|=|\mathcal F(v^*)+p^*|= O(\delta_S),
		\end{equation}
		which together with $|\mathcal F(v^S)-\mathcal F(v^*)|=O(1)|v^S-v^*|=O(\delta_S)$ yields that
		\[
		\left| \frac{\theta^S_y}{v^S_y} + p^*\right| \leqslant C\delta_S,\qquad \forall y\in\mathbb R.
		\]
		Thus \eqref{theta-s} is proved.
		
		Now it remains to prove \eqref{2ndorder}.   We need the following three lemmas. The first one comes from \cite{Kawashima-Matsumura-Nishida}. 
		
		\begin{lemma}\label{mu and kappa}
			For the hard sphere model,  the viscosity coefficient $\mu(\t)$ and the heat conductivity coefficient $\k(\t)$ have the explicit formula
			$$\mu(\t)=A_1\t^{\f12},~~~~\k(\t)=A_2\t^{\f12}, $$
			for positive constants $A_1$ and $A_2$.
		\end{lemma}
		
		The next two lemmas are used to estimate the microscopic part $\Pi^S_1$. 
		\begin{lemma}\label{eta4property}
			For the function $\o$ defined in \eqref{omega}, it holds that
			$$\left| \o\right|\leqslant C\delta_S^2,~~~ \left| \o_y\right| \leqslant C\delta_S\left| \eta_y\right| , ~~~\left| \o_{yy}\right|\leqslant C\delta_S\left| \eta_{yy}\right|+C(\eta_y)^2,$$
			and
			\begin{equation}\label{eta4property1}
				\frac{\di\o_{yyy}\eta_y-\o_{yy}\eta_{yy}}{(\eta_y)^3}=O(\delta_S).
			\end{equation}
		\end{lemma}
		{\it Proof.}  Here we only prove \eqref{eta4property1} and the other properties can be derived from \eqref{omega} directly, and we omit it. Differentiating \eqref{omega} with respect to $y$ yields that
		$$\o_y=\Psi(\eta,\o)(2\eta-\eta_--\eta_+)\eta_y+(\eta-\eta_-)(\eta-\eta_+)\left( \Psi_\eta(\eta,\o)\eta_y+\Psi_\o(\eta,\o)\o_y\right), $$
		\begin{align*}
			\o_{yy}=&\left( \Psi(\eta,\o)(\eta-\eta_-)(\eta-\eta_+)\right)_{yy}\\
			=&\Psi(\eta,\o)\left( (2\eta-\eta_--\eta_+)\eta_{y}\right)_{y}+2\left(\Psi(\eta,\o) \right)_y (2\eta-\eta_--\eta_+)\eta_{y}\\
			&+\left(\Psi(\eta,\o) \right)_{yy}(\eta-\eta_-)(\eta-\eta_+),
		\end{align*}
		and
		\begin{align*}
			\o_{yyy}=&\Psi(\eta,\o) (2\eta-\eta_--\eta_+)\eta_{yyy}+4\Psi(\eta,\o) \eta_y\eta_{yy}+2\Psi(\eta,\o) \eta_y\eta_{yy}\\
			&+3\left(\Psi(\eta,\o) \right)_y (2\eta-\eta_--\eta_+)\eta_{yy}+6\left(\Psi(\eta,\o) \right)_y(\eta_y)^2\\
			&+3\left(\Psi(\eta,\o) \right)_{yy}(2\eta-\eta_--\eta_+)\eta_{y}+\left(\Psi(\eta,\o) \right)_{yyy}(\eta-\eta_-)(\eta-\eta_+).
		\end{align*}
		Substituting the above three relations into the left-hand side of \eqref{eta4property1} implies that
		\begin{align*}
			\o_{yyy}\eta_y-\o_{yy}\eta_{yy}=&\Psi(\eta,\o) (2\eta-\eta_--\eta_+)\big(\eta_{yyy}\eta_y-(\eta_{yy})^2\big)\\
			&+(\eta-\eta_-)(\eta-\eta_+)\Psi_\eta(\eta,\o)\big(\eta_{yyy}\eta_y-(\eta_{yy})^2\big)\\
			&+(\eta-\eta_-)(\eta-\eta_+)\Psi_\o(\eta,\o)(\o_{yyy}\eta_y-\o_{yy}\eta_{yy})\\
			&+O(\delta_S)(\eta_y)^3,
		\end{align*}
		which together with \eqref{etasharp} completes the proof of \eqref{eta4property1}.
		
		\
		
		\begin{lemma}\label{Piestimate}
			For the collision invariants $\varphi_i ~(i=1, 4)$, it holds that
			\begin{equation}\label{Piestimate1}
				\displaystyle \frac{\di \int\xi_1\varphi_i(\xi) \Pi^S_{1yy}d\xi v^S_y-\int\xi_1\varphi_i(\xi) \Pi^S_{1y}d\xi v^S_{yy}}{\left( v^S_y\right)^3 }=O(\delta_S).
			\end{equation}
		\end{lemma}
		{\it Proof.} We first consider the case when $i=1$. Substituting \eqref{Pi1Sprecise} into the left-hand side of \eqref{Piestimate1} and using  \eqref{eta4property1}, we have
		\begin{equation*}
			\begin{array}{ll}
				\di \frac{\di\int\xi_1^2 \Pi^S_{1yy}d\xi v^S_y-\int\xi_1^2 \Pi^S_{1y}d\xi v^S_{yy}}{\left( v^S_y\right)^3 }=O(1)\frac{\di\o_{yyy}\eta_y-\o_{yy}\eta_{yy}}{(\eta_y)^3}+O(\d_S)=O(\d_S).
			\end{array}
		\end{equation*} 
		The case when $i=4$ is similar, and we complete the proof of Lemma \ref{Piestimate}.
		
		\
		
		Now it is ready to prove \eqref{2ndorder}. By Taylor expansion (cf. \cite{Kang-Vasseur-Wang-2024}), it holds that
		\begin{equation}\label{Taylor}
			\begin{array}{ll}
				\di \frac{p^S-p_+}{v^S-v_+}-\frac{p^S-p^*}{v^S-v^*}=\frac{1}{3v^*} \frac{d^2\theta^S}{d(v^S)^2}\bigg|_{v^S=v^*} (v_+-v^*) + \frac{5 p^*}{3(v^*)^2}(v_+-v^*)+O(\delta_S^2).
			\end{array}
		\end{equation}
		Then we compute $\di\frac{d^2\t^S}{d(v^S)^2}\bigg|_{v^S=v^*}$. Due to Lemma \ref{mu and kappa}, it holds that $\di\frac{\mu(\t^S)}{\kappa(\t^S)}=\frac{\mu(\t^*)}{\kappa(\t^*)}$ is a  positive constant. Differentiating \eqref{first-der0} with respect to $v^S$, we have
		\begin{equation*}
			\begin{array}{ll}
				\di 
				\frac{d^2\theta^S}{d(v^S)^2}=\frac{4\s^2\mu(\t^*)}{3\kappa(\t^*)}\cdot\frac{\frac{d\theta^S}{dv^S}+p^*-\sigma^2(v^S-v^*)-\frac{1}{\s}(v^S_y)^{-1}\left( \di \int\xi_1\left( \frac{1}{2}|\xi|^2-u^S_1 \xi_1\right) \Pi^S_1d\xi\right)_y }{(p^S-p^*)+\sigma^2(v^S-v^*)+\di \int\xi_1^2\Pi_1^Sd\x}\\[8mm]
				\di \qquad\quad-\frac{4\s^2\mu(\t^*)}{3\kappa(\t^*)}\cdot\frac{(\theta^S-\theta^*)+p^*(v^S -v^*)-\frac12\sigma^2(v^S-v^*)^2-\frac{1}{\s}\di \int\xi_1\left( \frac{1}{2}|\xi|^2-u^S_1 \xi_1\right) \Pi^S_1d\xi}{\left[ (p^S-p^*)+\sigma^2(v^S-v^*)+\di \int\xi_1^2\Pi_1^Sd\x\right] ^2}\\[10mm]
				\di\qquad\qquad\quad
				\cdot\left(\frac{2\frac{d\theta^S}{dv^S}}{3v^S}-\frac{2\theta^S}{3(v^S)^2}+\s^2+(v^S_y)^{-1}\di \int\xi_1^2\Pi_{1y}^Sd\x \right) \\[10mm]
				\di
				\qquad \quad\ =\frac{4\s^2\mu(\t^*)}{3\kappa(\t^*)}\cdot\frac{\frac{d\theta^S}{dv^S}+p^*-\sigma^2(v^S-v^*)-\frac{1}{\s}(v^S_y)^{-1} \left( \di \int\xi_1\left( \frac{1}{2}|\xi|^2-u^S_1 \xi_1\right) \Pi^S_1d\xi\right)_y}{(p^S-p^*)+\sigma^2(v^S-v^*)+\di \int\xi_1^2\Pi_1^Sd\x}\\[6mm]
				\di\qquad\qquad\quad -\frac{\frac{d\theta^S}{dv^S}}{(p^S-p^*)+\sigma^2(v^S-v^*)+\di \int\xi_1^2\Pi_1^Sd\x}\cdot\left(\frac{2\frac{d\theta^S}{dv^S}}{3v^S}-\frac{2\theta^S}{3(v^S)^2}+\s^2+(v^S_y)^{-1}\di \int\xi_1^2\Pi_{1y}^Sd\x \right)\\[10mm]
				\di\qquad \quad\  =-\frac{2}{3v^S  \Big[(p^S-p^*)+\sigma^2(v^S-v^*)+\di \int\xi_1^2\Pi_1^Sd\x\Big]}\bigg\{\Big(\frac{d\theta^S}{dv^S}\Big)^2+\frac{d\theta^S}{dv^S}\Big[-\f32(p^S-\s^2v^S)\\[10mm]
				\di\qquad\qquad\quad -2 \frac{\mu(\t^*) }{\kappa(\t^*)}\s^2v^S+\frac{3v^S}{2v^S_y}\int\xi_1^2\Pi_{1y}^Sd\x\Big]-2 \frac{\mu(\t^*) }{\kappa(\t^*)}\s^2v^S\Big[ p^*-\sigma^2(v^S-v^*)\\[8mm]    
				\di \qquad\qquad\quad -\frac{1}{\s}(v^S_y)^{-1}  \Big( \int\xi_1\left( \frac{1}{2}|\xi|^2-u^S_1 \xi_1\right) \Pi^S_1d\xi\Big)_y\Big]  
				\bigg\}.
			\end{array}
		\end{equation*}
		By using \eqref{l*} and the fact that
		$$
		\begin{array}{ll}
			\di -2\frac{\mu(\t^*)}{\kappa(\t^*)}\s^2v^S\left(p^*-\s^2(v^S-v^*) \right)  +\f43\frac{\mu(\t^*)}{\kappa(\t^*)}\s^2\theta^*&\di =-2\frac{\mu(\t^*)}{\kappa(\t^*)}\s^2\Big[v^S\big(p^*-\s^2(v^S-v^*\big)-\frac 23 \theta^*\Big]\\[4mm]
			&\di
			=-2\frac{\mu(\t^*)}{\kappa(\t^*)}\s^2\left( p^*-\s^2 v^S\right) \left( v^S-v^*\right) ,
		\end{array}$$
		we have
		\begin{equation}\label{(6.21)}
			\begin{array}{ll}
				\di 
				-\f32v^S\left[ (p^S-p^*)+\sigma^2(v^S-v^*)+\int\xi_1^2\Pi_1^Sd\x\right]\frac{d^2\theta^S}{d(v^S)^2}=\left(\frac{d\t^S}{dv^S} \right)^2-(l^*)^2 \\[4mm]
				\di\qquad
				-\left[ \f32(p^*-\s^2v^*)+2\frac{\mu(\t^*)}{\kappa(\t^*)}\s^2 v^*\right]\left( \frac{d\t^S}{dv^S} -l^* \right)-2\frac{\mu(\t^*)}{\kappa(\t^*)}\s^2\left( p^*-\s^2 v^S\right) \left( v^S-v^*\right)\\[4mm] 
				\di\qquad
				+\left[-\f32\left((p^S-p^*)-\s^2(v^S-v^*)\right)-2\frac{\mu(\t^*)}{\kappa(\t^*)}\s^2(v^S-v^*)   \right]\frac{d\t^S}{dv^S}   \\[4mm]
				\di \qquad 
				+\f32\Big(  v^S(v^S_y)^{-1}\di \int\xi_1^2\Pi_{1y}^Sd\x\frac{d\t^S}{dv^S}-v^*S_1l^*\Big)  \\[4mm]
				\di \qquad
				+2\frac{\mu(\t^*)}{\kappa(\t^*)}\s^2\left[ \frac{1}{\s}v^S\left( \int\xi_1\left( \frac{1}{2}|\xi|^2-u^S_1 \xi_1\right) \Pi^S_1d\xi\right)_y(v^S_y)^{-1}-v^*S_2\right].
			\end{array}
		\end{equation}
		
		By Lemma \ref{Piestimate}, one has
		\begin{equation*}
			\begin{array}{ll}
				\di\quad 
				\lim_{v^S\to v^*}\frac{1}{v^S-v^*}\bigg(  v^S (v^S_y)^{-1}\di \int\xi_1^2\Pi_{1y}^Sd\x\frac{d\t^S}{dv^S}-v^*S_1l^*\bigg) \\[4mm]
				\di
				=\lim_{v^S\to v^*}\left(v^S(v^S_y)^{-1}\di \int\xi_1^2\Pi_{1y}^Sd\x\frac{d\t^S}{dv^S} \right)_y(v^S_y)^{-1}\\[4mm]
				\di
				=\lim_{y\to -\i}\Bigg(\di (v^S_y)^{-1} \di\int\xi_1^2\Pi_{1y}^Sd\x\frac{d\t^S}{dv^S}+v^S\frac{d\t^S}{dv^S}(v^S_y)^{-3}\bigg(\int\xi_1^2\Pi^S_{1yy}d\xi v^S_y-\int\xi_1^2 \Pi^S_{1y}d\xi v^S_{yy}\bigg)\\
				\di \qquad\qquad \quad +v^S(v^S_y)^{-1}\di \int\xi_1^2\Pi_{1y}^Sd\x\frac{d^2\theta^S}{d(v^S)^2}\Bigg) \\[4mm]
				\di =O(\d_S).
			\end{array}
		\end{equation*}
		Similarly, one has
		$$\lim_{v^S\to v^*}\frac{1}{v^S-v^*}\left[ \frac{1}{\s}v^S\left( \int\xi_1\left( \frac{1}{2}|\xi|^2-u^S_1 \xi_1\right) \Pi^S_1d\xi\right)_y(v^S_y)^{-1}-v^*S_2\right]=O(\d_S).$$
		Dividing both sides of \eqref{(6.21)} by $\di v^S-v^*$ and then taking the limit $v^S\to v^*$, we have
		\begin{equation*}
			\begin{array}{ll}
				\di
				-\f32v^*\left[\frac{\f23l^*-p^*}{v^*}+\sigma^2+O(\d_S)\right]  \frac{d^2\t^S}{d(v^S)^2}\bigg|_{v^S=v^*}\\[7mm]
				\di =\left[\lim_{v^S\to v^*}\left( \frac{d\t^S}{dv^S}+l^*\right) - \f32(p^*-\s^2v^*)-2\frac{\mu(\t^*)}{\kappa(\t^*)}\s^2 v^*\right]\lim_{v^S\to v^*}\frac{\frac{d\t^S}{dv^S}-l^*}{v^S-v^*}\\[7mm]
				\di\quad
				-2\frac{\mu(\t^*)}{\kappa(\t^*)}\s^2\left( p^*-\s^2 v^*\right)-\f32\left( \frac{\f23l^*-p^*}{v^*}-\s^2\right)l^* -2\frac{\mu(\t^*)}{\kappa(\t^*)}\s^2l^*+O(\d_S)\\[7mm]
				\di =\left[ 2l^*- \f32(p^*-\s^2v^*)-2\frac{\mu(\t^*)}{\kappa(\t^*)}\s^2 v^*\right]\frac{d^2\t^S}{d(v^S)^2}\bigg|_{v^S=v^*}\\[7mm]
				\di\quad
				-2\frac{\mu(\t^*)}{\kappa(\t^*)}\s^2\left( p^*-\s^2 v^*\right)-\f32\left( \frac{\f23l^*-p^*}{v^*}-\s^2\right)l^* -2\frac{\mu(\t^*)}{\kappa(\t^*)}\s^2l^*+O(\d_S),
			\end{array}
		\end{equation*}
		which together with \eqref{sm1} and \eqref{l*andl*} implies
		\begin{equation}\label{A.19}
			\di\frac{d^2\t^S}{d(v^S)^2}\bigg|_{v^S=v^*}=\frac{10\mu(\t^*)-9\k(\t^*)}{10\mu(\t^*)+3\k(\t^*)}\frac{ 5p^*}{3v^*} +O(\d_S).
		\end{equation}
		Substituting \eqref{A.19} into \eqref{Taylor} completes the proof of \eqref{2ndorder}.
		
		\subsection*{A.2 \quad Proof of Proposition \ref{localexistence}}
		
		In this subsection, we prove Proposition \ref{localexistence} for the local-in time existence of the solution to the
		Boltzmann equation \eqref{tran-Lag-B} in Lagrangian coordinates. Motivated by the iteration scheme in Eulerian coordinates (cf. \cite{Duan-Huang-Wang-Yang,Ukai-Yang-Zhao}), we define the iteration sequence $f^n(t, x, \xi)~( n=0, 1,2, \cdots$) in Lagrangian coordinates by
		\begin{equation}\label{iteration-f}
			\begin{cases}\di f_t^{n+1}+\left(\frac{\xi_1}{v^n}-\frac{u_1^n}{v^n}-\s \right) f^{n+1}_y+f^{n+1}\int_{\mathbb R^3}\int_{\mathbb S^2_+}|(\xi-\xi_*)\cdot \Omega|f^n(t,y,\xi_*)d\Omega d \xi_* =Q_+(f^n,f^n),
				\\[3mm]
				f^{n+1}(0,y,\xi)=f_0(y,\xi)\geqslant0,
				\\[2mm]
				f^0(t,y,\xi)=\hat{\mathbf{M}}(t,y,\xi),\end{cases}
		\end{equation}
		where $\hat{\mb M}(t,y,\xi):=\mb M_{[\hat v,\hat u,\hat \t](t,y)}(\xi)$ is defined in \eqref{hat-v} and \eqref{hat-M}.
		
		Set 
		\begin{equation*}
			g^n(t,y,\xi):=f^n(t,y,\xi)-\hat{\mb M}(t,y,\xi),
		\end{equation*}
		then $g^n(t,y,\xi)$ satisfies the equation
		\begin{equation}\label{iteration-g}
			\begin{cases}\di g_t^{n+1}+\left(\frac{\xi_1}{v^n}-\frac{u_1^n}{v^n}-\s \right) g^{n+1}_y+g^{n+1}\int_{\mathbb R^3}\int_{\mathbb S^2_+}|(\xi-\xi_*)\cdot \Omega|\left( g^n(t,y,\xi_*)+\hat{\mb M}(t,y,\xi_*)\right)d\Omega d \xi_* \\[4mm]
				~~~~~~~~~\di=\sqrt{\hat{\mb M}}\mb K_{\hat{\mb M}}\left(\frac{g^n}{\sqrt{\hat{\mb M}}} \right) +Q_+(g^n,g^n)-\left( \hat{\mathbf{M}}_t+\left(\frac{\xi_1}{v^n}-\frac{u_1^n}{v^n}-\s \right)\hat{\mathbf{M}}_y\right),\\[4mm]
				g^{n+1}(0,y,\xi)=g_0(y,\xi),\\[1mm]
				g^0(t,y,\xi)=0,\end{cases}
		\end{equation}
		where $\mb K_{\hat{\mb M}}$ is the compact $L^2$-operator defined in \eqref{LM} and \eqref{4.2} and  $g_0(y,\xi):=f_0(y,\xi)-\hat{\mb M}(0,y,\xi)$. For $T>0$, the solution function space is defined by
		$$
		\left.\mathcal{S}_N([0,T]):=\left\{g(t,y,\xi)\begin{array}{c}\left|\begin{array}{c}\di\frac{\partial^\alpha_y g(t,y,\xi)}{\sqrt{\mathbf{M}_{_\#}(\xi)}}\in C\left([0,T],L_{y,\xi}^2\left(\mathbb{R}\times\mathbb{R}^3\right)\right) ,~ 0\leqslant|\a|\leqslant N,\\
				\di\|g\|_{\mathcal{S}_N}\leqslant\Xi,\end{array}\right.\end{array}\right.\right\}
		$$
		for some small positive constant $\Xi$ and 
		$$
		\begin{aligned}\|g\|_{\mathcal{S}_N}^{2}:=\sum_{0\leqslant|\alpha|\leqslant N}\left(\sup_{t\in[0,T]}\iint\frac{|\partial^\a_y g(t,y,\xi)|^2}{\mathbf{M}_{_\#}(\xi)}d\xi dy+\int_0^T\iint\frac{\nu_{\hat{\mb M}}(\xi)|\partial^\a_y g(t,y,\xi)|^2}{\mathbf{M}_{_\#}(\xi)}d\xi dydt\right).\end{aligned}
		$$
		
		We want to prove that under the the initial condition $f_0(y,\xi)\geqslant 0$ and
		$$\sum_{0\leqslant|\alpha|\leqslant2}\iint\frac{|\partial^\a_y g_0(y,\xi)|^2}{\mathbf{M}_{_\#}(\xi)}d\xi dy< \f1{4}\Xi^2<\ve_2^2,$$
		where $\ve_2>0$ is a small positive constant to be determined later, the induction on $n$ that 
		\begin{equation}\label{iteration-proposition}
			\begin{array}{ll}
				\mbox{the}~\mbox{assumption}~f^n(t,y,\xi)\geqslant0,~ \rho^n(t,y)>0,~\mbox{and}~ \|g^n\|_{\mathcal{S}_2}^2\leqslant\Xi^2\\[3mm]
				\mbox{implies}~\mbox{that}~f^{n+1}(t,y,\xi)\geqslant0,~ \rho^{n+1}(t,y)>0,~\mbox{and}~ \|g^{n+1}\|_{\mathcal{S}_2}^2\leqslant\Xi^2,
			\end{array}
		\end{equation}
		holds true provided that $\ve_2,\d_2$ and $T$ are chosen suitably small. Therefore, by the condition when $n=0$ that
		$$
		f^0(t,y,\xi)\geqslant 0,\quad \rho^0(t,y)=\frac{1}{\hat v(t,y)}>0, \ \mbox{and}\ \|g^0\|_{\mathcal{S}_2}^2\equiv0\leqslant\Xi^2
		$$ 
		and the induction \eqref{iteration-proposition}, we have $\forall n=0, 1, 2, \cdots,$
		$$
		f^n(t,y,\xi)\geqslant0,\quad \rho^n(t,y)>0,~\mbox{and}~ \|g^n\|_{\mathcal{S}_2}^2\leqslant\Xi^2.
		$$
		
		Now we prove the induction \eqref{iteration-proposition}. Assume that $f^n(t,y,\xi)\geqslant0$, $\rho^n(t,y)>0$, and $\|g^n\|_{\mathcal{S}_2}^2\leqslant\Xi^2$ for some $n=0, 1, 2, \cdots$. 
		
		First we prove $f^{n+1}(t,y,\xi)\geqslant0$. The backward characteristic line $Y(\tau;t,y,\xi)$ passing through the point $(t,y,\xi)\in \mathbb{R}_+\times\mathbb{R}\times\mathbb{R}^3$  to the linear equation \eqref{iteration-f} can be defined by
		\begin{equation*}\label{characteristic}
			\begin{cases}\di \frac{dY(\tau;t,y,\xi)}{d\tau}=\left(\frac{\xi_1}{v^n}-\frac{u_1^n}{v^n}-\s \right)(\tau,Y(\tau);t,y,\xi),
				\\Y(\tau;t,y,\xi)|_{\tau=t}=y.\end{cases}
		\end{equation*}
		Set 
		$$k^n(\tau,z,\xi):=\int_{\mathbb R^3}\int_{\mathbb S^2_+}|(\xi-\xi_*)\cdot \Omega|f^n(\tau,z,\xi_*)d\Omega d \xi_* .$$
		Then the solution of the linear equation \eqref{iteration-f} can be expressed as
		$$
		\begin{aligned}
			f^{n+1}(t,y,\xi) =&\di e^{\di-\int_{0}^{t} k^n(\tau,Y(\tau;t,y,\xi),\xi) d\tau}\cdot f_{0}\left( Y(0;t,y,\xi),\xi\right)  \\
			&\di  +\int_0^te^{\di-\int_s^t k^n(\tau,Y(\tau;t,y,\xi),\xi)d\tau}\cdot Q_+(f^n,f^n)\left( s,Y(s;t,y,\xi),\xi\right) ds,
		\end{aligned}
		$$
		which implies that  $f^{n+1}(t,y,\xi)\geqslant0$ since  both $f_0(y,\xi)\geqslant0$ and $f^{n}(t,y,\xi)\geqslant0$.

		Next we prove $\|g^{n+1}\|_{\mathcal{S}_2}^2\leqslant\Xi^2$. Multiplying the equation (\ref{iteration-g})$_1$ by $\frac{g^{n+1}}{\mb M_{_\#}}$ gives
		\begin{equation}\label{0orderg}
			\begin{aligned}
				&\left[\frac{(g^{n+1})^2}{2\mb M_{_\#}} \right]_t+\left[\left(\frac{\xi_1}{v^n}-\frac{u_1^n}{v^n}-\s \right)\frac{(g^{n+1})^2}{2\mb M_{_\#}} \right]_y+\frac{\nu_{\hat{\mb M}}(\xi)\left(g^{n+1} \right) ^2}{\mb M_{_\#}}  \\[3mm]
				&=\di\frac{g^{n+1}}{\mb M_{_\#}}\left[\sqrt{\hat{\mb M}}\mb K_{\hat{\mb M}}\left(\frac{g^n}{\sqrt{\hat{\mb M}}} \right) +Q_+(g^n,g^n)-Q_-(g^{n+1},g^n)-\left( \hat{\mathbf{M}}_t+\left(\frac{\xi_1}{v^n}-\frac{u_1^n}{v^n}-\s \right)\hat{\mathbf{M}}_y\right) \right]\\[3mm]
				&\quad \di-\left(\frac{\xi_1}{v^n}-\frac{u_1^n}{v^n}-\s \right)_y\frac{(g^{n+1})^2}{2\mb M_{_\#}}\\[3mm]
				& =:\sum_{i=1}^5 A
				_1^i.
			\end{aligned}
		\end{equation}
		Since $\|g^n\|_{\mathcal{S}_2}^2\leqslant\Xi^2$, by Sobolev inequality we have
		\begin{equation*}
			\begin{array}{ll}
				\di\left\|\int\f{|g^n|^2}{\mb{M}_{_\#}}d \xi\right\|_{L^\i_y} \leqslant
				C\left(\iint\f{|g^n|^2}{\mb{M}_{_\#}}d \xi
				dy\right)^{\f{1}{2}}\cdot\left(\iint\f{|g^n_y|^2}{\mb{M}_{_\#}}d \xi dy\right)^{\f{1}{2}}\leqslant
				C\Xi^2,
			\end{array}
		\end{equation*}
		and
		\begin{equation*}
			\begin{array}{ll}
				\di\left\|\int\f{|g^n_y|^2}{\mb{M}_{_\#}}d \xi\right\|_{L^\i_y} \leqslant
				C\left(\iint\f{|g^n_y|^2}{\mb{M}_{_\#}}d \xi
				dy\right)^{\f{1}{2}}\cdot\left(\iint\f{|g^n_{yy}|^2}{\mb{M}_{_\#}}d \xi dy\right)^{\f{1}{2}}\leqslant
				C\Xi^2.
			\end{array}
		\end{equation*}
		By using Lemma \ref{Lemma 4.1} and \ref{Property of K} and the smallness of $\Xi$, it holds that
		\begin{equation*}
			\begin{array}{ll}
				\di \quad \int_0^T\iint (A^1_1+A_1^2+A_1^3)d\xi dydt \\[3mm]
				\di \leqslant  \f1{16}\int_0^T\iint\frac{\nu_{\hat{\mb M}}(\xi)|g^{n+1}| ^2}{\mb M_{_\#}} d\xi dydt+C\int_0^T\iint\frac{|g^{n}| ^2}{\mb M_{_\#}} d\xi dydt\\[4mm]
				\di\qquad
				+C\int_0^T\int\left(\int\frac{|g^{n} | ^2}{\mb M_{_\#}}d\xi\int\frac{\nu_{\hat{\mb M}}(\xi)|g^{n}| ^2}{\mb M_{_\#}}d\xi \right) dydt\\[4mm]
				\di\qquad
				+C\int_0^T\int\left(\int\frac{|g^{n} | ^2}{\mb M_{_\#}}d\xi\int\frac{\nu_{\hat{\mb M}}(\xi)|g^{n+1} | ^2}{\mb M_{_\#}}d\xi \right) dydt\\[4mm]
				\di
				\leqslant \f1{8}\int_0^T\iint\frac{\nu_{\hat{\mb M}}(\xi)|g^{n+1} | ^2}{\mb M_{_\#}} d\xi dydt+C\int_0^T\iint\frac{|g^{n} | ^2}{\mb M_{_\#}} d\xi dydt+C\Xi^2\int_0^T\iint\frac{\nu_{\hat{\mb M}}(\xi)|g^{n}| ^2}{\mb M_{_\#}} d\xi dydt\\[3mm]
				\di
				\leqslant \f1{8}\int_0^T\iint\frac{\nu_{\hat{\mb M}}(\xi)|g^{n+1} | ^2}{\mb M_{_\#}} d\xi dydt+C\Xi^2T+C\Xi^4.
			\end{array}
		\end{equation*}
		Since
		$$\left|\f1{v^n} \right|=|\rho^n|=\left| \int\left(g^n+\hat{\mb M} \right)d\xi\right|\leqslant \hat{\rho}+C\Big(\int  \frac{|g^n|^2}{\mb M_{_\#}}d\xi\Big)^{\frac12}\leqslant\hat{\rho}+C\Xi, $$
		$$\left|\f{u_1^n}{v^n} \right|=|\rho^nu^n_1|=\left| \int\xi_1\left(g^n+\hat{\mb M} \right)d\xi\right|\leqslant |\hat{\rho}\hat{u}_1|+C\Big(\int  \frac{|g^n|^2}{\mb M_{_\#}}d\xi\Big)^{\frac12}\leqslant|\hat{\rho}\hat{u}_1|+C\Xi, $$
		we can deduce that $\f1{v^n}$ and $\f{u_1^n}{v^n}$ are uniformly bounded. Furthermore,
		$$|v^n_y|\leqslant C|\rho^n_y|=C\left| \int\left(g^n_y+\hat{\mb M}_y \right)d\xi\right|\leqslant C\left( \int  \frac{|g^n_y|^2}{\mb M_{_\#}}d\xi\right) ^{\frac12}+C\left( \int  \frac{|\hat{\mb M}_y|^2}{\mb M_{_\#}}d\xi\right) ^{\frac12}\leqslant C(\Xi+\d_2),$$
		and similar estimate holds for $u^n_{1y}$. Thus,
		$$\int_0^T\iint A_1^4d\xi dydt\leqslant\f1{16}\int_0^T\iint\frac{\nu_{\hat{\mb M}}(\xi)|g^{n+1} | ^2}{\mb M_{_\#}} d\xi dydt+C\delta_2^2 T,$$
		$$\int_0^T\iint A_1^5d\xi dydt\leqslant C(\Xi+\d_2)\int_0^T\iint\frac{\nu_{\hat{\mb M}}(\xi)|g^{n+1} | ^2}{\mb M_{_\#}} d\xi dydt.$$
		Integrating \eqref{0orderg} with respect to $\xi$, $y$, and $t$ and using the above estimates, we have
		\begin{equation}\label{estimategn+1}
			\begin{array}{ll}
				\di \iint\frac{|g^{n+1}|^2}{\mathbf{M}_{_\#}}d\xi dy+\int_0^T\iint\frac{\nu_{\hat{\mb M}}(\xi)|g^{n+1}|^2}{\mathbf{M}_{_\#}}d\xi dydt  \leqslant  \iint\frac{|g_0|^2}{\mathbf{M}_{_\#}}d\xi dy+C\Xi^4+C(\Xi^2+\delta_2^2) T
			\end{array}
		\end{equation}
		provided that $\Xi$ and $\d_2$ are chosen suitably small.
		
		To estimate $g_y^{n+1}$ and $g_{yy}^{n+1}$, noticing that
		\begin{equation*}
			\begin{array}{ll}
				&\di \quad \left[ \sqrt{\hat{\mb M}}\mb K_{\hat{\mb M}}\left(\frac{h}{\sqrt{\hat{\mb M}}} \right)\right]_y\di=\left( \mb L_{\hat{\mb M}}h+\nu_{\hat{\mb M}}h \right)  _y \\[5mm]
				&=\di \mb L_{\hat{\mb M}}h_y+Q(\hat{\mb M}_y,h)+Q(h,\hat{\mb M}_y)+\nu_{\hat{\mb M}}h_y+(\nu_{\hat{\mb M}})_yh\\[3mm]
				&=\di \sqrt{\hat{\mb M}}\mb K_{\hat{\mb M}}\left(\frac{h_y}{\sqrt{\hat{\mb M}}} \right)+Q(\hat{\mb M}_y,h)+Q(h,\hat{\mb M}_y)+h\int_{{\mathbb R}^3}\!\!\int_{{\mathbb S}^2_+} |(\xi-\xi_*)\cdot \Omega|\hat{\mb M}_y(\xi_*)d \xi_* d\Omega\\[5mm]
				&=\di \sqrt{\hat{\mb M}}\mb K_{\hat{\mb M}}\left(\frac{h_y}{\sqrt{\hat{\mb M}}} \right)+Q(\hat{\mb M}_y,h)+Q_+(h,\hat{\mb M}_y),
			\end{array}
		\end{equation*}
		then differentiating $\eqref{iteration-g}_1$ with respect to $y$ and utilizing the similar techniques as in deriving \eqref{estimategn+1}, for suitably small positive constants $\Xi$ and $\d_2$, we have 
		\begin{equation}\label{estimategn+1y}
			\begin{array}{ll}
				\di \quad \iint\frac{|g^{n+1}_y|^2}{\mathbf{M}_{_\#}}d\xi dy+\int_0^T\iint\frac{\nu_{\hat{\mb M}}(\xi)|g^{n+1}_y|^2}{\mathbf{M}_{_\#}}d\xi dydt \\[3mm]
				\di \leqslant  \iint\frac{|g_{0y}|^2}{\mathbf{M}_{_\#}}d\xi dy+C\Xi^4+C(\Xi^2+\delta_2^2) T+C\d_2^2\Xi^2\\[3mm]
				~~\di+C(\Xi+\delta_2)\int_0^T\iint\frac{\nu_{\hat{\mb M}}(\xi)|g^{n+1}|^2}{\mathbf{M}_{_\#}}d\xi dydt
			\end{array}
		\end{equation}
		and
		\begin{equation}\label{estimategn+1yy}
			\begin{array}{ll}
				\di \quad \iint\frac{|g^{n+1}_{yy}|^2}{\mathbf{M}_{_\#}}d\xi dy+\int_0^T\iint\frac{\nu_{\hat{\mb M}}(\xi)|g^{n+1}_{yy}|^2}{\mathbf{M}_{_\#}}d\xi dydt \\[3mm]
				\di \leqslant  \iint\frac{|g_{0yy}|^2}{\mathbf{M}_{_\#}}d\xi dy+C\Xi^4+C(\Xi^2+\delta_2^2) T+C\d_2^2\Xi^2\\[3mm]
				~~\di+C(\Xi+\delta_2)\int_0^T\iint\frac{\nu_{\hat{\mb M}}(\xi)(|g^{n+1}|^2+|g^{n+1}_y|^2)}{\mathbf{M}_{_\#}}d\xi dydt.
			\end{array}
		\end{equation}
		
		Combining \eqref{estimategn+1}, \eqref{estimategn+1y}, and \eqref{estimategn+1yy}, and choosing $\Xi,\delta_2$ and $T$ suitably small, we can get
		$$\|g^{n+1}\|_{\mathcal{S}_2}^2\leqslant \sum_{0\leqslant|\alpha|\leqslant2}\iint\frac{|\partial^\a_y g_0|^2}{\mathbf{M}_{_\#}}d\xi dy+C\Xi^4+C(\Xi^2+\delta^2_2) T+C\d_2^2\Xi^2\leqslant\Xi^2.$$
		
		Since
		\begin{equation*}
			\begin{array}{ll}
				\di \quad \rho^{n+1}(t,y)=\int f^{n+1}(t,y,\xi)d\xi&=\di\int \left( g^{n+1}(t,y,\xi)+\hat{\mb M}(t,y,\xi)\right) d\xi\\[3mm]
				&\di \geqslant \hat \rho(t,y)-C\Big(\int\frac{|g^{n+1}(t,y,\xi)|^2}{\mathbf{M}_{_\#}(\xi)}d\xi\Big)^{\frac 12}\\[3mm]
				&\di \geqslant \hat \rho(t,y)-C\Xi>0,
			\end{array}
		\end{equation*}
		we have $ \rho^{n+1}(t,y)>0$ for all $(t, y)\in \mathbb{R}_+\times\mathbb{R}$.
		
		Next we prove that $\{g^n(t,y,\xi)\}$ is a Cauchy sequence in $\mathcal{S}_0([0,T])$. Set
		$$h^n(t,y,\xi):=g^{n+1}(t,y,\xi)-g^n(t,y,\xi),~~~n=0, 1, 2, \cdots,$$
		then $h^n(t,y,\xi)$ solves the following equation
		\begin{equation}\label{iteration-h}
			\begin{cases}\di h_t^{n}+\left[  \left(\frac{\xi_1}{v^n}-\frac{u_1^n}{v^n} \right)-\left(\frac{\xi_1}{v^{n-1}}-\frac{u_1^{n-1}}{v^{n-1}} \right)\right]g^{n+1}_y+\left(\frac{\xi_1}{v^{n-1}}-\frac{u_1^{n-1}}{v^{n-1}}-\s \right) h^{n}_y\\[3mm]
				~~~~~~~\di+h^{n}\int_{\mathbb R^3}\int_{\mathbb S^2_+}|(\xi-\xi_*)\cdot \Omega|\left( g^n(t,y,\xi_*)+\hat{\mb M}(t,y,\xi_*)\right) d\Omega d \xi_*\\[3mm]
				~~~~~~~\di+g^n\int_{\mathbb R^3}\int_{\mathbb S^2_+}|(\xi-\xi_*)\cdot \Omega|h^{n-1} d\Omega d \xi_* \\[3mm]
				~~~\di=\sqrt{\hat{\mb M}}\mb K_{\hat{\mb M}}\left(\frac{h^{n-1}}{\sqrt{\hat{\mb M}}} \right) +Q_+(g^{n-1},h^{n-1})+Q_+(h^{n-1},g^{n})\\[4mm]
				~~~~~~~\di-\left[  \left(\frac{\xi_1}{v^n}-\frac{u_1^n}{v^n} \right)-\left(\frac{\xi_1}{v^{n-1}}-\frac{u_1^{n-1}}{v^{n-1}} \right)\right] \hat{\mathbf{M}}_y ,\\[4mm]
				h^{n}(0,y,\xi)=0.\end{cases}
		\end{equation}
		Multiplying the equation $\eqref{iteration-h}$ by $\frac{h^{n}}{\mb M_{_\#}}$ gives
		\begin{equation*}
			\begin{aligned}
				&\left[\frac{(h^{n})^2}{2\mb M_{_\#}} \right]_t+\left[\left(\frac{\xi_1}{v^n}-\frac{u_1^n}{v^n}-\s \right)\frac{(h^{n})^2}{2\mb M_{_\#}} \right]_y+\frac{\nu_{\hat{\mb M}}(\xi)\left(h^{n} \right) ^2}{\mb M_{_\#}}  \\
				=&\di\frac{h^{n}}{\mb M_{_\#}}\left[\sqrt{\hat{\mb M}}\mb K_{\hat{\mb M}}\left(\frac{h^{n-1}}{\sqrt{\hat{\mb M}}} \right) +Q_+(g^{n-1},h^{n-1})+Q_+(h^{n-1},g^{n})-Q_-(h^{n},g^{n})-Q_-(g^{n},h^{n-1}) \right.\\
				&~\di-  \left.\left(\frac{1}{v^n}-\frac{1}{v^{n-1}} \right)\xi_1g^{n+1}_y+  \left(\frac{u^n_1}{v^n}-\frac{u_1^{n-1}}{v^{n-1}} \right)g^{n+1}_y   -  \left(\frac{1}{v^n}-\frac{1}{v^{n-1}} \right)\xi_1\hat{\mb M}_y+  \left(\frac{u^n_1}{v^n}-\frac{u_1^{n-1}}{v^{n-1}} \right)\hat{\mb M}_y\right]\\[2mm]
				&~-\left(\frac{\xi_1}{v^n}-\frac{u_1^n}{v^n}-\s \right)_y\frac{(h^{n})^2}{2\mb M_{_\#}}\\
				=&\!:\sum_{i=1}^{10}A_2^i.
			\end{aligned}
		\end{equation*}
		By Lemma \ref{Property of K} and Cauchy inequality, we have
		\begin{equation*}
			\begin{array}{ll}
				&\di \quad \int_0^T\iint A_2^1d\xi dydt\\[5mm]
				&\di\leqslant\f1{16}\int_0^T\iint\frac{\nu_{\hat{\mb M}}(\xi)|h^{n} |  ^2}{\mb M_{_\#}} d\xi dydt+C\int_0^T\iint\frac{|h^{n-1} | ^2}{\mb M_{_\#}} d\xi dydt\\[5mm]
				&\di\leqslant  \f1{16}\int_0^T\iint\frac{\nu_{\hat{\mb M}}(\xi)|h^{n} | ^2}{\mb M_{_\#}} d\xi dydt+CT\sup_{t\in[0,T]}\iint\frac{| h^{n-1}|^2}{\mathbf{M}_{_\#}}d\xi dy.
			\end{array}
		\end{equation*}
		The terms $A_2^i(i=2, \cdots, 5)$ can be estimated directly by Lemma \ref{Lemma 4.1}. For $A_2^6$, we have
		\begin{equation*}
			\begin{array}{ll}
				&\di \quad \int_0^T\iint A_2^6d\xi dydt\\[3mm]
				&\di\leqslant\f1{16}\int_0^T\iint\frac{\nu_{\hat{\mb M}}(\xi)|h^{n}| ^2}{\mb M_{_\#}} d\xi dydt+C\int_0^T\int (\rho^n-\rho^{n-1})^2 \int\frac{\nu_{\hat{\mb M}}(\xi)|g^{n+1} | ^2}{\mb M_{_\#}} d\xi dydt\\[3mm]
				&\di\leqslant\f1{16}\int_0^T\iint\frac{\nu_{\hat{\mb M}}(\xi)|h^{n}| ^2}{\mb M_{_\#}} d\xi dydt+C\int_0^T\int \left( \int h^{n-1}d\xi\right) ^2 \int\frac{\nu_{\hat{\mb M}}(\xi)|g^{n+1} | ^2}{\mb M_{_\#}} d\xi dydt\\[3mm]
				&\di\leqslant\f1{16}\int_0^T\iint\frac{\nu_{\hat{\mb M}}(\xi)|h^{n} | ^2}{\mb M_{_\#}} d\xi dydt+C\int_0^T\iint \frac{|h^{n-1} | ^2}{\mb M_{_\#}}d\xi \int\frac{\nu_{\hat{\mb M}}(\xi)|g^{n+1}| ^2}{\mb M_{_\#}} d\xi dydt\\[3mm]
				&\di\leqslant  \f1{16}\int_0^T\iint\frac{\nu_{\hat{\mb M}}(\xi)|h^{n} | ^2}{\mb M_{_\#}} d\xi dydt+C\Xi^2\sup_{t\in[0,T]}\iint\frac{| h^{n-1}|^2}{\mathbf{M}_{_\#}}d\xi dy.
			\end{array}
		\end{equation*}
		Similarly, we can estimate $A_2^i(i=7,\cdots,10)$. Consequently, it holds that
		\begin{equation}\label{estimatehn}
			\begin{array}{ll}
				&\di \quad \sup_{t\in[0,T]}\iint\frac{|h^{n}|^2}{\mathbf{M}_{_\#}}d\xi dy+\int_0^T\iint\frac{\nu_{\hat{\mb M}}(\xi)|h^{n}|^2}{\mathbf{M}_{_\#}}d\xi dydt\\[3mm]
				&\di\leqslant C(\Xi^2+T)
				\left[\sup_{t\in[0,T]}\iint\frac{| h^{n-1}|^2}{\mathbf{M}_{_\#}}d\xi dy+\int_0^T\iint\frac{\nu_{\hat{\mb M}}(\xi)|h^{n-1}|^2}{\mathbf{M}_{_\#}}d\xi dydt \right].
			\end{array}
		\end{equation}
		provided that $\Xi$, $\d_2$, $T$ are chosen suitably small. \eqref{estimatehn} implies that $\{g^n\}$ is a Cauchy sequence in $\mathcal{S}_0([0,T])$ and the limit function $g\in \mathcal{S}_2([0,T])$.  Therefore, $f(t,y,\xi):=g(t,y,\xi)+\hat{\mb M}(t,y,\xi)\geqslant0$ is a solution to Cauchy problem \eqref{tran-Lag-B}.
		
		We prove the uniqueness of the solution  to Cauchy problem \eqref{tran-Lag-B} by contradiction. Assume that $g,\tilde{g}\in\mathcal{S}_2([0,T])$ are two solutions to the equivalent equation
		\begin{equation*}
			\begin{cases}\di g_t+\left(\frac{\xi_1}{v}-\frac{u_1}{v}-\s \right) g_y=\mb L_{\hat{\mb M}}g+Q(g,g)-\left( \hat{\mathbf{M}}_t+\left(\frac{\xi_1}{v^n}-\frac{u_1}{v}-\s \right)\hat{\mathbf{M}}_y\right),\\[3mm]
				g(0,y,\xi)=g_0(y,\xi).\end{cases}
		\end{equation*}
		Using the similar arguments as in obtaining \eqref{estimatehn}, we can deduce
		\begin{equation*}
			\begin{array}{ll}
				&\di \quad \sup_{t\in[0,T]}\iint\frac{|g-\tilde{g}|^2}{\mathbf{M}_{_\#}}d\xi dy+\int_0^T\iint\frac{\nu_{\hat{\mb M}}(\xi)|g-\tilde{g}|^2}{\mathbf{M}_{_\#}}d\xi dydt\\[3mm]
				&\di\leqslant C(\Xi^2+T)
				\left[\sup_{t\in[0,T]}\iint\frac{|g-\tilde{g}|^2}{\mathbf{M}_{_\#}}d\xi dy+\int_0^T\iint\frac{\nu_{\hat{\mb M}}(\xi)|g-\tilde{g}|^2}{\mathbf{M}_{_\#}}d\xi dydt \right],
			\end{array}
		\end{equation*}
		which implies that $g-\tilde{g}\equiv0$ provided that $\Xi$, $\d_2$, $T$ are chosen to be suitably small. Hence the proof of Proposition \ref{localexistence} is completed.

	\end{appendix}

\end{document}